\documentclass[11pt]{amsart}
\usepackage[T1]{fontenc}
\usepackage[utf8]{inputenc}
\usepackage{lmodern,amsmath,amssymb,amsthm,mathtools,mathrsfs,mathptmx}
\usepackage[margin=1.05in]{geometry}
\usepackage[colorlinks=true,linkcolor=blue,citecolor=blue,urlcolor=blue]{hyperref}

\numberwithin{equation}{section}

\theoremstyle{plain}
\newtheorem{theorem}{Theorem}[section]
\newtheorem{lemma}[theorem]{Lemma}
\newtheorem{proposition}[theorem]{Proposition}
\newtheorem{corollary}[theorem]{Corollary}
\newtheorem{introtheorem}{Theorem}

\theoremstyle{definition}

\newtheorem{remark}[theorem]{Remark}
\newtheorem{definition}[theorem]{Definition}

\newcommand{\A}{\mathbb A}
\newcommand{\ZF}{\mathbb Z\mathrm F}
\newcommand{\Fr}{\mathrm{Fr}}
\newcommand{\Sm}{\mathrm{Sm}}
\newcommand{\Nis}{\mathrm{Nis}}
\newcommand{\Spec}{\operatorname{Spec}}
\newcommand{\pr}{\operatorname{pr}}
\newcommand{\id}{\mathrm{id}}
\newcommand{\perf}{\mathrm{perf}}
\newcommand{\tr}{\mathrm{tr}}
\newcommand{\GL}{\mathrm{GL}}
\newcommand{\SL}{\mathrm{SL}}
\newcommand{\diag}{\operatorname{diag}}
\newcommand{\sd}{\operatorname{sd}}

\renewcommand{\leq}{\leqslant}
\renewcommand{\geq}{\geqslant}

\title{Frobenius and perfection in framed motives}
\author{Grigory Garkusha}
\address{Department of Mathematics, Swansea University, Fabian Way, Swansea SA1 8EN, UK}
\email{g.garkusha@swansea.ac.uk}
\urladdr{\url{https://ggarkusha.github.io}}
\date{}



\thanks{Explore interactive visualisations on framed correspondences and motives at \url{https://ggarkusha.github.io/math-atlas/}}

\begin{document}

\maketitle

\footskip30pt

\thispagestyle{empty} \pagestyle{plain}

\begin{abstract}
We develop explicit Frobenius and perfection arguments for framed
motives over imperfect fields. In positive characteristic,
power framing homotopies and local framed transfers make relative
Frobenius invertible after inverting the characteristic. Descent
through successive Frobenius twists then gives a proof of strict
homotopy invariance for homotopy invariant quasi-stable additive framed presheaves
and identifies the canonical comparison with the perfect closure
as a stable Nisnevich local equivalence of framed motives.
We extend the comparison to $E$-framed motives for every motivic Thom spectrum $E$.
It also yields the cancellation theorem and the cone theorem.
For commutative symmetric Thom ring spectra,
we identify triangulated categories of effective $E$-framed motives and effective
$E$-modules and prove perfection invariance of these triangulated
categories. All comparisons are made after inverting the
exponential characteristic. 
\end{abstract}

\section*{Introduction}

Framed motives introduced and studied in \cite{GPfm} compute the stable motivic homotopy type of
smooth schemes by Voevodsky's framed correspondences~\cite{Voe}. The principal results
of the theory are proved over (infinite) perfect fields:
\begin{itemize}
\item strict $\A^1$-invariance of framed presheaves \cite{GP};
\item the cancellation and cone theorems \cite{AGP,GNP};
\item their extension to motivic Thom spectra \cite{GN};
\item the description of $\mathsf{SH}(k)$ and of $E$-modules by spectral
categories of framed correspondences \cite{Corr}.
\end{itemize}
This paper removes the perfectness assumption after inverting the
exponential characteristic. The method is to compare framed
constructions over a field with the corresponding constructions over its
perfect closure by explicit Frobenius arguments that keep track of the
canonical maps defined by base change of framed correspondences. The
main comparison takes place for presheaves of $S^1$-spectra before
$T$-stabilization.

Throughout, $k$ is a field and $p$ is its exponential characteristic. We put $R=\mathbb Z[1/p]$ and
$K=k^{\perf}$, and write $\tau\colon\Spec K\to\Spec k$ for the canonical
morphism. We write $\Sm/k$ for smooth separated schemes of finite type
over $k$, and $\ZF_*(k)$ for the graded additive category of linear
framed correspondences. For spectra, $[1/p]$ denotes localization
obtained by inverting multiplication by $p$.

\subsection*{Frobenius transfers and descent}
Let $p>1$ and let $U$ admit an \'{e}tale map to $\A_k^d$. We construct a
framed transfer $\tr_U$ of level $d$ for the relative Frobenius
$F_{U/k}\colon U\to U^{(p)}$. Its two composites with relative Frobenius
satisfy
\[
 F_{U/k}\circ\tr_U\simeq\Theta_{p,d,U^{(p)}},
 \qquad
 \tr_U\circ F_{U/k}\simeq\Theta_{p,d,U}.
\]
Here $\Theta_{p,d}$ is the correspondence with power framing
$(t_1^p,\ldots,t_d^p)$, and $\simeq$ denotes equality modulo additive
$\A^1$-homotopies at the indicated framing level. Explicit polynomial
homotopies make the power framing operators invertible on homotopy
invariant quasi-stable additive framed presheaves with $R$-coefficients.
Further homotopies compare taking powers of framing functions with
taking powers of ambient coordinates. Together, these calculations show
that pullback along relative Frobenius is a natural isomorphism on the
associated framed Nisnevich sheaves (Theorem~\ref{prop:Fiso}).

The second ingredient is descent to the perfect closure. Write
 $K=\bigcup_{n\geq0}k^{1/p^n}$. Following the transfer along field isomorphisms, the successive
inclusions in this tower become Frobenius base change, so the Frobenius
comparisons give compatible identifications at all finite stages $n$.
All correspondence data and their defining relations descend to some
finite stage. The extension $k^{1/p^n}/k$ need not be finite, and no
assumption on $[k:k^p]$ is made.

\subsection*{Strict homotopy invariance}
The first application is the following theorem.

\begin{introtheorem}[Strict homotopy invariance]\label{thm:main}
Let $k$ be an imperfect field and let $\mathcal F$ be an
$\A^1$-invariant quasi-stable additive $\ZF_*(k)$-presheaf of abelian
groups. Then $\mathcal A=(\mathcal F\otimes_{\mathbb Z}R)_{\Nis}$ is
strictly $\A^1$-invariant: for every $X\in\Sm/k$ and $i\geq0$,
\[
 H^i_{\Nis}(X,\mathcal A)\xrightarrow{\simeq}
 H^i_{\Nis}(X\times\A^1,\mathcal A).
\]
\end{introtheorem}

The proof carries the sheaf over to $K$, applies strict homotopy
invariance over $K$ \cite[Theorem 1.1]{GP}, and compares Nisnevich
cohomology by invariance of the small Nisnevich site under universal
homeomorphisms. The comparison between cohomology on the big and small
sites is made functorial in the smooth scheme. In characteristic $2$ the
formulation of \cite{GP} with $\mathbb Z[1/2]$-coefficients suffices.

We follow the method of Suslin's work on motivic complexes over
imperfect fields \cite{Sus17}. Suslin proves strict $\mathbb A^1$-invariance
by passage to the perfect closure and comparison of Nisnevich cohomology
\cite[Theorems 2.13 and 3.4]{Sus17}. For framed correspondences the
required transfer is constructed locally in \'{e}tale coordinates.
Its composites with Frobenius are the power framing operators above.
Druzhinin \cite{Dru} states strict homotopy invariance with integral
coefficients over any field. Our proof of Theorem~\ref{thm:main}
uses $R$-coefficients and does not use \cite{Dru}. It supplies the
explicit Frobenius comparisons needed for the later results.

\subsection*{Perfection invariance of framed motives}
For a field $F$ and $Y\in\Sm/F$, let $M_{\mathrm{fr},F}(Y)$ denote the
framed motive of $Y$, a presheaf of $S^1$-spectra on $\Sm/F$. For a
presheaf $\mathcal E$ on $\Sm/K$ put $(\tau_*\mathcal E)(V)=\mathcal
E(V_K)$, $V\in\Sm/k$. The following theorem is a framed counterpart of
Suslin's base change formula \cite[Theorem 4.12(3)]{Sus17} and of its
consequence, invariance under purely inseparable extensions
\cite[Corollary 4.13]{Sus17}.

\begin{introtheorem}[Theorem~\ref{thm:local-perfection}]\label{thm:intro-perfection}
Let $k$ be any field of exponential characteristic $p$, and put
$K=k^{\perf}$. For every $X\in\Sm/k$, the canonical morphism induced by
base change of framed correspondences,
\[
 b_X[1/p]\colon M_{\mathrm{fr},k}(X)[1/p]\longrightarrow
 \tau_*M_{\mathrm{fr},K}(X_K)[1/p],
\]
is a stable Nisnevich local equivalence of presheaves of $S^1$-spectra
on $\Sm/k$.
\end{introtheorem}

Thus $b_X[1/p]$ identifies the Nisnevich sheaves of stable homotopy
groups. At the base field it is already a stable equivalence
$M_{\mathrm{fr},k}(X)(\Spec k)[1/p]\simeq M_{\mathrm{fr},K}(X_K)(\Spec
K)[1/p]$. After stable Nisnevich fibrant replacement, formed on $\Sm/k$
and on $\Sm/K$ respectively, there are natural zigzags of stable
equivalences
\[
 M_{\mathrm{fr},k}(X)[1/p]_f(V)\simeq M_{\mathrm{fr},K}(X_K)[1/p]_f(V_K),
 \qquad V\in\Sm/k,
\]
(see Corollary~\ref{cor:local-perfection-fibrant}). The proof works with the
linear complexes underlying framed motives. The local Frobenius
transfers act on cycles and boundaries at chosen finite framing levels,
the power framing identities give isomorphisms on stabilized homology,
and descent in the target extends the computation from schemes with
\'{e}tale coordinates to all smooth schemes. Working at finite framing
levels keeps track of the order of framing coordinates when transfers
are composed with suspensions.

The same method applies to motivic Thom spaces and to Thom spectra. For a
closed subset $Z\subset X$, not necessarily smooth, put $Q=X/(X-Z)$.
Theorem~\ref{thm:pairs} proves the analogue of
Theorem~\ref{thm:intro-perfection} for $M_{\mathrm{fr},k}(Q)$. For every
motivic Thom $T$-spectrum $E$ in the sense of \cite{GN},
Theorem~\ref{thm:E-perfection} shows that
\[
 M_{E,k}(Q)[1/p]\longrightarrow\tau_*M_{E_K,K}(Q_K)[1/p],
 \qquad E_K=\tau^*E,
\]
is a stable Nisnevich local equivalence, with the same consequences
after fibrant replacement. This holds for every motivic Thom spectrum,
without bounding, symmetry or multiplicativity assumptions. Its proof
uses the Thom layer descriptions of \cite{GN}.

\subsection*{Cancellation and the cone theorem}
The perfection comparisons carry over the main theorems of
\cite{AGP,GNP,GN} from the perfect closure. Over every infinite field,
after inverting $p$, we obtain:
\begin{itemize}
\item the cancellation theorem for framed motives
(Theorem~\ref{cor:cancellation});
\item motivic fibrancy of the associated $(S^1,\mathbb G)$-bispectra
(Corollary~\ref{bispectrum});
\item the stable cone theorem for framed motives of relative motivic
spheres (Theorem~\ref{thm:gnp});
\item cancellation for $E$-framed motives
(Theorem~\ref{thm:E-cancellation}).
\end{itemize}
Imperfect fields are infinite, so over imperfect fields these results
hold without restriction. Over infinite perfect fields they are the
results of \cite{AGP,GNP,GN} with $p$ inverted.

\subsection*{Effective $E$-framed motives}
The last part of the paper concerns the triangulated categories defined
by the spectral categories of \cite[\S\S5--6]{Corr}. Here $E$ is a
commutative symmetric Thom ring $T$-spectrum over $k$ with bounding
constant $d\leq1$ and contractible alternating group action in the sense of \cite{GN}, as required
for the symmetric monoidal module constructions of \cite{Corr}.
\begin{itemize}
\item $\mathscr D_{E,R}(k)$ denotes the Nisnevich localization of the
category of right modules over the spectral category of $E$-framed
correspondences with $R$-coefficients.
\item $D_E^{\mathrm{eff}}(k;R)$ is the full subcategory of
$\mathscr D_{E,R}(k)$ of objects whose Nisnevich sheaves of stable
homotopy groups are $\A^1$-invariant.
\item $\operatorname{Mod}_{\mathsf{SH}(k)[1/p]}E[1/p]$ is the homotopy
category of the model category of $E$-modules in symmetric $T$-spectra
with multiplication by $p$ inverted.
\item $\operatorname{Mod}^{\mathrm{eff}}_{\mathsf{SH}(k)[1/p]}E[1/p]$ is
its localizing subcategory generated by
$E[1/p]\wedge\Sigma_T^\infty X_+$, $X\in\Sm/k$.
\end{itemize}
The following theorem extends \cite[Theorem 6.12 and Corollary
6.13]{Corr} to arbitrary fields after inverting $p$.

\begin{introtheorem}[Theorem~\ref{thm:E-effective}]\label{thm:intro-effective}
Let $k$ be any field of exponential characteristic $p$, let $E$ be as
above, and put $E_K=\tau^*E$.
\begin{enumerate}
\item Suslin localization $C_*$ is left adjoint to the inclusion
$D_E^{\mathrm{eff}}(k;R)\hookrightarrow\mathscr D_{E,R}(k)$. It
identifies the motivic localization of $\mathscr D_{E,R}(k)$ with
$D_E^{\mathrm{eff}}(k;R)$, and
$D_E^{\mathrm{eff}}(k;R)\simeq\mathsf{SH}^{\Nis}_{S^1}(\mathcal
O^E_{\Delta,R,k})$.
\item Canonical base change is a triangulated equivalence
\[
 \mathbf L\tau^*\colon D_E^{\mathrm{eff}}(k;R)
 \xrightarrow{\simeq}D_{E_K}^{\mathrm{eff}}(K;R).
\]
\item Motivic realization induces a triangulated equivalence
\[
 D_E^{\mathrm{eff}}(k;R)
 \simeq\operatorname{Mod}^{\mathrm{eff}}_{\mathsf{SH}(k)[1/p]}E[1/p].
\]
\end{enumerate}
These categories are compactly generated with no boundedness
restriction on their objects. The compact generators on the framed side
are $m_{E,k}(X)=C_*\mathcal O^E_{R,k}(-,X)$, $X\in\Sm/k$. Base change
sends $m_{E,k}(X)$ to $m_{E_K,K}(X_K)$, and realization sends it to
$E[1/p]\wedge\Sigma_T^\infty X_+$.
\end{introtheorem}

Full faithfulness in (2) follows from the comparison of spectral
representables, a consequence of Theorem~\ref{thm:E-perfection}.
Essential surjectivity combines descent of smooth schemes to a finite
stage of $K/k$ with relative Frobenius, keeping the coefficient spectrum
$E_K$ fixed. For perfect $k$, assertion (3) is \cite[Corollary
6.13]{Corr}. For imperfect $k$, the equivalence is obtained by transferring
from $K$ along (2) and the perfection invariance of $E$-modules below.

Two further inputs enter the proof of (3). First, canonical pullback is
a symmetric monoidal equivalence
\[
 \tau^*\colon\mathsf{SH}(k)[1/p]\xrightarrow{\simeq}\mathsf{SH}(K)[1/p]
\]
(see Proposition~\ref{prop:SH-perfection}). This is a special case of the
theorem of Elmanto and Khan that $\mathsf{SH}(-)[1/p]$ is invariant under
universal homeomorphisms of schemes of exponential characteristic $p$
\cite[Corollary 2.1.5]{EK}. We give an independent proof based on the
Frobenius transfers of Section~\ref{thma}, the $S^1$-stable motivic homotopy
category and continuity of stable motivic categories
\cite[Proposition C.12(4)]{Hoy}. The result of \cite{EK} is not used as
an input. Second, derived pullback is an equivalence
\[
 \tau^*\colon\operatorname{Mod}_{\mathsf{SH}(k)[1/p]}E[1/p]
 \xrightarrow{\simeq}\operatorname{Mod}_{\mathsf{SH}(K)[1/p]}E_K[1/p]
\]
which preserves effective modules
(Corollary~\ref{cor:E-module-perfection}). Its proof uses the model
categories of $E$-modules of \cite{SS00}, whose existence assumes the
monoid axiom. We deduce the monoid axiom for the Brown--Gersten stable
model structure from Jardine's results on motivic symmetric spectra
\cite{Jar00}.

Finally, we extend \cite[Corollary 6.9]{Corr} to arbitrary fields. Let
$\mathcal C_k^{\mathrm{st}}$ be the $\mathbb G_m$-stabilization of
$\mathsf{SH}^{\Nis}_{S^1}(\mathcal O^E_{\Delta,R,k})\simeq
D_E^{\mathrm{eff}}(k;R)$, formed with respect to the twist
$A(1)=A\boxtimes^{\mathbf L}\mathbb G_m^{\wedge1}$ as in
\cite[Section~6]{Corr}.

\begin{introtheorem}[Theorem~\ref{cor:E-stabilized}]
\label{thm:intro-stabilized}
Let $k$ be any field of exponential characteristic $p$, and let $E$ be
as in Theorem~\ref{thm:intro-effective}. There is a triangulated
equivalence
\[
 \overline{\mathcal R}_k\colon\mathcal C_k^{\mathrm{st}}
 \xrightarrow{\simeq}\operatorname{Mod}_{\mathsf{SH}(k)[1/p]}E[1/p].
\]
The suspension functor $\Sigma^\infty_{\mathbb G_m}$ is fully faithful,
and the cancellation isomorphism
$\operatorname{Hom}(A,B[j])\cong\operatorname{Hom}(A(1),B(1)[j])$ holds for
all $A,B\in\mathsf{SH}^{\Nis}_{S^1}(\mathcal O^E_{\Delta,R,k})$ and
$j\in\mathbb Z$. The essential image of
$\overline{\mathcal R}_k\Sigma^\infty_{\mathbb G_m}$ is
$\operatorname{Mod}^{\mathrm{eff}}_{\mathsf{SH}(k)[1/p]}E[1/p]$, and
$\overline{\mathcal R}_k$ is compatible with base change to $K$.
\end{introtheorem}

For perfect $k$ this is \cite[Corollary 6.9]{Corr} with
$R$-coefficients. An imperfect field is infinite, so cancellation over
$k$ follows from Theorem~\ref{thm:E-cancellation}, and the remaining
assertions are transferred from $K$.

\subsection*{Conventions}
A framed presheaf is contravariant, so $(b\circ a)^*=a^*b^*$. An
additive framed presheaf is a presheaf on $\ZF_*(k)$. In particular, it
respects the additivity relation for disjoint support components.
Quasi-stability means that $\sigma_X^*$ is an automorphism, and does not
require it to be the identity. For abelian presheaves, inverting $p$
means tensoring with $R$.

We compose explicit correspondences by placing the source framing first.
Thus, for $a=(Z,V,\varphi,g)\in\Fr_n(X,Y)$ and
$b=(Z',V',\psi,g')\in\Fr_m(Y,S)$, the composite $b\circ a$ has ambient
space $\A_X^{n+m}$, \'{e}tale neighbourhood $V\times_YV'$, framing
$(\varphi,\psi)$ and target map $g'\circ\pr_{V'}$. The computations in
the proof of Theorem~\ref{thm:main} are performed at fixed framing levels
with this convention. The coordinate permutations used later to
establish quasi-stability are treated by explicit elementary matrix homotopies.

\subsection*{Organization}
Section~\ref{thma} proves Theorem~\ref{thm:main}. It contains
framing computations, the power homotopies, the framed naturality of
relative Frobenius, the local Frobenius transfers
(Theorem~\ref{prop:Fiso}) and the descent of the sheaf to the perfect
closure. Section~\ref{sec:perfection-framed-motives} proves
Theorem~\ref{thm:intro-perfection} and deduces cancellation and motivic
fibrancy of bispectra. Section~\ref{sec:pairs} treats framed motives of
motivic Thom spaces, and Section~\ref{sec:gnp} the cone theorem.
Section~\ref{sec:E-perfection} proves perfection invariance and
cancellation for $E$-framed motives. Section~\ref{sec:E-effective}
proves perfection invariance of $\mathsf{SH}(-)[1/p]$,
Theorem~\ref{thm:intro-effective}, perfection invariance of $E$-modules
and Theorem~\ref{thm:intro-stabilized}.

\section{Proof of Theorem~A}\label{thma}

\subsection{Two elementary framing computations}

For $a\in\mathcal O(X)^\times$ write
\[
 \langle a\rangle_X=(\{0\}\times X,\A^1_X,at,\pr)
 \in\ZF_1(X,X),\qquad \langle1\rangle_X=\sigma_X,
\]
and for $l\geq1$ put
\[
 \theta_{l,X}=(\{0\}\times X,\A^1_X,t^l,\pr),
 \qquad
 \Theta_{l,n,X}=(\{0\}\times X,\A^n_X,(t_1^l,\ldots,t_n^l),\pr),
\]
with $\Theta_{l,0,X}=\id_X$. With the composition convention above,
$\Theta_{l,n,X}=\theta_{l,X}^{\circ n}$. We write $\simeq$ for
equality modulo additive $\A^1$-homotopies at the indicated level.
Every $\A^1$-invariant additive framed presheaf sends such equalities
to equal homomorphisms.

The following lemma is due to Druzhinin \cite[Lemmas~C.3-C.4]{Dru}.

\begin{lemma}\label{lem:elem}
Let \(t=(t_1,\ldots,t_n)^{\mathsf T}\) be the column vector of standard coordinates on \(\mathbb A^n_X\).
Let $M,N\in\GL_n(\mathcal O(X))$ and suppose $N$ is a product of
elementary matrices. Then the linear framings $Mt$ and $MNt$ on
$\A^n_X$, both taken with support $\{0\}\times X$ and target map
$\pr$, define $\A^1$-homotopic elements of $\Fr_n(X,X)$. In
particular, for $a,b\in\mathcal O(X)^\times$,
\begin{equation}\label{eq:units}
 \langle a\rangle_X\circ\langle b\rangle_X
 \simeq \sigma_X\circ\langle ab\rangle_X.
\end{equation}
\end{lemma}

\begin{proof}
Write $N=N_1\cdots N_r$ with $N_i=I+c_ie_{j_ik_i}$ elementary. For a
single elementary factor $E=I+c\,e_{jk}$ the family
$E_\lambda=I+(1-\lambda)c\,e_{jk}$, $\lambda\in\A^1$, consists of
elementary, hence invertible, matrices. Consequently, for any
$P\in\GL_n(\mathcal O(X))$ the framing $PE_\lambda t$ on
$\A^n_{X\times\A^1_\lambda}$ has zero locus
$\{0\}\times X\times\A^1_\lambda$ throughout, and
$(\{0\}\times X\times\A^1_\lambda,\A^n_{X\times\A^1_\lambda},
PE_\lambda t,\pr)$ is a framed homotopy from $PEt$ to $Pt$. Removing
the factors $N_r,N_{r-1},\ldots,N_1$ one at a time in this way gives a
chain of $r$ framed homotopies from $MNt$ to $Mt$.

For \eqref{eq:units}, the left-hand side is the linear framing of
$\A^2_X$ with matrix $\diag(b,a)$ and the right-hand side the one with
matrix $\diag(ab,1)$, both with support $\{0\}\times X$ and target
$\pr$. Since
\[
 \diag(b,a)^{-1}\diag(ab,1)=\diag(a,a^{-1}),
\]
which is a product of elementary matrices by the Whitehead lemma, the
first assertion applies.
\end{proof}

\begin{lemma}\label{lem:unit}
Suppose $\operatorname{char}k=p>0$. Then
\begin{equation}\label{eq:thetap}
 \theta_{p,X}\simeq p\,\langle-1\rangle_X
 \qquad\text{in }\ZF_1(X,X).
\end{equation}
Consequently, for every $\A^1$-invariant quasi-stable additive
$\ZF_*(k)$-presheaf $\mathcal M$ of $R$-modules, $\theta_{p,X}^*$ is
an automorphism of $\mathcal M(X)$, and so is $\Theta_{p,n,X}^*$ for
every $n\geq0$.
\end{lemma}

\begin{proof}
The polynomial $f_\lambda=t^p-\lambda t\in\mathcal O(X)[\lambda][t]$
is monic of degree $p$ in $t$, so its vanishing locus in
$\A^1_{X\times\A^1_\lambda}$ is finite locally free of degree $p$ over
$X\times\A^1_\lambda$. Therefore
$(Z(f_\lambda),\A^1_{X\times\A^1_\lambda},f_\lambda,\pr)\in
\mathrm{Fr}_1(X\times\mathbb A^1_\lambda,X)$ is a framed
homotopy. At $\lambda=0$ it is $\theta_{p,X}$. At $\lambda=1$ it is
$(\mathbb F_p\times X,\A^1_X,t^p-t,\pr)$, because
$t^p-t=\prod_{c\in\mathbb F_p}(t-c)$ and $\mathbb F_p\subset k$.

The support $\mathbb F_p\times X$ is a disjoint union of $p$ sections,
so additivity splits this correspondence into $p$ summands: the $c$-th
of them is
\[
 \tau_c=\bigl(\{c\}\times X,\
 \A^1_X\smallsetminus((\mathbb F_p\smallsetminus\{c\})\times X),\
 t^p-t,\ \pr\bigr),\qquad c\in\mathbb F_p .
\]
For $c\neq0$ the supports of $\tau_c$ and $\tau_0$ differ, so these are
different elements of $\Fr_1(X,X)$. We connect them by a translation
homotopy. For this consider on $\A^1_{X\times\A^1_\lambda}$ the closed subset
$D_c=\bigcup_{a\in\mathbb F_p^\times}\{t=\lambda c+a\}$ (the pieces are pairwise disjoint) and the function
\[
 f_{c,\lambda}=(t-\lambda c)^p-(t-\lambda c)
 =\prod_{a\in\mathbb F_p}(t-\lambda c-a).
\]
On the open subscheme $\A^1_{X\times\A^1_\lambda}\smallsetminus D_c$ the
zero locus of $f_{c,\lambda}$ is exactly the graph $\{t=\lambda c\}$,
which is finite over $X\times\A^1_\lambda$. Hence
\[
 \bigl(\{t=\lambda c\},\ \A^1_{X\times\A^1_\lambda}\smallsetminus D_c,\
 f_{c,\lambda},\ \pr\bigr)\in \mathrm{Fr}_1(X\times\mathbb A^1_\lambda,X)
\]
is a homotopy. At $\lambda=0$ it is $\tau_0$. At $\lambda=1$ it is $\tau_c$. Indeed, we
have $D_c=(\mathbb F_p\smallsetminus\{c\})\times X$ as 
$\{c+a:a\in\mathbb F_p^\times\}=\mathbb F_p\setminus\{c\}$, $c^p=c$, $f_{c,1}=t^p-t$. 
We see that $\tau_c\simeq\tau_0$
for every $c\in\mathbb F_p$. It remains to identify
\[
 \tau_0=\bigl(\{0\}\times X,\
 \A^1_X\smallsetminus(\mathbb F_p^\times\times X),\
 u^p-u,\ \pr\bigr),
 \qquad u^p-u=(-u)(1-u^{p-1}),
\]
where we write $u$ for the coordinate $t$.
The factor $w=1-u^{p-1}$ vanishes exactly on
$\mathbb F_p^\times\times X$, hence is a unit on this neighbourhood,
and $w|_{\{0\}\times X}=1$. 
Consider
\[
 w_\mu=1-\mu u^{p-1},\qquad
 V=\mathbb A^1_X\setminus(\mathbb F_p^\times\times X),
 \qquad
 W=D_{V\times\mathbb A^1_\mu}(w_\mu).
\]
Then $w_\mu$ is invertible on $W$. Since its restriction
to $\{0\}\times X\times\mathbb A^1_\mu$ is $1$, the open
subscheme $W$ contains this zero section. On $W$ we have an equality of ideals
$(w_\mu\cdot(-u))=(u)$, so the framing $w_\mu\cdot(-u)$
has zero scheme exactly
$\{0\}\times X\times\mathbb A^1_\mu$.
Its projection to $X\times\mathbb A^1_\mu$ is an
isomorphism, hence finite.
So
$$H=\bigl(
\{0\}\times X\times\mathbb A^1_\mu,\;
W,\;
-u(1-\mu u^{p-1}),\;
\operatorname{pr}_X
\bigr)
\in\mathrm{Fr}_1(X\times\mathbb A^1_\mu,X)$$
is a homotopy from $\langle-1\rangle_X$ to $\tau_0$. Since
$\theta_{p,X}\simeq\sum_{c\in\mathbb F_p}\tau_c$ and
$\tau_c\simeq\tau_0\simeq\langle-1\rangle_X$ for every $c$, we obtain
\eqref{eq:thetap}.

Now put $s=\sigma_X^*$ and $v=\langle-1\rangle_X^*$. Lemma
\ref{lem:elem} with $a=b=-1$ gives
$\langle-1\rangle_X\circ\langle-1\rangle_X\simeq\sigma_X\circ\sigma_X$,
whence $v^2=s^2$. Since $s$ is invertible, so is $v$, with
$v^{-1}=s^{-2}v$. As $p\in R^\times$, the endomorphism
$\theta_{p,X}^*=p\,v$ of $\mathcal M(X)$ is invertible. Finally
$\Theta_{p,n,X}=\theta_{p,X}^{\circ n}$, so
$\Theta_{p,n,X}^*=(p\,v)^n$ is invertible as well.
The framed homotopy is established integrally. The assumption 
that \(p\) is invertible is used only to deduce that \(\theta_{p,X}^{*}\) is an automorphism.
\end{proof}

\subsection{Moving a power between two coordinates}
For an integer $l\geq1$ define
\[
 D_l(x,y)=(x^l,y),\qquad
 B_u=\begin{pmatrix}1&u\\0&1\end{pmatrix}
     \begin{pmatrix}1&0\\-u&1\end{pmatrix}
     \begin{pmatrix}1&u\\0&1\end{pmatrix},\qquad
 H_u=B_u^{-1}\circ D_l\circ B_u.
\]
Explicitly $B_u=\left(\begin{smallmatrix}1-u^2&2u-u^3\\
-u&1-u^2\end{smallmatrix}\right)\in\SL_2(k[u])$, with
$B_0=I_2$ and $B_1(x,y)=(y,-x)$. $B_u$ gives rise to the ``rotational homotopy''. Consequently
\begin{equation}\label{eq:H}
 H_0(x,y)=(x^l,y),\qquad H_1(x,y)=(x,y^l).
\end{equation}

The family $H_u$ defines a morphism
\[
 \mathcal H\colon\A^2\times\A^1_u\longrightarrow
 \A^2\times\A^1_u,\qquad
 (x,y,u)\longmapsto(H_u(x,y),u),
\]
over $\A^1_u$. The morphism
$\mathcal B(x,y,u)=(B_u(x,y),u)$ is an automorphism, and
\[
 \mathcal H=\mathcal B^{-1}\circ
 (D_l\times\id_{\A^1_u})\circ\mathcal B.
\]
The morphism $D_l\times\id_{\A^1_u}$ is finite locally free
of degree $l$: its coordinate ring homomorphism
\[
 k[u,s,t]\longrightarrow k[u,x,y],\qquad
 u\longmapsto u,\quad s\longmapsto x^l,\quad t\longmapsto y,
\]
makes $k[u,x,y]$ a free $k[u,s,t]$-module with basis
$1,x,\ldots,x^{l-1}$. So $\mathcal H$ is also
finite locally free of degree $l$.

Moreover, write $H_u=(h_1,h_2)$ and
$(v,w)^{\mathsf T}=B_u(x,y)^{\mathsf T}$.
Since
\[
 (h_1,h_2)^{\mathsf T}=B_u^{-1}(v^l,w)^{\mathsf T}
\]
and $B_u$ is invertible over $k[u]$, the ideals satisfy
\[
 (h_1,h_2)=(v^l,w),\qquad (v,w)=(x,y).
\]
Consequently,
\begin{equation}\label{eq:zero}
 \sqrt{(h_1,h_2)}
 =\sqrt{(v^l,w)}
 =(v,w)
 =(x,y)
 \qquad\text{in }k[u,x,y].
\end{equation}
Thus the zero locus of the family $H_u$ has underlying
closed subset $\{(0,0)\}\times\A^1_u$.

\begin{lemma}\label{lem:power}
Let $l\geq1$, let $b=(Z,V,\varphi,g)\in\Fr_n(X,Y)$, and define
\[
 b^{[l]}=(Z,V,(\varphi_1^l,\ldots,\varphi_n^l),g),\qquad
 \mu_l:\A^n_X\longrightarrow\A^n_X,
 \quad t_i\longmapsto t_i^l.
\]
There are framed $\A^1$-homotopies at level $2n$
\begin{align}
 \sigma_Y^n\circ b^{[l]}&\simeq \Theta_{l,n,Y}\circ b,
 \label{eq:framepower}\\
 \sigma_Y^n\circ\mu_l^*b&\simeq \Theta_{l,n,Y}\circ b.
 \label{eq:ambientpower}
\end{align}
\end{lemma}

\begin{proof}
To check \eqref{eq:framepower}, introduce coordinates
$s_1,\ldots,s_n$ on $\A_s^n$ and a homotopy parameter $u$.
If $\pi\colon V\to\A_X^n$ is the original
\'{e}tale neighbourhood map, then
\[
 \pi\times\id_{\A_s^n}\times\id_{\A_u^1}\colon
V\times\A_s^n\times\A_u^1
 \longrightarrow
 \A_X^n\times\A_s^n\times\A_u^1
 =\A^{2n}_{X\times\A_u^1}
\]
is \'{e}tale. Write
$H_u(x,y)=(h_1(x,y,u),h_2(x,y,u))$ and equip this
neighbourhood with the ordered framing
\[
\begin{aligned}
 \Phi_u=\bigl(
 &h_1(\varphi_1,s_1,u),\ldots,h_1(\varphi_n,s_n,u),
 h_2(\varphi_1,s_1,u),\ldots,h_2(\varphi_n,s_n,u)
 \bigr)
\end{aligned}
\]
and target map $g\circ\pr_V$.
By \eqref{eq:zero}, applied to each pair $(\varphi_i,s_i)$,
the common zero support of these functions is
$Z\times\{0\}\times\A_u^1$, which is closed in
$\A^{2n}_{X\times\A_u^1}$ and finite over $X\times\A_u^1$.
These data therefore define a framed homotopy of level $2n$.
By \eqref{eq:H} its endpoints have framings
$(\varphi_1^l,\ldots,\varphi_n^l,s_1,\ldots,s_n)$ and
$(\varphi_1,\ldots,\varphi_n,s_1^l,\ldots,s_n^l)$, which are the two
correspondences in \eqref{eq:framepower}.

To check \eqref{eq:ambientpower}, start with the composite
$\sigma_Y^n\circ b$. If $\pi\colon V\to\A_X^n$ is the
\'{e}tale neighbourhood map of $b$, this composite is represented by
\[
 \bigl(
 Z\times\{0\},\;
 W:=V\times\A_s^n,\;
 (\varphi_1,\ldots,\varphi_n,s_1,\ldots,s_n),\;
 g\circ\pr_V
 \bigr).
\]
Its neighbourhood map is
$ \rho=\pi\times\id_{\A_s^n}\colon
 W\longrightarrow
 \A_X^n\times\A_s^n=\A_X^{2n}$
with right hand side coordinates ordered as
$(t_1,\ldots,t_n,s_1,\ldots,s_n)$.

Put $S=X\times\A_u^1$ and
$P=\A_X^{2n}\times\A_u^1=\A_S^{2n}$.
Writing $H_u(x,y)=(h_1(x,y,u),h_2(x,y,u))$, define
$\mathcal H\colon P\to P$ over $S$ by the ordered
coordinate functions
\[
 \bigl(
 h_1(t_1,s_1,u),\ldots,h_1(t_n,s_n,u),
 h_2(t_1,s_1,u),\ldots,h_2(t_n,s_n,u)
 \bigr).
\]
This is the product over $S$ of $n$ copies of the
finite locally free morphism defined by $H_u$,
hence it is finite locally free of degree $l^n$.

Extend the correspondence constantly over $\A_u^1$ and
pull back its data along $\mathcal H$. Explicitly, put
\[
 \widetilde W=
 P\times_{\mathcal H,P,\rho\times\id_{\A_u^1}}
 (W\times\A_u^1),
 \qquad
 \widetilde Z=
 \mathcal H^{-1}(Z\times\{0\}\times\A_u^1),
\]
where $\widetilde Z$ is understood as a closed subset.
Let $r\colon\widetilde W\to P$ and
$q\colon\widetilde W\to W\times\A_u^1$ be the projections.
The map $r$ is \'{e}tale by base change. Equip
$\widetilde W$ with the framing
\[
 (\varphi_1\circ q,\ldots,\varphi_n\circ q,
   s_1\circ q,\ldots,s_n\circ q)
\]
and target map $g\circ\pr_V\circ q$, with all functions
and projections extended constantly over $\A_u^1$
(on \(W\times\mathbb A^1_u\), \(\varphi_i\) means \(\varphi_i\circ\pr_V\circ\pr_W\), while \(s_i\) 
means \(s_i\circ\pr_W\); with this convention, the notation \(\varphi_i\circ q\) and \(s_i\circ q\) is well-defined).

The support $\widetilde Z$ is finite over $S$. Indeed, the
morphism
$\widetilde Z\longrightarrow Z\times\{0\}\times\mathbb A_u^1$
is finite, since it is obtained by base change from the finite
morphism $\mathcal H$, and its target is finite over $S$.
Moreover, base change preserves the \'etale neighbourhood map and
pulls back its identification over
$Z\times\{0\}\times\mathbb A_u^1$ to an identification over
$\widetilde Z$. These data therefore define a framed homotopy
of level $2n$.

At $u=0$ equation \eqref{eq:H} gives $\mathcal H_0=\mu_l\times\id_{\A_s^n}$.
Put
\[
 V_l=\A_X^n\times_{\mu_l,\A_X^n,\pi}V,
\]
and let $r_l\colon V_l\to\A_X^n$ and
$q_l\colon V_l\to V$ be the projections.
The pulled back neighbourhood is canonically
$V_l\times\A_s^n$ with neighbourhood map
$r_l\times\id_{\A_s^n}$.
The endpoint correspondence is therefore
\[
 \bigl(
 \mu_l^{-1}(Z)\times\{0\},\;
 V_l\times\A_s^n,\;
 (\varphi_1\circ q_l,\ldots,\varphi_n\circ q_l,
  s_1,\ldots,s_n),\;
 g\circ q_l
 \bigr),
\]
where $q_l$ is composed with the projection to $V_l$
when used on $V_l\times\A_s^n$.
By definition, this is $\sigma_Y^n\circ\mu_l^*b$.

At $u=1$, equation \eqref{eq:H} gives
$\mathcal H_1=\id_{\A_X^n}\times\mu_l$, where
$\mu_l(s_1,\ldots,s_n)=(s_1^l,\ldots,s_n^l)$.
Writing $s'$ for the coordinates on the affine space factor
of the original neighbourhood $W=V\times\A_{s'}^n$, its pullback
along $\mathcal H_1$ is canonically identified with
$V\times\A_s^n$: the fibre product relation is $s'=s^l$.
Under this identification, the neighbourhood map to $\mathbb A_X^n\times\mathbb A_s^n
\cong \mathbb A_X^{2n}$  is
\[
(v,s)\longmapsto(\pi(v),s),
\]
and the projection to the original neighbourhood is
\[
(v,s)\longmapsto(v,s^l),
\qquad s^l=(s_1^l,\ldots,s_n^l).
\]
Thus the endpoint framing is
$(\varphi_1,\ldots,\varphi_n,s_1^l,\ldots,s_n^l)$,
the target map is $g\circ\pr_V$, and the support,
as a closed subset, is $Z\times\{0\}$.
These are precisely the data of
$\Theta_{l,n,Y}\circ b$.
This proves \eqref{eq:ambientpower}.
For $n=0$, both \eqref{eq:framepower} and
\eqref{eq:ambientpower} reduce to identities.
\end{proof}

\subsection{Framed naturality of relative Frobenius}\label{naturality}

Assume $\operatorname{char}k=p>0$, set $S=\Spec k$. Following~\cite[Section~33.36]{Stacks} write
\[
 T(X)=X^{(p)}=X\times_{S,F_S}S,\qquad
 F_X=F_{X/S}:X\longrightarrow X^{(p)}.
\]
Base change makes $T$ a level-preserving additive functor on
$\ZF_*(k)$. It preserves suspensions, $\A^1$-projections,
fibre products and Nisnevich coverings. Thus $\mathcal M\circ T$
is a quasi-stable $\A^1$-invariant framed presheaf whenever
$\mathcal M$ is. If $\mathcal M$ is a Nisnevich sheaf, so is
$\mathcal M\circ T$.

\begin{proposition}\label{prop:naturality}
For every $\A^1$-invariant quasi-stable additive $\ZF_*(k)$-presheaf
$\mathcal M$ of $R$-modules, the maps $F_X^*$ constitute a morphism of
framed presheaves $\mathcal M\circ T\longrightarrow\mathcal M$.
\end{proposition}

\begin{proof}
For every $X\in\Sm/k$, the morphism
$F_X:X\to X^{(p)}$ induces an $R$-linear homomorphism
\[
 F_X^*:\mathcal M(X^{(p)})\longrightarrow\mathcal M(X).
\]
We must show that these homomorphisms commute with the action
of framed correspondences. Thus, for every
$c\in\Fr_n(X,Y)$, we must prove
\[
 c^*\circ F_Y^*=F_X^*\circ T(c)^*
 \colon \mathcal M(Y^{(p)})\longrightarrow\mathcal M(X).
\]

Choose an explicit representative
$c=(Z,V,\varphi,g)$, and denote its \'etale neighbourhood map by
$e:V\to P:=\A_X^n$.
Put
\[
 c_1=F_Y\circ c,\qquad c_2=T(c)\circ F_X
 \quad\text{in }\Fr_n(X,Y^{(p)}).
\]
By contravariance, the required equality is precisely
$c_1^*=c_2^*$.

The correspondence $c_1$ is represented by
$c_1=(Z,V,\varphi,F_Y\circ g)$,
and hence $c_1^{[p]}=\bigl(Z,V,(\varphi_1^p,\ldots,\varphi_n^p),F_Y\circ g\bigr)$.
We shall identify this correspondence with the pullback
$\mu_p^*c_2$, where
\[
 \mu_p:P\to P,\qquad
 (t_1,\ldots,t_n)\longmapsto(t_1^p,\ldots,t_n^p).
\]

Using the standard identification
$P^{(p)}=\A_{X^{(p)}}^n$, write $\rho=F_X\times\id_{\A^n}:P\to P^{(p)}$.
The Frobenius twist of $c$ is represented by
$T(c)=\bigl(Z^{(p)},V^{(p)},\varphi^{(p)},g^{(p)}\bigr)$.
Precomposition with $F_X$ pulls back this correspondence along
$\rho$. More explicitly, set
\[
 V_2=V^{(p)}\times_{P^{(p)},\rho}P,
 \qquad
 Z_2=\rho^{-1}(Z^{(p)}),
\]
and let $q_2:V_2\to V^{(p)}$ be the projection. Then
\[
 c_2=
 \bigl(Z_2,V_2,\varphi^{(p)}\circ q_2,
                   g^{(p)}\circ q_2\bigr).
\]

The relative Frobenius of $P$ factors as $F_{P/S}=\rho\circ\mu_p$.
Indeed, on the factor $X$ it acts by $F_X$, and on each affine
coordinate it acts by $t_i\mapsto t_i^p$. Consequently, 
\[
 \mu_p^*c_2
 =
 \mu_p^*\rho^*T(c)
 =
 F_{P/S}^*T(c),
\]
where the fibre products are identified canonically.

The \'etale neighbourhood of the last correspondence is
\[
 \widetilde V
 =
 V^{(p)}\times_{P^{(p)},F_{P/S}}P.
\]
Naturality of relative Frobenius gives $e^{(p)}\circ F_{V/S}=F_{P/S}\circ e$,
so there is a canonical morphism
\[
 \alpha=(F_{V/S},e):V\longrightarrow\widetilde V.
\]
This morphism is an isomorphism because $e$ is \'etale.
To explain the application of the \'etale Frobenius theorem~\cite[Lemma~41.14.3]{Stacks}, let
$F_P^{\mathrm{abs}}$ be the absolute Frobenius of $P$.
There are canonical identifications
\[
 V^{(p)}
 \cong V\times_{P,\pr_1}P^{(p)},
 \qquad
 \widetilde V
 \cong V\times_{P,F_P^{\mathrm{abs}}}P,
\]
since $\pr_1\circ F_{P/S}=F_P^{\mathrm{abs}}$.
Under the second identification, $\alpha$ is the relative
Frobenius $F_{V/P}$, which is an isomorphism for an \'etale
morphism $V\to P$.

Let $q:\widetilde V\to V^{(p)}$ and
$r:\widetilde V\to P$ be the projections. By construction,
\[
 q\circ\alpha=F_{V/S},\qquad r\circ\alpha=e.
\]
Thus $\alpha$ identifies the pulled-back \'etale neighbourhood
with the original neighbourhood $V$ over $P$.
For each framing function, naturality of relative Frobenius
gives
\[
 \varphi_i^{(p)}\circ F_{V/S}=\varphi_i^p,
\]
where $(\A^1)^{(p)}$ is identified with $\A^1$ using its
standard coordinate. Similarly,
\[
 g^{(p)}\circ F_{V/S}=F_Y\circ g.
\]
The common zero set of the pulled-back framing is therefore
\[
 Z(\varphi_1^p,\ldots,\varphi_n^p)
 =
 Z(\varphi_1,\ldots,\varphi_n)=Z
\]
as a closed subset of $V$. Since $\alpha$ is an isomorphism
over $P$, this also identifies the ambient supports.
Hence, under the canonical identification of neighbourhoods,
all the correspondence data agree, proving
\begin{equation}\label{eq:exact}
 c_1^{[p]}=\mu_p^*c_2.
\end{equation}

Put $Y'=Y^{(p)}$ and
$\Theta=\Theta_{p,n,Y'}$.
Apply Lemma~\ref{lem:power} with $l=p$:
equation~\eqref{eq:framepower}, applied to $c_1$, gives
\[
 \sigma_{Y'}^n\circ c_1^{[p]}
 \simeq \Theta\circ c_1,
\]
whereas equation~\eqref{eq:ambientpower}, applied to $c_2$,
gives
\[
 \sigma_{Y'}^n\circ\mu_p^*c_2
 \simeq \Theta\circ c_2.
\]
By \eqref{eq:exact}, the left-hand sides agree. Thus
$\Theta\circ c_1\simeq\Theta\circ c_2$.
Applying $\mathcal M$, one gets
$c_1^*\circ\Theta^*=c_2^*\circ\Theta^*$.

By Lemma~\ref{lem:unit}, the endomorphism
$\Theta^*:\mathcal M(Y')\to\mathcal M(Y')$
is invertible. This uses quasi-stability of $\mathcal M$
and the invertibility of $p$ in $R$.
Composing the preceding equality on the right with
$(\Theta^*)^{-1}$ gives $c_1^*=c_2^*$.
By the definitions of $c_1$ and $c_2$, this is
\[
 c^*\circ F_Y^*=F_X^*\circ T(c)^*.
\]
The argument also applies when $n=0$, with
$\Theta_{p,0,Y'}=\id_{Y'}$.
Finally, additivity extends the equality from explicit
framed correspondences to all linear framed correspondences.
Thus the maps $F_X^*$ define the asserted morphism
$\mathcal M\circ T\to\mathcal M$ of framed presheaves.
\end{proof}

\subsection{The local Frobenius transfer}

\begin{lemma}\label{lem:target}
Suppose $e:Y\to\A^d_k$ is \'{e}tale.
Let $c_j=(Z,V,\varphi,g_j)\in\Fr_n(X,Y)$, $j=0,1$,
have exactly the same neighbourhood and framing, and suppose
$g_0|_{Z_{\mathrm{red}}}=g_1|_{Z_{\mathrm{red}}}$.
Then $c_0$ and $c_1$ are connected by a framed
$\A^1$-homotopy at level $n$.
\end{lemma}

\begin{proof}
Write $\pi:V\longrightarrow\A_X^n$
for the \'etale neighbourhood map. We regard
$Z_{\mathrm{red}}$ as a closed subscheme of $V$ through its
specified lift to this neighbourhood. By assumption,
\[
 g_Z:=g_0|_{Z_{\mathrm{red}}}
     =g_1|_{Z_{\mathrm{red}}}
 \colon Z_{\mathrm{red}}\longrightarrow Y.
\]

Put $B=V\times\A^1_u$, and define
\[
 h:B\longrightarrow\A^d_k,\qquad
 h(v,u)=(1-u)e(g_0(v))+u\,e(g_1(v)).
\]
Here the linear combination is taken coordinatewise in
$\A^d_k$. Form the fibre product
$W=B\times_{h,\A^d_k,e}Y$
with projections $q:W\rightarrow B$, $G:W\rightarrow Y$.
The morphism $q$ is separated and \'etale, since it is a base
change of the \'etale separated morphism $e$.

Let $S=Z_{\mathrm{red}}\times\A^1_u\subset B$.
The equality of the restrictions of $g_0$ and $g_1$ gives
$h|_S=e\circ g_Z\circ\pr_{Z_{\mathrm{red}}}$.
Consequently, the universal property of the fibre product
defines a morphism
\[
 s:S\longrightarrow W,\qquad
 (z,u)\longmapsto\bigl((z,u),g_Z(z)\bigr),
\]
whose composite with $q$ is the inclusion $S\hookrightarrow B$.
Equivalently, $s$ is a section of
\[
 q_S:W_S:=W\times_B S\longrightarrow S.
\]
Since $q_S$ is separated and \'etale, its section $s$ is an
open and closed immersion~\cite[Proposition~41.6.1]{Stacks}.
Thus $s(S)$ is open and closed in $W_S$.
Its complement $D=W_S\setminus s(S)$
is closed in $W_S$, hence closed in $W$, because $S$ is closed
in $B$. Set $W'=W\setminus D$.
Then $q|_{W'}$ is \'etale and
$W'\times_B S=s(S)\cong S$.
Together with the original neighbourhood map $\pi$, this makes
$W'$ an \'etale neighbourhood of
$Z_{\mathrm{red}}\times\A^1_u$ in
$\A^n_{X\times\A^1_u}$.

Pull back the framing functions from $V$:
$\Phi_i
 =\varphi_i\circ\pr_V\circ q|_{W'}$, $1\leq i\leq n$.
Their common zero set is exactly $s(S)$. Indeed, the common
zero set of the functions $\varphi_i\circ\pr_V$ on $B$ is
$S$ as a closed subset, and the construction of $W'$ removes
all of its inverse image except the chosen section.
The support $Z_{\mathrm{red}}\times\A^1_u$ is finite over
$X\times\A^1_u$, since $Z_{\mathrm{red}}$ is finite over $X$.
Therefore
\[
 \bigl(
 Z\times\A^1_u,\,
 W',\,
 (\Phi_1,\ldots,\Phi_n),\,
 G|_{W'}
 \bigr)
\]
defines a framed correspondence from $X\times\A^1_u$ to $Y$.

For $j=0,1$, we have
$h|_{V\times\{j\}}=e\circ g_j$.
Hence $g_j$ defines a section
\[
 s_j:V\longrightarrow W_j:=W\times_{\A^1_u}\{j\},
 \qquad
 v\longmapsto\bigl((v,j),g_j(v)\bigr)
\]
of $q_j:W_j\to V$.
On $Z_{\mathrm{red}}$, this is the section used to define
$W'$. Thus $s_j(V)$ avoids $D$ and lies in
$W'_j=W'\times_{\A^1_u}\{j\}$.
Since $q_j$ is \'etale, $s_j$ identifies $V$ with an open
subscheme of $W'_j$ containing the support.

Restricting the endpoint correspondence to this open
neighbourhood, its ambient map is $\pi$, its framing is
$\varphi$, and its target map is $g_j$, because
$\pr_V\circ q\circ s_j=\id_V$, $G\circ s_j=g_j$.
The endpoint at $u=j$ therefore represents
$(Z,V,\varphi,g_j)=c_j$.
The constructed correspondence is a framed
$\A^1$-homotopy from $c_0$ to $c_1$.
\end{proof}

\begin{theorem}\label{prop:Fiso}
If $\mathcal M$ is an $\A^1$-invariant quasi-stable additive framed
Nisnevich sheaf of $R$-modules, then
$F^*:\mathcal M\circ T\to\mathcal M$ is an isomorphism
of framed Nisnevich sheaves.
\end{theorem}

\begin{proof}
By Proposition~\ref{prop:naturality}, the maps $F_X^*$ define
a morphism of framed presheaves
$F^*:\mathcal M\circ T\rightarrow\mathcal M$.
Both presheaves are Nisnevich sheaves, hence Zariski sheaves.
It therefore suffices to prove that
\[
 F_U^*:\mathcal M(U^{(p)})\longrightarrow\mathcal M(U)
\]
is an isomorphism on a basis for the Zariski topology.

Every smooth $k$-scheme has a basis of affine open subschemes
$U$ admitting an \'etale morphism
\[
 e=(e_1,\ldots,e_d):U\longrightarrow\A_k^d
\]
by \cite[Tag 054L]{Stacks}.
Fix such $U$ and $e$.

The naturality square
\begin{equation}\label{square}
\begin{array}{ccc}
 U & \xrightarrow{\ F_{U/k}\ } & U^{(p)}\\
 {\scriptstyle e}\downarrow &&
 \downarrow{\scriptstyle e^{(p)}}\\
 \A_k^d & \xrightarrow{\ \mu_p\ } & \A_k^d
\end{array}
\end{equation}
is cartesian because $e$ is \'etale
\cite[Tag 0EBS]{Stacks} (also see the proof of Proposition~\ref{prop:naturality}).
Here we use the standard identification
$(\A_k^d)^{(p)}=\A_k^d$, under which relative Frobenius is
$\mu_p(t_1,\ldots,t_d)=(t_1^p,\ldots,t_d^p)$.
So,
\[
 j=(F_{U/k},e):U\longrightarrow
 P:=\A^d_{U^{(p)}}=U^{(p)}\times_k\A_k^d
\]
identifies $U$ with the closed subscheme defined by
$ \psi_i=t_i^p-e_i^{(p)}$, $1\leq i\leq d$.
Write $Z=j(U)\subset P$.

The projection $Z\to U^{(p)}$, identified with $F_{U/k}$, is finite
locally free of degree $p^d$. To see this explicitly, write
$U^{(p)}=\Spec B$. Its coordinate algebra is
\[
 B[t_1,\ldots,t_d]/
 (t_1^p-e_1^{(p)},\ldots,t_d^p-e_d^{(p)}),
\]
which is free over $B$ with basis $t_1^{a_1}\cdots t_d^{a_d}$, $0\leq a_i<p$.
Moreover, $Z$ is reduced, since $Z\cong U$ and $U\in\Sm/k$.

We next construct an \'etale neighbourhood of $Z$ in $P$
with a target map to $U$. Form
\[
 V_0=P\times_{\pr_{\A^d},\A_k^d,e}U
 \cong U^{(p)}\times_k U.
\]
Under this identification, the projection
$\nu_0:V_0\rightarrow P$
is given by $(u',v)\mapsto(u',e(v))$.
It is \'etale and separated: it is a base change of $e$,
and $e$ is separated because $U$ is separated over $k$.

Over $Z\cong U$, the morphism $\nu_0$ has the section
\[
 s:U\longrightarrow V_0,\qquad
 u\longmapsto(F_U(u),u),
\]
since $\nu_0\circ s=j$.
The image of this section is open and closed in
$V_0\times_P Z$ by \cite[Tag 024T]{Stacks}.
Delete its complementary closed subscheme from $V_0$.
The resulting open subscheme $V\subset V_0$ satisfies
$V\times_P Z\cong Z$.
Thus the restriction $\nu:V\to P$ is an \'etale neighbourhood
of $Z$.

Let $g:V\to U$ be the restriction of the second projection
$U^{(p)}\times_k U\to U$. By construction, $g\circ s=\id_U$.
Composing the functions $\psi_i$ with $\nu$, define
\[
 \tr_U=
 \bigl(Z,V,(\psi_1\circ\nu,\ldots,\psi_d\circ\nu),g\bigr)
 \in\Fr_d(U^{(p)},U).
\]
It is a framed correspondence as $Z$ is finite over
$U^{(p)}$ and its framing has common zero scheme
$V\times_P Z\cong Z$.

The composite $F_U\circ\tr_U\in\Fr_d(U^{(p)},U^{(p)})$ has the same support, neighbourhood and framing
as $\tr_U$, with target map $F_U\circ g:V\to U^{(p)}$. Consider the morphism
$a=\pr_{U^{(p)}}\circ\nu:V\rightarrow U^{(p)}$.
On the support, these two target maps agree, since
$(F_U\circ g)\circ s=F_U=a\circ s$.
The target $U^{(p)}$ admits the \'etale coordinate map
$e^{(p)}:U^{(p)}\to\A_k^d$.
Lemma~\ref{lem:target} therefore gives
\[
 F_U\circ\tr_U
 \simeq
 \bigl(Z,V,\psi\circ\nu,a\bigr).
\]

Since $a=\pr_{U^{(p)}}\circ\nu$, the explicit correspondence
$\bigl(Z,V,\psi\circ\nu,a\bigr)$
is equivalent to $ \bigl(Z,P,\psi,\pr_{U^{(p)}}\bigr)$.
Indeed, the two étale neighbourhoods of $Z\subset P$ are
$\nu:V\to P$ and $\id_P:P\to P$. Their fibre product
$V\times_P P$ is canonically $V$, with projections
$\id_V$ and $\nu$, and contains the chosen lift of $Z$.
On this common étale neighbourhood, the framing maps agree,
since $\psi\circ\nu$ is the framing of the first
correspondence, and the target maps agree, since
$a=\pr_{U^{(p)}}\circ\nu$. Thus the two tuples represent
the same framed correspondence.

Consider now the family of framings on
$P\times\A^1_\lambda$
\[
 \psi_{i,\lambda}=t_i^p-\lambda e_i^{(p)},
 \qquad 1\leq i\leq d,
\]
with target map the projection to $U^{(p)}$.
Its common zero scheme is finite locally free of degree
$p^d$ over $U^{(p)}\times\A^1_\lambda$, by the same
monomial basis argument as above.
It therefore defines a framed homotopy.
At $\lambda=1$ it gives the preceding correspondence,
and at $\lambda=0$ it gives $\Theta_{p,d,U^{(p)}}$.
Thus
\[
 F_U\circ\tr_U\simeq\Theta_{p,d,U^{(p)}}.
\]

Next, we want to describe $\tr_U\circ F_U$ up to homotopy.
Precomposition with $F_U$ is base change on the source.
Its ambient space is
\[
 P_1=\A_U^d
 \xrightarrow{\ \rho=F_U\times\id_{\A^d}\ }P,
\]
and its \'etale neighbourhood is $V_1=V\times_{P,\rho}P_1$.
Write $\beta:V_1\to V$ and $\nu_1:V_1\to P_1$ for the
projections. Its target map is $g_1=g\circ\beta:V_1\to U$.
By naturality of relative Frobenius~\eqref{square}, $e_i^{(p)}\circ F_U=e_i^p$.
Consequently, its framing is the pullback along $\nu_1$ of
\[
 t_i^p-e_i^p=(t_i-e_i)^p,
 \qquad 1\leq i\leq d.
\]

The scheme theoretic inverse image $Z_1=\rho^{-1}(Z)$
identifies with $U\times_{U^{(p)}}U$,
where the first factor is the source of the correspondence.
The relative Frobenius $F_U=F_{U/k}$ is a universal
homeomorphism by~\cite[Tag 0CCB]{Stacks}. In particular,
it is radicial. Since $F_U$ is radicial, \cite[Tag 01S4]{Stacks} implies
that the diagonal $\Delta:U\rightarrow U\times_{U^{(p)}}U$
is surjective on underlying topological spaces.
It is also a closed immersion as $F_U$ is finite, hence separated.
Since $U$ is reduced, $\Delta$
identifies $U$ with $(Z_1)_{\mathrm{red}}$.

Under the identification $Z_1\cong U\times_{U^{(p)}}U$,
the closed immersion $Z_1\hookrightarrow P_1=U\times_k\A_k^d$
is induced by 
\[
 (\pr_1,e\circ\pr_2):
 U\times_{U^{(p)}}U\longrightarrow U\times_k\A_k^d.
\]
Its composite with the diagonal $\Delta$ is
$(\id_U,e)$. Since $\Delta$ identifies $U$ with
$(Z_1)_{\mathrm{red}}$, the reduced support in $P_1$
is therefore the graph of $e$, defined by
$t_i-\pr_U^*(e_i)=0$, $1\leq i\leq d$,
where $\pr_U:P_1\to U$ is the projection.

Let $a_1=\pr_U\circ\nu_1:V_1\to U$ be the structural
projection, and let $s_1:Z_1\to V_1$ be the chosen lift
of the support. Under the identification
$Z_1\cong U\times_{U^{(p)}}U$, the construction gives
\[
 a_1\circ s_1=\pr_1,
 \qquad
 g_1\circ s_1=\pr_2.
\]
Indeed, $a_1$ projects to the new source, whereas $g_1$
restricts to the point of the original support
$Z\cong U$, on which $g$ is the identity.
Since the diagonal $\Delta:U\to Z_1$ identifies $U$
with $(Z_1)_{\mathrm{red}}$, and
\[
 \pr_1\circ\Delta=\id_U=\pr_2\circ\Delta,
\]
we obtain
\[
 a_1\circ s_1\circ\Delta
 =\id_U
 =g_1\circ s_1\circ\Delta.
\]
Thus $a_1$ and $g_1$ agree on the reduced support
in $V_1$. Applying Lemma~\ref{lem:target} with target $U$ gives
   $$\tr_U\circ F_U\simeq(|Z_1|,V_1,(t_i^p-e_i^p)_{i=1}^d\circ\nu_1,a_1).$$
As in the first composite, the framing and target
map are pulled back from $P_1$. Thus this correspondence is
equivalent to the one represented on all of $P_1$ with
framing $(t_i^p-e_i^p)_{i=1}^d$ and target map $\pr_U$.

On $P_1\times\A^1_\lambda=\A^d_{U\times\A^1_\lambda}$,
consider the framing $\chi_\lambda
 =(t_1^p-\lambda e_1^p,\ldots,t_d^p-\lambda e_d^p)$
and the target map $\pr_U$. Write $U=\Spec A$.
The common zero scheme $\mathcal Z$ has coordinate ring
\[
 A[\lambda,t_1,\ldots,t_d]/
 (t_1^p-\lambda e_1^p,\ldots,t_d^p-\lambda e_d^p),
\]
which is free over $A[\lambda]$ with basis
$t_1^{a_1}\cdots t_d^{a_d}$, $0\leq a_i<p$.
Thus $\mathcal Z$
is finite locally free of degree $p^d$ over
$U\times\A^1_\lambda$. Taking $P_1\times\A^1_\lambda$
as the \'etale neighbourhood, we obtain a framed homotopy.
At $\lambda=1$ its restriction is
$\bigl(Z_1,P_1,(t_i^p-e_i^p)_{i=1}^d,\pr_U\bigr)$,
whereas at $\lambda=0$ its framing is
$(t_1^p,\ldots,t_d^p)$ and its support is the zero section.
The latter endpoint is precisely $\Theta_{p,d,U}$.

We have therefore proved, at framing level $d$,
\begin{equation}\label{eq:composites}
 F_U\circ\tr_U\simeq\Theta_{p,d,U^{(p)}},
 \qquad
 \tr_U\circ F_U\simeq\Theta_{p,d,U}.
\end{equation}

Set
\[
 f=F_U^*:\mathcal M(U^{(p)})\to\mathcal M(U),
 \qquad
 t=\tr_U^*:\mathcal M(U)\to\mathcal M(U^{(p)}).
\]
Since $\mathcal M$ is additive and $\A^1$-invariant,
equations~\eqref{eq:composites} imply
\[
 t\circ f=\Theta_{p,d,U^{(p)}}^*,
 \qquad
 f\circ t=\Theta_{p,d,U}^*.
\]
Both right-hand sides are invertible by
Lemma~\ref{lem:unit}, using quasi-stability and
$p\in R^\times$. Hence
\[
 \bigl(\Theta_{p,d,U^{(p)}}^*\bigr)^{-1}\circ t
 \quad\text{and}\quad
 t\circ\bigl(\Theta_{p,d,U}^*\bigr)^{-1}
\]
are respectively a left and a right inverse of $f$.
They coincide, so $F_U^*$ is an isomorphism.
\end{proof}

\subsection{Descent of the sheaf to the perfect closure}\label{descent}

We make the indexing and the change of base field explicit.
Put $K=k^{\perf}$ and $k_n=k^{1/p^n}\subset K$, so that
$k_0=k$ and $K=\bigcup_{n\geq0}k_n$. For each $n$, the map
\[
 \sigma_n:k\xrightarrow{\sim}k_n,\qquad
 a\longmapsto a^{1/p^n},
\]
is a field isomorphism with inverse $b\mapsto b^{p^n}$.
Descent to a finite stage means descent to $k_n$ for some
finite index $n$. It does not mean that $k_n/k$ is a finite
extension. No finiteness assumption on $[k:k^p]$ is made.

Let $E_n:\ZF_*(k)\longrightarrow\ZF_*(k_n)$
be the base change equivalence along $\sigma_n$, and let
$\rho_n$ be its canonical quasi-inverse, given by base change
along $\sigma_n^{-1}$. Write $i_n:k_n\hookrightarrow k_{n+1}$
for the inclusion and
$b_n:\ZF_*(k_n)\to\ZF_*(k_{n+1})$
for the corresponding base-change functor. Recall that $T$
is base change along the Frobenius homomorphism
$F_k:k\to k$, $a\mapsto a^p$. There are canonical natural
isomorphisms of functors
\begin{equation}\label{eq:transport}
 \rho_{n+1}b_n\cong T\rho_n,\qquad
 \rho_n(X_{k_n})\cong T^nX.
\end{equation}

For the first isomorphism, observe that
\[
 (\sigma_{n+1}^{-1}\circ i_n\circ\sigma_n)(a)
 =\bigl(a^{1/p^n}\bigr)^{p^{n+1}}
 =a^p.
\]
Consequently, $\sigma_{n+1}^{-1}\circ i_n
 =F_k\circ\sigma_n^{-1}$,
and associativity of base change gives
$\rho_{n+1}b_n\cong T\rho_n$.

For the second isomorphism, let $j_n:k\hookrightarrow k_n$
be the inclusion used to define $X_{k_n}$. This inclusion
sends $a$ to $a$, whereas $\sigma_n$ sends $a$ to
$a^{1/p^n}$. Hence $(\sigma_n^{-1}\circ j_n)(a)=a^{p^n}$.
Thus $\rho_n(X_{k_n})$ is canonically the base change of $X$
along $F_k^n$, namely $T^nX$. All these identifications
apply to the support, étale neighbourhood, framing and
target map of a framed correspondence. They are compatible
with composition and the additivity relations, and hence
give the asserted natural isomorphisms on the categories
of linear framed correspondences.

\begin{proposition}\label{prop:descent}
Suppose that $\mathcal A$ is
an $\A^1$-invariant quasi-stable additive framed
Nisnevich sheaf of $R$-modules over $k$. Then there exists an $\A^1$-invariant quasi-stable
additive framed Nisnevich sheaf $\mathcal A_K$
of $R$-modules over $K$, together with a natural
isomorphism of framed Nisnevich sheaves 
$\mathcal A\xrightarrow{\sim}\mathcal A_K\circ B_K$ over $k$,
where $B_K:\ZF_*(k)\to\ZF_*(K)$ is base extension. In particular, there are isomorphisms
\begin{equation}\label{eq:sections}
 \mathcal A(X)\xrightarrow{\sim}\mathcal A_K(X_K),
 \qquad X\in\Sm/k,
\end{equation}
compatible with all linear framed correspondences over $k$.
\end{proposition}

\begin{proof}
For each $n\geq0$, define $\mathcal B_n=\mathcal A\circ\rho_n$
on $\ZF_*(k_n)$. Nisnevich covers, products with the affine line, disjoint unions,
and suspension correspondences are preserved along the field isomorphism
$\sigma_n^{-1}:k_n\to k$. Consequently, each $\mathcal B_n$ is an $\A^1$-invariant
quasi-stable additive framed Nisnevich sheaf.

Using \eqref{eq:transport} and
Theorem~\ref{prop:Fiso}, relative Frobenius gives
natural isomorphisms of framed presheaves
\[
 \alpha_n:
 \mathcal B_{n+1}\circ b_n
 \xrightarrow{\sim}\mathcal A\circ T\circ\rho_n
 \xrightarrow{F^*}\mathcal A\circ\rho_n
 =\mathcal B_n.
\]
Here the component of the second arrow at
$Y\in\Sm/k_n$ is $F_{\rho_n(Y)/k}^*:\mathcal A(T\rho_n(Y))\xrightarrow{\sim}\mathcal A(\rho_n(Y))$.

For $m\geq n$, let $b_{n,m}$ denote base extension
along $k_n\hookrightarrow k_m$, with
$b_{n,n}=\mathrm{Id}$. Set
$\alpha_{n,n}=\mathrm{id}_{\mathcal B_n}$ and define
recursively
\[
 \alpha_{n,m+1,Y}
 =
 \alpha_{n,m,Y}\circ\alpha_{m,b_{n,m}Y}.
\]
This gives natural isomorphisms of framed presheaves
$\alpha_{n,m}:\mathcal B_m\circ b_{n,m}\xrightarrow{\sim}\mathcal B_n$.
For $n\leq m\leq r$, they satisfy
\[
 \alpha_{n,r,Y}
 =
 \alpha_{n,m,Y}\circ
 \alpha_{m,r,b_{n,m}Y}.
\]
Indeed, both sides are the same ordered composite
of the consecutive comparisons from stage $r$
down to stage $n$. Naturality ensures compatibility
with the canonical base-change identifications.

Since $K=\varinjlim_n k_n$, schemes of
finite presentation over $K$, morphisms between
them, and equalities between such morphisms descend
to a common stage --- see \cite[Tag 01ZM]{Stacks}.
Thus every $Y\in\Sm/K$ admits a model
$Y_n\in\Sm/k_n$, together with an isomorphism
$(Y_n)_K\cong Y$. Smoothness and separatedness of the descended model
follow by faithfully flat descent along $k_n\subset K$,
since both properties are fpqc local on the target ---
see \cite[Tags 02VL and 02KU]{Stacks}, respectively.

The same argument applies to explicit framed
correspondences. Namely, let $c=(Z,V,\varphi,g)\in\Fr_r^K(X,Y)$.
The support $Z$ is finite over $X$, hence
quasi-compact. Restricting $V$ to a quasi-compact
open neighbourhood of the chosen lift of $Z$
does not change the represented correspondence.
We may therefore assume that $V$ is quasi-compact.
Give $Z\subset\A_X^r$ its reduced induced scheme
structure. The schemes involved are then of finite
presentation over $K$.

The closed immersion $Z\hookrightarrow\A_X^r$,
the étale neighbourhood map $V\to\A_X^r$, the
chosen lift $Z\to V$, the framing
$\varphi:V\to\A_K^r$, and the target map $g:V\to Y$
descend, together with their commutative diagrams,
to a common stage. These properties can be checked after the faithfully
flat quasi-compact base change $k_n\subset K$.
Indeed, finiteness, being a closed immersion,
étaleness, and being an isomorphism are fpqc local
on the target; see
\cite[Tags 02LA, 02L6, 02VN, and 02L4]{Stacks},
respectively. Applying these results to
\[
 Z_n\to X_n,\qquad
 Z_n\to\A^r_{X_n},\qquad
 V_n\to\A^r_{X_n},\qquad
 Z_n\to V_n\times_{\A^r_{X_n}}Z_n
\]
verifies the required properties at the chosen stage.

The zero set condition can also be imposed at a
finite stage. On a finite affine cover of $V$, let
$I$ be the ideal of the lifted support and put
$J=(\varphi_1,\ldots,\varphi_r)$. The condition is $\sqrt I=\sqrt J$.
Writing $I=(f_1,\ldots,f_s)$, this equality is
implied by finitely many relations $f_a^{N_a}\in J,
 \varphi_i^{M_i}\in I$ for suitable positive integers $N_a$ and $M_i$.
The coefficients expressing these relations
descend after increasing the stage. The resulting
framing therefore has exactly the prescribed
closed support. Thus $c$ descends to an explicit
framed correspondence over some $k_n$. If models
of $X$ and $Y$ have already been fixed, this
descent can be carried out over their base changes
to a sufficiently large common stage.

Equivalence of explicit representatives also
descends. A common étale neighbourhood expressing
equivalence may be restricted around the support
to be quasi-compact. Its maps, and the equalities
of the two framing and target maps on it, then
descend to a common stage. Finite disjoint
decompositions of supports descend as well.
Consequently, the additivity relations defining
linear framed correspondences descend. Since
every linear combination and every expression
in these relations is finite, linear framed
correspondences and equalities between them
descend to a sufficiently large common stage.
These constructions are compatible with
composition, which is defined using fibre
products and the corresponding framing and
target maps.

Now choose a model $Y_n$ of $Y\in\Sm/K$ and put
$Y_m=b_{n,m}(Y_n)$ for $m\geq n$. Define
\[
 \mathcal A_K(Y)
 =
 \varinjlim_{m\geq n}\mathcal B_m(Y_m),
\]
with transition maps $\alpha_{m,Y_m}^{-1}:
 \mathcal B_m(Y_m)\xrightarrow{\sim}
 \mathcal B_{m+1}(Y_{m+1})$.
The coherence identity for the $\alpha_{n,m}$
shows that the transition from stage $m$ to
stage $r$ is $\alpha_{m,r,Y_m}^{-1}$.
In particular, every transition map is an
isomorphism.

This definition is independent of the chosen
model up to canonical isomorphism. Indeed, for
two models equipped with identifications of
their base changes with $Y$, the induced
isomorphism over $K$, its inverse, and their
inverse identities descend to a common stage.
Applying the corresponding $\mathcal B_m$ and
using naturality of the comparisons identifies
the two colimits. Any two descents of this
isomorphism agree at a later stage, so the
identification is independent of the descended
comparison. For three models these identifications satisfy the cocycle identity,
since the corresponding identity of morphisms holds after passage to a common stage.

We next define the framed action. Let
$c\in\ZF_r^K(X,Y)$. Choose models of $X$ and $Y$
and, after increasing the stage, a descent $c_n\in\ZF_r^{k_n}(X_n,Y_n)$.
For $m\geq n$, write $c_m=b_{n,m}(c_n)$. Naturality of $\alpha_m$ gives
\[
 \alpha_{m,X_m}^{-1}\circ c_m^*
 =
 c_{m+1}^*\circ\alpha_{m,Y_m}^{-1}.
\]
Hence the maps $c_m^*:\mathcal B_m(Y_m)\to
 \mathcal B_m(X_m)$ induce an $R$-linear map
$c^*:\mathcal A_K(Y)\longrightarrow\mathcal A_K(X)$.
The descent of representatives and relations
proved above shows that this map is independent
of all choices. Identities, composition, and
addition may be checked at a common finite
stage, where they follow from the corresponding
properties of $\mathcal B_m$. Thus $\mathcal A_K$
is a linear framed presheaf. It is additive, because
disjoint unions admit disjoint union models,
and filtered colimits commute with finite direct sums.

We verify the Nisnevich sheaf condition.
Let $\{U_i\to Y\}_{i=1}^r$ be a finite Nisnevich
cover over $K$. Descend the schemes and maps to
a common stage, obtaining étale morphisms
$U_{i,n}\to Y_n$. Their pairwise fibre products recover those of
the given cover after base change to $K$.
This family is already a Nisnevich cover over
$k_n$. To see this, fix $y\in Y_n$. Since $K/k_n$
is purely inseparable, $Y\to Y_n$ is a universal
homeomorphism. Let $y'\in Y$ be the unique point above $y$. The extension
$\kappa(y')/\kappa(y)$ is purely inseparable.
Choose $u'\in U_i$ above $y'$ such that
$\kappa(y')\to\kappa(u')$ is an isomorphism,
and let $u\in U_{i,n}$ be its image.
The extension $\kappa(u)/\kappa(y)$ is finite
separable because $U_{i,n}\to Y_n$ is étale.
Moreover, there is an embedding
\[
 \kappa(u)\hookrightarrow\kappa(u')
 \cong\kappa(y')
\]
over $\kappa(y)$. Since a finite separable
subextension of a purely inseparable extension
is trivial, $\kappa(y)\to\kappa(u)$ is an
isomorphism. This proves the claim.

For $m\geq n$, put
$U_{i,m}=(U_{i,n})_{k_m}$ and $Y_m=(Y_n)_{k_m}$.
Base change preserves Nisnevich covers, so the
sheaf condition for $\mathcal B_m$ gives
\[
 \mathcal B_m(Y_m)\longrightarrow
 \prod_{i=1}^r\mathcal B_m(U_{i,m})
 \rightrightarrows
 \prod_{i,j=1}^r
 \mathcal B_m(U_{i,m}\times_{Y_m}U_{j,m})
\]
as an equalizer diagram. These diagrams are
compatible with the transition maps by
naturality of $\alpha_m$. Filtered colimits of
$R$-modules commute with finite limits.
Taking the colimit therefore gives the sheaf
equalizer for $\mathcal A_K$ and the original
cover. Finite covers suffice for the Nisnevich
topology here. So the preceding
equalizer argument proves that
$\mathcal A_K$ is a Nisnevich sheaf.

To check $\A^1$-invariance, choose a model $Y_n$
of $Y$ and use $Y_n\times\A^1_{k_n}$ as a model
of $Y\times\A^1_K$. The projection maps induce
isomorphisms $\mathcal B_m(Y_m)\xrightarrow{\sim}
 \mathcal B_m(Y_m\times\A^1_{k_m})$
at every stage. Their colimit is the required
isomorphism $\mathcal A_K(Y)\xrightarrow{\sim}
 \mathcal A_K(Y\times\A^1_K)$.
Likewise, the suspension correspondence
$\sigma_Y$ is the base change of
$\sigma_{Y_m}$ at each stage. Since
$\sigma_{Y_m}^*$ is an automorphism of
$\mathcal B_m(Y_m)$ and these automorphisms
commute with the transition maps,
$\sigma_Y^*$ is an automorphism of
$\mathcal A_K(Y)$. Thus $\mathcal A_K$ is
quasi-stable.

Finally, for $Y=X_K$ choose the model
$X/k=k_0$. The initial term of the defining
colimit is canonically
$\mathcal B_0(X)\cong\mathcal A(X)$.
Since every transition is an isomorphism,
the canonical map from this initial term
gives~\eqref{eq:sections}. If
$c\in\ZF_r^k(X,Y)$, its base change to $K$
is represented throughout the construction
by the base changes of $c$. Hence \eqref{eq:sections} commutes with the action
of every linear framed correspondence over $k$,
proving the asserted isomorphism of framed
Nisnevich sheaves.
\end{proof}

\subsection{Completion of the proof}
The following lemma compares big and small Nisnevich cohomology.

\begin{lemma}
Put $\mathscr C_k=(\Sm/k)_{\Nis}$. Let $\mathcal B$ be an abelian
sheaf on $\mathscr C_k$. For $X\in\Sm/k$, let $\mathcal B_X$ be
its restriction to the small Nisnevich site of $X$, defined on
the basis of separated \'{e}tale $X$-schemes of finite type by
$\mathcal B_X(U\to X)=\mathcal B(U)$.
There is a canonical isomorphism in $D(\mathrm{Ab})$
\begin{equation}\label{eq:big-small-comparison}
 \chi_X:\ R\Gamma_{\mathscr C_k}(X,\mathcal B)
 \xrightarrow{\simeq}R\Gamma(X_{\Nis},\mathcal B_X),
\end{equation}
natural in $\mathcal B$ and contravariantly functorial in
$X\in\Sm/k$.

More precisely, for any $k$-morphism $g:Y\to X$, base change
defines a morphism of small Nisnevich topoi $g_{\Nis}$ and a
canonical morphism of coefficient sheaves
$\gamma_g:g_{\Nis}^{-1}\mathcal B_X\longrightarrow\mathcal B_Y$.
The following square commutes in $D(\mathrm{Ab})$:
\[
\begin{array}{ccc}
 R\Gamma_{\mathscr C_k}(X,\mathcal B)
 &\xrightarrow{\ g^*\ }&R\Gamma_{\mathscr C_k}(Y,\mathcal B)\\[3pt]
 {\scriptstyle\chi_X}\,\big\downarrow
 &&\big\downarrow\,{\scriptstyle\chi_Y}\\[3pt]
 R\Gamma(X_{\Nis},\mathcal B_X)
 &\xrightarrow{\ g^*\ }&R\Gamma(Y_{\Nis},\mathcal B_Y).
\end{array}
\]
The lower arrow is pullback followed by the map induced by
$\gamma_g$. These maps are compatible with identities and
composition. In particular, \eqref{eq:big-small-comparison}
identifies the cohomology groups in every degree and respects
pullback along $X\times\A^1\to X$.
\end{lemma}

\begin{proof}
For a site $\mathscr C$, write $\operatorname{Ab}(\mathscr C)$
for its category of abelian sheaves. All sheaf functors
in this proof are functors between these abelian categories.
When applying references stated for sheaves of modules,
we take the coefficient ring sheaf to be
$\underline{\mathbb Z}$. Indeed, abelian sheaves on
a site are canonically the same as modules over its
constant sheaf of integers (see
\cite[Tag 079X, observation (3)]{Stacks}).

Put $\mathscr D_X=\mathscr C_k/X$. An object
of this slice site is a $k$-morphism $a:U\to X$, with
$U\in\mathrm{Sm}/k$. Consider the restriction functor
\[
j_X^{-1}:\operatorname{Ab}(\mathscr C_k)
\longrightarrow\operatorname{Ab}(\mathscr D_X),
\]
defined by $(j_X^{-1}\mathcal F)(U\xrightarrow{a}X)=\mathcal F(U)$.
This is restriction to a localized site in the sense of
\cite[Tags 00XZ and 00Y0]{Stacks}, applied to abelian sheaves.

The functor $j_X^{-1}$ is exact by \cite[Tags 00XR and 00YV]{Stacks}.
It has a left adjoint
$j_{X!}:\operatorname{Ab}(\mathscr D_X)
\longrightarrow\operatorname{Ab}(\mathscr C_k)$.
Its existence and exactness follow from
\cite[Tags 03DI and 03DJ]{Stacks}, applied to the
ringed site $(\mathscr C_k,\underline{\mathbb Z})$
and its localization at $X$. Here restriction of
$\underline{\mathbb Z}$ to the slice site is again
the constant sheaf of integers. Consequently,
$j_X^{-1}$ preserves injectives. Indeed, for an injective
$\mathcal I$, the adjunction gives
\[
\operatorname{Hom}_{\operatorname{Ab}(\mathscr D_X)}
(-,j_X^{-1}\mathcal I)
\cong
\operatorname{Hom}_{\operatorname{Ab}(\mathscr C_k)}
(j_{X!}(-),\mathcal I),
\]
and the right hand functor is exact because $j_{X!}$
is exact and $\mathcal I$ is injective.

We compute on the small Nisnevich site using its basis
of separated \'etale $X$-schemes of finite type.
For this proof, denote this basis by $\mathscr E_X$.
We identify $\operatorname{Ab}(\mathscr E_X)$ with
$\operatorname{Ab}(X_{\mathrm{Nis}})$ by restriction
to this basis. Consider the inclusion
\[
u_X:\mathscr E_X\longrightarrow\mathscr D_X.
\]
It is fully faithful, preserves fibre
products and the terminal object $(X,\mathrm{id}_X)$.
The category $\mathscr E_X$ also has equalizers, and $u_X$ preserves them. 
Moreover, $u_X$ is plainly continuous and cocontinuous. 

Write $r_X:\operatorname{Ab}(\mathscr D_X)
\longrightarrow\operatorname{Ab}(X_{\mathrm{Nis}})$
for restriction along $u_X$. Since $u_X$ is continuous and cocontinuous,
$r_X$ is the inverse image functor
of the associated morphism of topoi \cite[Tag 00XR]{Stacks}. It is therefore exact
on abelian sheaves by \cite[Tag 00YV]{Stacks}.
By \cite[Tag 04BG]{Stacks} restriction on abelian
sheaves has a left adjoint
\[
L_X:\operatorname{Ab}(X_{\mathrm{Nis}})
\longrightarrow\operatorname{Ab}(\mathscr D_X).
\]
By \cite[Tag 04BH]{Stacks} $L_X$ is exact. Thus $r_X$ preserves injectives,
by the same argument used for $j_X^{-1}$.

Choose an injective resolution
$\mathcal B\to I^\bullet$ on $\mathscr C_k$, and put $I_X^\bullet=r_Xj_X^{-1}I^\bullet$.
By the preceding arguments,
$\mathcal B_X\to I_X^\bullet$ is an injective
resolution on $X_{\mathrm{Nis}}$. Evaluation at the
terminal object $(X,\mathrm{id}_X)$ gives an equality
of complexes $\Gamma(X_{\mathrm{Nis}},I_X^\bullet)=I^\bullet(X)$.
The complex on the right represents
$R\Gamma_{\mathscr C_k}(X,\mathcal B)$, while the one
on the left represents
$R\Gamma(X_{\mathrm{Nis}},\mathcal B_X)$.
This gives the canonical comparison isomorphism
$\chi_X$ in \eqref{eq:big-small-comparison}.
The construction is independent of the chosen
injective resolution and natural in $\mathcal B$
by the comparison theorem for injective resolutions.
This is the comparison argument of~\cite[Tag 03YU]{Stacks} in our setting.

We now verify functoriality in $X$. Let $g:Y\to X$
be a $k$-morphism. For an object $U\to X$ of
$\mathscr E_X$, the projection $p_U:U\times_XY\to U$
induces $\mathcal B(U)\xrightarrow{p_U^*}\mathcal B(U\times_XY)$.
These maps define a morphism
$\delta_g:\mathcal B_X\rightarrow g_{\mathrm{Nis},*}\mathcal B_Y$,
because
\[
(g_{\mathrm{Nis},*}\mathcal B_Y)(U\to X)
=
\mathcal B(U\times_XY).
\]
The coefficient morphism $\gamma_g$ in the statement
is the adjoint of $\delta_g$.

Use the same resolution $I^\bullet$ for both $X$ and
$Y$. Applying the preceding construction in each
degree gives a morphism of complexes
$I_X^\bullet\longrightarrow
g_{\mathrm{Nis},*}I_Y^\bullet$.
The complex on the right represents
$Rg_{\mathrm{Nis},*}\mathcal B_Y$, since
$I_Y^\bullet$ is an injective resolution.
The displayed morphism represents the composite
\[
\mathcal B_X\xrightarrow{\delta_g}
g_{\mathrm{Nis},*}\mathcal B_Y
\longrightarrow Rg_{\mathrm{Nis},*}\mathcal B_Y.
\]
Under derived adjunction, this is the morphism
corresponding to $\gamma_g$
($g_{\mathrm{Nis}}^{-1}$ is exact by \cite[Tag 04JC]{Stacks}).

The direct image functor $g_{\mathrm{Nis},*}$
preserves injectives, since its left adjoint
$g_{\mathrm{Nis}}^{-1}$ is exact. Taking global
sections therefore shows that the lower arrow
in the statement is represented by $I^\bullet(X)\longrightarrow I^\bullet(Y)$,
whose map in degree $q$ is restriction along $g$.
This also represents the upper arrow, because
$R\Gamma_{\mathscr C_k}(X,-)$ is the derived functor
of evaluation at $X$. Under the two comparison
isomorphisms, the square in the statement is thus
represented by a commutative square of complexes.

Finally, restriction maps of sheaves satisfy
$\mathrm{id}^*=\mathrm{id}$ and
$(g\circ h)^*=h^*\circ g^*$ for
$Z\xrightarrow{h}Y\xrightarrow{g}X$.
The same identities therefore hold for the
derived maps just constructed. Taking cohomology
proves the asserted comparisons in every degree,
including compatibility with pullback along
$X\times\mathbb A^1\to X$.
\end{proof}

\begin{lemma}\label{polezno}
Let $f:Y\to X$ be a universal homeomorphism of schemes. Then base change induces
an equivalence of small Nisnevich sites. As a consequence, for every
abelian sheaf $\mathcal B$ on $X_{\Nis}$ and every $i\geq0$,
there is a natural isomorphism
\[
 H^i_{\Nis}(X,\mathcal B)
 \xrightarrow{\simeq}H^i_{\Nis}(Y,f_{\Nis}^*\mathcal B),
\]
where $f_{\Nis}^*$ is the inverse-image functor on sheaves.
\end{lemma}

\begin{proof}
A universal homeomorphism is integral, surjective and
radicial by~\cite[Tags 04DF and 01S4]{Stacks}.
By~\cite[Theorem 18.1.2]{EGA} or
\cite[Tag 04DZ]{Stacks} base change induces an equivalence
\[
 \operatorname{Et}(X)\xrightarrow{\sim}\operatorname{Et}(Y),
 \qquad U\longmapsto U\times_XY,
\]
between the categories of schemes étale over $X$ and $Y$.
We want to verify that this equivalence preserves and reflects
Nisnevich coverings.

Let $U\to X$ be étale, let $u\in U$ have image $x\in X$,
and let $y\in Y$ be the unique point above $x$. Put
\[
 F=\kappa(x),\qquad L=\kappa(u),\qquad E=\kappa(y).
\]
Since $U\to X$ is étale, $L/F$ is finite separable
by~\cite[Tag 02GL]{Stacks}. Since $f$ is radicial,
$E/F$ is algebraic purely inseparable.

Separable algebraic extensions and purely inseparable
algebraic extensions are linearly disjoint by
\cite[Chapter 10, 10.7(b)]{MilneAG}.
Thus after embedding $L$ and $E$ into a common algebraic
closure of $F$, the multiplication homomorphism
\[
 L\otimes_F E\longrightarrow LE,
 \qquad \ell\otimes e\longmapsto\ell e,
\]
is injective. In particular, $L\otimes_F E$ is an integral
domain. Moreover, extension of scalars carries an
$F$-basis of $L$ to an $E$-basis of $L\otimes_F E$, and so $\dim_E(L\otimes_F E)=[L:F]$.
A finite-dimensional algebra over a field that is an
integral domain is a field. Consequently, $[L\otimes_F E:E]=[L:F]$.

By~\cite[Tag 01JT]{Stacks} points of
$U\times_XY$ lying above both $u$ and $y$ correspond
to prime ideals of $L\otimes_F E$. Since $ L\otimes_F E$ is a field, there is exactly one such point
$u'$ with $\kappa(u')\cong L\otimes_F E$.
Every point above $u$ necessarily maps to $y$, because
$y$ is the unique point of $Y$ above $x$. Hence $u'$
is the unique point above $u$. The degree formula gives
\[
 [\kappa(u'):\kappa(y)]
 =
 [\kappa(u):\kappa(x)].
\]
In particular, 
$\kappa(x)\to\kappa(u)$ is an isomorphism if and only if
$\kappa(y)\to\kappa(u')$ is.

Applying this observation to each member of an étale
family $\{U_i\to X\}_{i\in I}$ shows that the family
satisfies the Nisnevich lifting condition at $x$ if and
only if its base change satisfies that condition at $y$.
Since $Y\to X$ is bijective on points, the family is a
Nisnevich covering of $X$ if and only if
$\{U_i\times_XY\to Y\}_{i\in I}$ is a Nisnevich
covering of $Y$.
The same argument applies to coverings of any object
$U\to X$ of $\operatorname{Et}(X)$: replace $f$ by its
base change $U\times_XY\longrightarrow U$,
which is again a universal homeomorphism. Therefore
the equivalence of étale categories preserves and
reflects all Nisnevich coverings. It is thus an
equivalence of the small Nisnevich sites. We refer the reader to~\cite[Tag~03SI]{Stacks}
for the same equivalence of small \'etale sites.

It follows that the inverse-image functor
$f_{\mathrm{Nis}}^*:
 \operatorname{Ab}(X_{\mathrm{Nis}})
 \xrightarrow{\sim}
 \operatorname{Ab}(Y_{\mathrm{Nis}})$
is an equivalence of abelian categories. In particular,
it is exact and preserves injective objects. Since base
change carries the final object $X$ to the final object
$Y$, this equivalence identifies global sections:
\[
 \Gamma(X_{\mathrm{Nis}},\mathcal B)
 \cong
 \Gamma(Y_{\mathrm{Nis}},f_{\mathrm{Nis}}^*\mathcal B).
\]
Applying these identifications to an injective resolution
of $\mathcal B$ gives the natural isomorphisms
\[
 H^i_{\mathrm{Nis}}(X,\mathcal B)
 \xrightarrow{\sim}
 H^i_{\mathrm{Nis}}(Y,f_{\mathrm{Nis}}^*\mathcal B),
 \qquad i\geq0.
\]
as required.
\end{proof}

We are now in a position to prove Theorem \ref{thm:main}.

\begin{proof}[Proof of Theorem \ref{thm:main}]
Put $\mathcal G=\mathcal F\otimes_{\mathbb Z}R$.
It is additive, quasi-stable and $\A^1$-invariant.
Since $k$ is imperfect, $p>1$ and $k$ is infinite.
By \cite[Theorem~1.1]{GP} $\mathcal A=\mathcal G_{\Nis}$ is a homotopy invariant quasi-stable framed sheaf. 
This sheafification result requires no perfectness assumption. In characteristic $2$ we use its
$\mathbb Z[1/2]$-coefficient version, which applies since $R=\mathbb Z[1/2]$. 
Proposition \ref{prop:descent} supplies $\mathcal A_K$ over
$K=k^{\perf}$. The field $K$ is infinite and perfect, so
\cite[Theorem 1.1]{GP} (again with $\mathbb Z[1/2]$-coefficients
in characteristic $2$) implies
\begin{equation}\label{eq:strictK}
 H^i_{\Nis}(Y,\mathcal A_K)\xrightarrow{\simeq}
 H^i_{\Nis}(Y\times\A^1_K,\mathcal A_K),
 \qquad Y\in\Sm/K.
\end{equation}
Note that $\mathcal A_K$ is already a Nisnevich sheaf, so there is no need for its sheafification.

For $X\in\Sm/k$, $\pi:X_K\to X$ is a universal
homeomorphism. The preceding lemma gives
$(X_K)_{\Nis}\simeq X_{\Nis}$, with an \'{e}tale $X$-scheme
$U$ corresponding to $U_K$. On the basis of separated \'{e}tale
$X$-schemes of finite type the natural isomorphisms
\eqref{eq:sections} therefore
identify the restrictions of $\mathcal A$ and $\mathcal A_K$
to these small sites. Combining the cohomological assertion of
the preceding lemma with \eqref{eq:big-small-comparison} over
$k$ and over $K$ gives
\[
 H^i_{\mathscr C_k}(X,\mathcal A)
 \cong H^i(X_{\Nis},\mathcal A_X)
 \cong H^i((X_K)_{\Nis},(\mathcal A_K)_{X_K})
 \cong H^i_{\mathscr C_K}(X_K,\mathcal A_K).
\]
These isomorphisms are natural in $X$. In particular, for
$q:X\times\A^1\to X$, functoriality of
\eqref{eq:big-small-comparison} and naturality of
\eqref{eq:sections} identify the cohomological pullback by
$q$ with that by $q_K:X_K\times\A^1_K\to X_K$.
Applying the comparisons to
$X$ and $X\times\A^1$, together with
\eqref{eq:strictK} and
$(X\times\A^1)_K=X_K\times\A^1_K$, gives the required
isomorphism over $k$ in every degree. This proves strict
$\A^1$-invariance on the big Nisnevich site
$\mathscr C_k=\Sm/k$.
\end{proof}

\section{Perfection and cancellation theorem for framed motives}\label{sec:perfection-framed-motives}

We apply Theorem \ref{thm:main} and the explicit Frobenius computations
above to framed motives. Throughout this section $k$ is a field of
exponential characteristic $p$, $R=\mathbb Z[1/p]$, $K=k^{\perf}$, and
$\tau\colon\Spec K\to\Spec k$ is the canonical morphism. For a field $F$
and $Y,V\in\Sm/F$ we write
\[
 \mathcal C_F(Y)(V)
 =N\bigl(\ZF_F(\Delta_F^\bullet\times V,Y)\bigr)\otimes R,
\]
where $N$ denotes normalized chains and $\ZF_F$ denotes stabilized
linear framed correspondences over $F$, $M_{\mathrm{fr},F}(Y)$ denotes
the framed motive of $Y$, a presheaf of $S^1$-spectra on $\Sm/F$. For a
presheaf $\mathcal E$ on $\Sm/K$ we write $\tau_*\mathcal E$ for the
presheaf $V\mapsto\mathcal E(V_K)$ on $\Sm/k$. For $X\in\Sm/k$, base
change of framed correspondences, $c\mapsto c_K$, defines natural
morphisms of presheaves on $\Sm/k$
\begin{equation}\label{eq:local-base-change}
 b_X\colon M_{\mathrm{fr},k}(X)\longrightarrow
 \tau_*M_{\mathrm{fr},K}(X_K),
 \qquad
 b_X^{\mathrm{lin}}\colon\mathcal C_k(X)\longrightarrow
 \tau_*\mathcal C_K(X_K).
\end{equation}
A morphism of presheaves of $S^1$-spectra on $\Sm/k$ is a \emph{stable
local equivalence\/} if it induces isomorphisms between the Nisnevich
sheaves associated with the presheaves of stable homotopy groups.

We first choose a finite pointed simplicial model of the circle
equipped with a map of degree $p$. Throughout this construction,
$S^1=\Delta^1/\partial\Delta^1$.
Choose an integer $r\geq0$ sufficiently large that
$D=\sd^rS^1$ admits pointed simplicial maps
\[
 u:D\longrightarrow S^1,
 \qquad
 v:D\longrightarrow S^1
\]
of degrees $1$ and $p$, respectively (we may choose $u$ to be the iterated last vertex map). Let $S^1_a$ and $S^1_b$
be two copies of $S^1$, and form the pointed double mapping
cylinder
\[
 \mathcal S^1
 =
 \bigl(
 S^1_a\vee(D\wedge\Delta^1_+)\vee S^1_b
 \bigr)\big/\bigl(
 (d,0)\sim u(d),\ (d,1)\sim v(d)
 \bigr).
\]
Here the endpoint identifications use $u$ with values in
$S^1_a$ and $v$ with values in $S^1_b$.

The pointed simplicial set $\mathcal S^1$ is finite.
Write
\[
 \delta_p:S^1=S^1_a\longrightarrow\mathcal S^1,
 \qquad
 e:S^1=S^1_b\longrightarrow\mathcal S^1
\]
for the canonical inclusions. Since $u$ is a weak equivalence,
the homotopy pushout property of the double mapping cylinder
implies that $e$ is a weak equivalence. Moreover, $\delta_p\circ u\simeq e\circ v$.
So under the identification
$\mathcal S^1\simeq S^1$ supplied by $e$, the map $\delta_p$
has degree $p$. Fix these choices.

Put $\mathcal S^0=S^0,\mathcal S^n=(\mathcal S^1)^{\wedge n}$, $n\geq1$.
Define $\mathbb S[1/p]$ to be the $S^1$-spectrum
\[
 \mathbb S[1/p]
 =
 (S^0,\mathcal S^1,\mathcal S^2,\ldots)
\]
whose structure maps are $\mathcal S^n\wedge S^1\xrightarrow{\id\wedge\delta_p}\mathcal S^n\wedge\mathcal S^1=\mathcal S^{n+1}$.

\begin{definition}\label{frmotive1/p}
Let $Y\in\Sm/k$. The \emph{framed motive with
$1/p$-coefficients} of $Y$ is the $S^1$-spectrum
\[
 M_{\mathrm{fr}}(Y)[1/p]
 =
 \bigl(
 C_*\Fr(-,Y),
 C_*\Fr(-,Y\otimes\mathcal S^1),
 C_*\Fr(-,Y\otimes\mathcal S^2),
 \ldots
 \bigr).
\]
Its $m$-th structure map is the composite
\[
\begin{aligned}
 C_*\Fr(-,Y\otimes\mathcal S^m)\wedge S^1
 &\longrightarrow
 C_*\Fr\bigl(-,Y\otimes(\mathcal S^m\wedge S^1)\bigr)\\
 &\xrightarrow{(\id\wedge\delta_p)_*}
 C_*\Fr(-,Y\otimes\mathcal S^{m+1}),
\end{aligned}
\]
where the first arrow is the assembly map.
\end{definition}

The levelwise maps $\delta_p^{\wedge n}\colon S^n\to\mathcal S^n$ are
compatible with the structure maps, so they define a map of
$S^1$-spectra $\mathbb S\to\mathbb S[1/p]$ and hence a natural map
$M_{\mathrm{fr}}(Y)\to M_{\mathrm{fr}}(Y)[1/p]$. It induces, for every
integer $i$, a natural isomorphism of presheaves
\[
 \pi_i\bigl(M_{\mathrm{fr}}(Y)\bigr)\otimes_{\mathbb Z}\mathbb Z[1/p]
 \xrightarrow{\ \cong\ }
 \pi_i\bigl(M_{\mathrm{fr}}(Y)[1/p]\bigr).
\]
Indeed, fix $U\in\Sm/k$ and let $X$ be the $\Gamma$-space
$K\mapsto C_*\Fr(U,Y\otimes K)$, which is special by the additivity
theorem~\cite{GPfm}. Prolonged $\Gamma$-spaces preserve weak
equivalences, so $e^{\wedge n}$ induces weak equivalences
$X(S^n)\simeq X(\mathcal S^n)$. Under these identifications, the
structure maps of $M_{\mathrm{fr}}(Y)[1/p](U)$ become those of the Segal
spectrum of $X$ followed by $X$ applied to a self map of $S^{n+1}$ of
degree~$p$. Since $X$ is special, such a map induces multiplication
by $p$ on $\pi_*X(S^{n+1})$. Hence
$\pi_iM_{\mathrm{fr}}(Y)[1/p](U)$ is the colimit of
$\pi_iM_{\mathrm{fr}}(Y)(U)\xrightarrow{p}\pi_iM_{\mathrm{fr}}(Y)(U)
\xrightarrow{p}\cdots$, and the map induced by $\mathbb S\to\mathbb S[1/p]$
is the canonical map into this colimit. Thus $M_{\mathrm{fr}}(Y)[1/p]$
is obtained from $M_{\mathrm{fr}}(Y)$ by inverting $p$.

We are now in a position to prove the main result of this section.

\begin{theorem}[Perfection invariance]\label{thm:local-perfection}
Let $k$ be any field, $p$ its exponential characteristic and
$K=k^{\perf}$. For every $X\in\Sm/k$ the morphism
\[
 b_X[1/p]\colon M_{\mathrm{fr},k}(X)[1/p]\longrightarrow
 \tau_*M_{\mathrm{fr},K}(X_K)[1/p]
\]
is a stable local equivalence of presheaves of $S^1$-spectra on
$\Sm/k$.
\end{theorem}

\begin{proof}
If $k$ is perfect, then $K=k$ and $b_X$ is the identity. Assume
henceforth that $k$ is imperfect. Then $p>1$ and $k$ is infinite.
The theorem follows once $b_X^{\mathrm{lin}}$ is known to
be a Nisnevich local quasi-isomorphism. For a field $F$ and
$Y,U\in\Sm/F$, the linearization theorem \cite[Theorem 1.2]{GNP}, proven without any
hypothesis on the base field, gives natural identifications
\[
 H_i\bigl(M_{\mathrm{fr},F}(Y)(U);R\bigr)
 \cong H_i\bigl(\mathcal C_F(Y)(U)\bigr).
\]
Here $H_i(-;R)$ is ordinary homology of a spectrum, computed by
smashing with the Eilenberg--Mac Lane spectrum $HR$. Applying
this over $k$ with $U=V$ and over $K$ with $U=V_K$, we see that the maps
\[
 H_i\bigl(M_{\mathrm{fr},k}(X)(V);R\bigr)\longrightarrow
 H_i\bigl(M_{\mathrm{fr},K}(X_K)(V_K);R\bigr)
\]
induced by $b_X$ are identified, naturally in $V\in\Sm/k$, with the
maps $H_i(b_X^{\mathrm{lin}})$.

Let $Y\in\Sm/k$. The stabilization maps defining $\ZF_k$ are
postcomposition with $\sigma_Y$, which commutes with precomposition.
Hence $V\mapsto\mathcal C_k(Y)(V)$ is a complex of additive
$\ZF_*(k)$-presheaves, and its homology presheaves
$V\mapsto H_i(\mathcal C_k(Y)(V))$ are $\A^1$-invariant, as for any
Suslin complex. The presheaves $V\mapsto H_i(\mathcal C_k(Y)(V))$ are 
quasi-stable, and Theorem \ref{thm:main} makes their
Nisnevich sheafifications strictly $\A^1$-invariant. The
Nisnevich hypercohomology spectral sequence then shows that the
Nisnevich local replacement of $\mathcal C_k(Y)$ is
homotopy invariant, and hence motivically fibrant. Convergence here follows from the finite Nisnevich
cohomological dimension of each smooth scheme on which it is
evaluated.

It follows from~\cite[Corollary~16.6]{GP} that for a finite Zariski cover $Y=\bigcup_{a=1}^r U_a$, the augmented
alternating \v{C}ech complex of presheaves
\[
 \cdots\longrightarrow
 \bigoplus_{a<b}\ZF_k(-,U_a\cap U_b)\otimes R
 \longrightarrow\bigoplus_a\ZF_k(-,U_a)\otimes R
 \longrightarrow\ZF_k(-,Y)\otimes R
\]
is Nisnevich locally exact. 

Let $P$ be the total complex of the non-augmented part of this
\v{C}ech complex, and put $Q=\ZF_k(-,Y)\otimes R$. Local exactness says
that the augmentation $P\to Q$ is a Nisnevich local
quasi-isomorphism. The \v{C}ech complex of a finite cover is bounded,
and $C_*$ commutes with finite direct sums and total complexes. Therefore
$C_*Q=\mathcal C_k(Y)$, and $C_*P$ is the total complex on the
left hand side of \eqref{eq:perfection-target-descent} below. We do not
claim that $C_*$ preserves Nisnevich local quasi-isomorphisms. Instead
we use motivic equivalences of complexes of
presheaves of $R$-modules. The map $P\to Q$ is a motivic equivalence,
and so are $P\to C_*P$ and $Q\to C_*Q$, because the
map of a complex of presheaves into its Suslin complex is an
$\A^1$-equivalence. By two-out-of-three, $C_*P\to C_*Q$ is a motivic
equivalence. By the preceding paragraph, applied to $Y$ and to every
finite intersection of the $U_a$, the complex $C_*Q$ and all terms of
$C_*P$ are $\A^1$-local. The double complex spectral sequence implies 
$C_*P$ is $\A^1$-local as well (alternatively, use~\cite[Lemma~3.2]{Gar}). A motivic equivalence
between $\A^1$-local complexes is a Nisnevich local quasi-isomorphism,
because motivic equivalences between motivically fibrant complexes are
sectionwise quasi-isomorphisms. This proves the Nisnevich local
quasi-isomorphism
\begin{equation}\label{eq:perfection-target-descent}
 \operatorname{Tot}\!\left[
 \cdots\longrightarrow
 \bigoplus_{a<b}\mathcal C_k(U_a\cap U_b)
 \longrightarrow\bigoplus_a\mathcal C_k(U_a)
 \right]
 \xrightarrow{\simeq_{\Nis}}\mathcal C_k(Y).
\end{equation}

We claim that for every $Y\in\Sm/k$ postcomposition with
the relative Frobenius $F_Y:Y\to Y^{(p)}$ is a Nisnevich
local quasi-isomorphism
\begin{equation}\label{eq:perfection-target-frobenius}
 (F_Y)_*:\mathcal C_k(Y)
 \xrightarrow{\simeq_{\Nis}}\mathcal C_k(Y^{(p)}).
\end{equation}

We first prove the stronger assertion of an objectwise
quasi-isomorphism when $Y=U$ is smooth affine and admits
an étale morphism $U\to\A_k^d$.
For a smooth source $V$ and a smooth target $W$, write
\[
 C^{(n)}_\bullet(V,W)
 =
 C_\bullet\bigl(\ZF_n(-,W)\otimes_{\mathbb Z}R\bigr)(V).
\]
Thus its group in chain degree $q$ is
$C^{(n)}_q(V,W)
 =
 \ZF_n(V\times\Delta_k^q,W)\otimes_{\mathbb Z}R$,
and
\[
 \mathcal C_k(W)(V)
 =
 \varinjlim_n C^{(n)}_\bullet(V,W),
\]
where the transition maps add suspension coordinates.

The transfer constructed in the proof of Theorem~\ref{prop:Fiso} is a correspondence
$\tr_U\in\ZF_d(U^{(p)},U)$ satisfying the framed homotopies~\eqref{eq:composites}
\[
 F_U\circ\tr_U\simeq\Theta_{p,d,U^{(p)}},
 \qquad
 \tr_U\circ F_U\simeq\Theta_{p,d,U}.
\]
At each fixed framing level $n$, postcomposition defines
chain maps
\[
 (F_U)_{*,n}:
 C^{(n)}_\bullet(V,U)
 \longrightarrow C^{(n)}_\bullet(V,U^{(p)})
\]
and
\[
 (\tr_U)_{*,n}:
 C^{(n)}_\bullet(V,U^{(p)})
 \longrightarrow C^{(n+d)}_\bullet(V,U).
\]

The displayed framed homotopies induce chain homotopies
between the corresponding composite chain maps and
postcomposition with the respective $\Theta$-correspondences.
To see this, let $H\in\ZF_r(W\times\A^1,W')$ be a framed homotopy and let
$c\in\ZF_n(V\times\Delta^\bullet,W)$. The composite
\[
 H\circ(c\boxtimes\id_{\A^1})
 \in
 \ZF_{n+r}(V\times\Delta^\bullet\times\A^1,W')
\]
induces a simplicial homotopy between the endpoint maps on simplicial abelian groups. 
Passing to their associated chain complexes gives a chain homotopy.
The construction is additive and therefore also
applies with $R$-coefficients.

We next describe the operations on stabilized homology
that occur in the transfer identities.
Appending a coordinate with framing $-t$ induces an
operation $\nu$ on stabilized homology. To verify its
compatibility with suspension, compare the last two
framing entries $(s,-t)$ and $(-s,t)$. Their framing
matrices differ by $\diag(-1,-1)$, which is a product
of elementary matrices. Lemma~\ref{lem:elem} therefore
gives a framed homotopy between these operations, and
hence they induce the same map on homology. The same
argument identifies appending two negatively framed
coordinates with appending two suspension coordinates. So $\nu^2=1$ on stabilized homology.

By Lemma~\ref{lem:unit} adding one coordinate with
framing $t^p$ induces $p\nu$. Applying this successively to $d$ coordinates gives
\begin{equation}\label{eq:perfection-scalar}
 (\Theta_{p,d})_*=(p\nu)^d.
\end{equation}
Here the left hand side means the operation obtained by
appending the $d$ indicated coordinates. Since $p$ is invertible in $R$ and $\nu^2=1$, this operation
is an automorphism.

We now prove that $(F_U)_*$ induces an isomorphism on
homology for every smooth source $V$. Postcomposition
with $F_U$ commutes with suspension and therefore defines a chain map on the stabilized complexes.
For surjectivity, let $\eta\in H_q\bigl(\mathcal C_k(U^{(p)})(V)\bigr)$.
Put $\eta'=(p\nu)^{-d}\eta$, and represent $\eta'$ by
a cycle $z\in C^{(n)}_q(V,U^{(p)})$
at some finite level $n$. Then $\tr_U\circ z$
is a cycle of framing level $n+d$. The first transfer
homotopy gives
\[
 \begin{aligned}
 (F_U)_*[\tr_U\circ z]
 =[(F_U\circ\tr_U)\circ z]
 =[\Theta_{p,d,U^{(p)}}\circ z]
 =(p\nu)^d\eta'=\eta.
 \end{aligned}
\]
Indeed, the framed homotopy induces a chain homotopy $K$
between the two postcomposition maps. Since $\partial z=0$,
\[
 (F_U\circ\tr_U)\circ z-\Theta_{p,d,U^{(p)}}\circ z=\partial K_q(z)+K_{q-1}(\partial z) =\partial K_q(z).
\]
The two cycles therefore represent the same homology class.

For injectivity, let $\xi\in H_q\bigl(\mathcal C_k(U)(V)\bigr)$
satisfy $(F_U)_*\xi=0$. After increasing the framing
level, we may choose a cycle
$z\in C^{(n)}_q(V,U)$ representing $\xi$ and a chain
$b\in C^{(n)}_{q+1}(V,U^{(p)})$
such that $(F_U)_{*,n}z=\partial b$.
So,
\[
 (\tr_U\circ F_U)\circ z
 =
 \partial(\tr_U\circ b).
\]
The second transfer homotopy therefore implies that
$\Theta_{p,d,U}\circ z$ is also a boundary. Passing to
stabilized homology and using
\eqref{eq:perfection-scalar}, we obtain $(p\nu)^d\xi=0$. Since $(p\nu)^d$ is invertible, $\xi=0$.

The choices of finite level cycles and boundary
equalities used above are justified by exactness of
filtered colimits of $R$-modules. In particular,
\[
 H_q\!\left(\varinjlim_n C^{(n)}_\bullet(V,W)\right)
 \cong
 \varinjlim_n H_q\bigl(C^{(n)}_\bullet(V,W)\bigr).
\]
Every homology class is therefore represented by a cycle
at some finite level. Also, a class that vanishes after
stabilization becomes a boundary at a sufficiently
large finite level. This completes the proof of the
objectwise quasi-isomorphism for $U$ admitting étale
coordinates.

Cover $Y$ by finitely many affine opens $U_1,\ldots,U_r$ admitting
\'{e}tale coordinates. Since $Y$ is separated (objects of $\Sm/k$ are
separated), every finite intersection $U_I=\bigcap_{a\in I}U_a$ is
affine. Restricting the \'{e}tale coordinates of any $U_a$ with
$a\in I$ gives \'{e}tale coordinates on $U_I$, which are separated
because $U_I$ is affine. So the local calculation applies to every
$U_I$. The Frobenius twists $U_a^{(p)}$ cover $Y^{(p)}$ with
$(U_I)^{(p)}=\bigcap_{a\in I}U_a^{(p)}$ and relative Frobenius is
compatible with all inclusions in the two covers. Hence the maps
$(F_{U_I})_*$ assemble into a morphism from the total complex in
\eqref{eq:perfection-target-descent} for the cover $(U_a)$ to the one
for the cover $(U_a^{(p)})$ of $Y^{(p)}$. By the local calculation it is
a sectionwise quasi-isomorphism, and under the augmentations it is
compatible with $(F_Y)_*$. Since both augmentations are
local quasi-isomorphisms by
\eqref{eq:perfection-target-descent}, this proves
\eqref{eq:perfection-target-frobenius}.

Given $V\in\Sm/k$ let $F_V\colon V\to T(V)=V^{(p)}$ be the relative
Frobenius as in \S\ref{naturality}. We will use the following observations:

\begin{enumerate}
\item[(a)] The projection $T(V)=V\times_{k,F_k}k\to V$ is a base change
of the absolute Frobenius of $\Spec k$, hence an integral universal
homeomorphism. For $U\to V$ \'{e}tale one has $T(U)=U\times_VT(V)$.
By \cite[IV, 18.1.2]{EGA} or~\cite[Tag 04DZ]{Stacks} $U\mapsto T(U)$ is an equivalence between
the small \'{e}tale sites of $V$ and of $T(V)$. It preserves and
reflects Nisnevich coverings by the residue field argument in the proof
of Theorem~\ref{thm:main}. Consequently, for every presheaf
$\mathcal G$ on $\Sm/k$, the restriction of $\mathcal G\circ T$ to the
small Nisnevich site of $V$ is identified with the restriction of
$\mathcal G$ to that of $T(V)$. Therefore precomposition with $T$
preserves Nisnevich sheaves and Nisnevich local isomorphisms, and
$(\mathcal G\circ T)_{\Nis}\cong\mathcal G_{\Nis}\circ T$.

\item[(b)] Let $\mathcal M$ be an $\A^1$-invariant quasi-stable additive
$\ZF_*(k)$-presheaf of $R$-modules. Then
$F^*\colon\mathcal M\circ T\to\mathcal M$ is a
Nisnevich local isomorphism. Indeed, by
\cite[Theorem~1.1]{GP} $\mathcal A=\mathcal M_{\Nis}$ is an
$\A^1$-invariant quasi-stable additive framed Nisnevich sheaf of
$R$-modules, so $F^*\colon\mathcal A\circ T\to\mathcal A$ is an
isomorphism by Theorem~\ref{prop:Fiso}. Since $F^*$ is natural in
the presheaf, the square formed by $F^*_{\mathcal M}$,
$F^*_{\mathcal A}$ and the maps $\mathcal M\circ T\to\mathcal A\circ T$
and $\mathcal M\to\mathcal A$ induced by sheafification commutes. The
latter two maps are Nisnevich local isomorphisms, the first by (a),
hence so is $F^*_{\mathcal M}$.
\end{enumerate}

For $Y\in\Sm/k$, applying $T$ to framed correspondences and using the
canonical identification $T(\Delta^q_k\times V)\cong\Delta^q_k\times
T(V)$ defines a morphism of complexes of presheaves on $\Sm/k$
\[
 \beta_Y\colon\mathcal C_k(Y)\longrightarrow
 \mathcal C_k(Y^{(p)})\circ T,\qquad c\longmapsto T(c),
\]
at $V=\Spec k$ it is the chain map induced by base change along the
absolute Frobenius $F_k\colon k\to k$. We claim that $\beta_Y$ is a
Nisnevich local quasi-isomorphism.

Let $V\in\Sm/k$, put $W=\Delta^q_k\times V$, and let
$c\in\Fr_n(W,Y)$. The proof of \eqref{eq:exact} applies to any smooth
source and gives
\[
 (F_Y\circ c)^{[p]}=\mu_p^*\bigl(T(c)\circ F_W\bigr).
\]
Write $\gamma(c)=T(c)\circ F_W$. Applying the two homotopies of
Lemma~\ref{lem:power} gives, after stabilization, for every homology
class $z$ represented at level $n$,
\[
 (\Theta_{p,n})_*\bigl((F_Y)_*z-\gamma_*z\bigr)=0.
\]
The constructions commute with pullback in the source.
Consequently, the parameterized framed homotopies respect
the face and degeneracy maps and induce simplicial homotopies,
natural in the smooth source $V$. Passing to the associated
chain complexes gives homotopies natural in $V$.
By~\eqref{eq:perfection-scalar}, the operation
$(\Theta_{p,n})_*=(p\nu)^n$ is invertible on stabilized
homology: its inverse is $p^{-n}\nu^n$, since $\nu^2=1$.
We may therefore cancel $(\Theta_{p,n})_*$ in the
preceding equality. It follows that $(F_Y)_*$ and
$\gamma$ induce the same map of homology presheaves
on $\Sm/k$.

We next compare $\gamma$ with the composite of
$\beta_Y$ and pullback along Frobenius on the smooth
source. Write
\[
 \delta_Y:
 \mathcal C_k(Y^{(p)})\circ T
 \longrightarrow
 \mathcal C_k(Y^{(p)})
\]
for the morphism whose component at $V$ is pullback
along $F_V:V\to T(V)$. In degree $q$,
it is given by
\[
 (\delta_Y)_{V,q}(a)
 =
 (\id_{\Delta_k^q}\times F_V)^*a.
\]
Thus
\[
 (\delta_Y\circ\beta_Y)_{V,q}(c)
 =
 T(c)\circ(\id_{\Delta_k^q}\times F_V).
\]
In contrast, the auxiliary map used above is
\[
 \gamma_{V,q}(c)=T(c)\circ F_W,
 \qquad W=\Delta_k^q\times V.
\]

Under the canonical identification
$T(\Delta_k^q)\cong\Delta_k^q$, let
$F_{\Delta^q}$ denote the relative Frobenius of the
algebraic simplex. In barycentric coordinates it is
the $k$-morphism
\[
 F_{\Delta^q}(t_0,\ldots,t_q)
 =
 (t_0^p,\ldots,t_q^p).
\]
The relative Frobenius of $W$ factors as
\[
 F_W
 =
 F_{\Delta^q}\times F_V
 =
 \bigl(F_{\Delta^q}\times\id_{T(V)}\bigr)
 \circ
 \bigl(\id_{\Delta_k^q}\times F_V\bigr).
\]
Consequently, if $P_V$ denotes pullback along
$F_{\Delta^\bullet}\times\id_{T(V)}$ on
$\mathcal C_k(Y^{(p)})(T(V))$, then
\[
 \gamma_V=(\delta_Y)_V\circ P_V\circ(\beta_Y)_V.
\]

We show that $P_V$ is chain homotopic to the identity,
naturally in $V$. Consider the family
\begin{equation}\label{eq:perfection-simplex-homotopy}
 h_q:\Delta_k^q\times\A^1\longrightarrow\Delta_k^q,
 \qquad
 (t_0,\ldots,t_q;u)\longmapsto
 \bigl((1-u)t_0+ut_0^p,\ldots,
       (1-u)t_q+ut_q^p\bigr).
\end{equation}
This is well-defined because
\[
 \sum_{j=0}^q\bigl((1-u)t_j+ut_j^p\bigr)
 =
 (1-u)\sum_{j=0}^q t_j
 +u\left(\sum_{j=0}^q t_j\right)^p
 =
 1.
\]
Its endpoints are $ h_q(-,0)=\id_{\Delta_k^q}$, $h_q(-,1)=F_{\Delta^q}$.

These families respect the cosimplicial structure.
Indeed, a map $\theta:[m]\to[q]$ in $\Delta$
induces the morphism
\[
 \theta_*:\Delta_k^m\longrightarrow\Delta_k^q,
 \qquad
 (s_0,\ldots,s_m)\longmapsto
 \left(\sum_{\theta(j)=0}s_j,\ldots,
       \sum_{\theta(j)=q}s_j\right).
\]
Since taking $p$-th powers commutes with finite sums
in characteristic $p$, we have
\[
 h_q\circ(\theta_*\times\id_{\A^1})
 =
 \theta_*\circ h_m.
\]

To obtain a simplicial homotopy explicitly, let
$\epsilon:[q]\to[1]$ be order-preserving and put
\[
 u_\epsilon(t)=\sum_{\epsilon(j)=1}t_j,
 \qquad
 h_{q,\epsilon}(t)=h_q(t,u_\epsilon(t)).
\]
For a correspondence
$ a\in
 \ZF_n(\Delta_k^q\times T(V),Y^{(p)})
 \otimes_{\mathbb Z}R$, use the pullback $(h_{q,\epsilon}\times\id_{T(V)})^*a$.
The identities
\[
 u_\epsilon\circ\theta_*=u_{\epsilon\circ\theta},
 \qquad
 h_{q,\epsilon}\circ\theta_*
 =
 \theta_*\circ h_{m,\epsilon\circ\theta}
\]
show that these pullbacks define a simplicial homotopy.
For $\epsilon$ constantly $0$ the resulting map is the
identity, and for $\epsilon$ constantly $1$ it is $P_V$.
The construction preserves framing levels and commutes
with suspension, so it passes to the stabilized
simplicial abelian groups. It is also natural in $V$,
because it acts only on the simplex coordinate.
Passing to the associated chain complexes therefore
gives a chain homotopy between $P_V$ and the identity.

It follows that $\gamma$ and
$\delta_Y\circ\beta_Y$ induce the same maps on homology.
Combining this with the preceding comparison of
$\gamma$ and $(F_Y)_*$, we obtain, for every $i$,
\begin{equation}\label{eq:twist-factorization}
 H_i\bigl((F_Y)_*\bigr)
 =
 F^*\circ H_i(\beta_Y):
 H_i\bigl(\mathcal C_k(Y)\bigr)
 \longrightarrow
 H_i\bigl(\mathcal C_k(Y^{(p)})\bigr),
\end{equation}
where
\[
 F^*=H_i(\delta_Y):
 H_i\bigl(\mathcal C_k(Y^{(p)})\bigr)\circ T
 \longrightarrow
 H_i\bigl(\mathcal C_k(Y^{(p)})\bigr).
\]
Here we use the objectwise identification
$H_i\bigl(\mathcal C_k(Y^{(p)})\circ T\bigr)= H_i\bigl(\mathcal C_k(Y^{(p)})\bigr)\circ T$.

The presheaf $H_i(\mathcal C_k(Y^{(p)}))$ is $\A^1$-invariant,
additive and stable. Item~(b) above therefore implies that
$F^*$ is a Nisnevich local isomorphism.
We have already shown that $H_i((F_Y)_*)$ is also a
Nisnevich local isomorphism. After Nisnevich sheafification,
\eqref{eq:twist-factorization} expresses the latter
isomorphism as the composite of the former with the
sheafification of $H_i(\beta_Y)$. Hence
$H_i(\beta_Y)$ is a Nisnevich local isomorphism for
every $i$. This proves that $\beta_Y$ is a
Nisnevich local quasi-isomorphism.

Next, let us use $k_n$, $\sigma_n$ and $\rho_n$ as in \S\ref{descent}. For $V\in\Sm/k$,
the isomorphism $\sigma_n\colon k\xrightarrow{\cong}k_n$ together with
\eqref{eq:transport} give identifications, natural in $V$,
\[
 \mathcal C_{k_n}(X_{k_n})(V_{k_n})\cong\mathcal C_k(T^nX)(T^nV),
\]
under which base change along $k_n\subset k_{n+1}$ becomes
$(\beta_{T^nX})_{T^nV}$, because
\[
 \sigma_{n+1}^{-1}\circ(k_n\hookrightarrow k_{n+1})\circ\sigma_n=F_k.
\]
Explicit framed correspondences, common refinements of their
\'{e}tale neighbourhoods, and additivity relations descend through
filtered unions of fields by finite presentation (see the proof of Proposition~\ref{prop:descent}). So
$\mathcal C_K(X_K)(V_K)\cong
\varinjlim_n\mathcal C_{k_n}(X_{k_n})(V_{k_n})$, naturally in $V$, and
$b_X^{\mathrm{lin}}$ is identified with the canonical map from the
initial term to the colimit of the sequence of complexes of presheaves
on $\Sm/k$
\[
 \mathcal C_k(X)\xrightarrow{\ \beta_X\ }\mathcal C_k(TX)\circ T
 \xrightarrow{\ \beta_{TX}\circ T\ }\mathcal C_k(T^2X)\circ T^2
 \longrightarrow\cdots.
\]
Since homology is computed objectwise,
$H_i(\beta_{T^nX}\circ T^n)=H_i(\beta_{T^nX})\circ T^n$.
As we have shown above, the map \(H_i(\beta_{T^nX})\) is a Nisnevich local isomorphism. Iterating item~(a), we have
$(\mathcal G\circ T^n)_{\mathrm{Nis}}\cong\mathcal G_{\mathrm{Nis}}\circ T^n$
naturally in every presheaf \(\mathcal G\). Therefore
\[
\bigl(H_i(\beta_{T^nX}\circ T^n)\bigr)_{\mathrm{Nis}}
\cong
\bigl(H_i(\beta_{T^nX})\bigr)_{\mathrm{Nis}}\circ T^n
\]
is an isomorphism.
We see that the $n$-th transition $\beta_{T^nX}\circ T^n$ is a Nisnevich local quasi-isomorphism. 
As homology commutes with filtered colimits, and
Nisnevich sheafification is exact and commutes with filtered colimits,
$b_X^{\mathrm{lin}}$ is a Nisnevich local quasi-isomorphism. This
proves the theorem.
\end{proof}

Next, we will record how Theorem~\ref{thm:local-perfection} reads after fibrant
replacement. On $\Sm/k$ and on $\Sm/K$ we use the (injective) stable
Nisnevich local model structure on presheaves of $S^1$-spectra
\cite{Jar00}, whose weak equivalences are the stable local equivalences,
and we write $E\mapsto E_f$ for a functorial fibrant replacement. Thus
$E\to E_f$ is a stable local equivalence and $E_f$ is fibrant. A
stable local equivalence between fibrant objects is a sectionwise
stable equivalence. A presheaf of $S^1$-spectra $E$ is \emph{Nisnevich
excisive\/} if $E(\emptyset)$ is contractible and $E$ carries every
elementary distinguished square to a homotopy cartesian square of
spectra. Fibrant objects are Nisnevich excisive. Conversely, if $E$ is
Nisnevich excisive, then $E\to E_f$ is a sectionwise stable
equivalence, by the Brown--Gersten criterion \cite[\S3.1]{MV} applied
levelwise to a sectionwise $\Omega$-spectrum model of $E$.

We should explain the subscripts. The target of $b_X$ is $\tau_*$ of a
presheaf on $\Sm/K$: its value at $Y\in\Sm/k$ is formed from framed
correspondences over $K$ with source $Y_K$. It is not a framed motive
over $k$ of the $K$-scheme $X_K$. Indeed, if $k$ is imperfect, a
nonempty smooth $k$-scheme admits no $k$-morphism to $\Spec K$: its
function fields are separable over $k$ and therefore contain no element
of $K\smallsetminus k$. Hence framed correspondences over $k$ from
objects of $\Sm/k$ to $X_K$ are trivial. Accordingly, the fibrant
replacement of $M_{\mathrm{fr},K}(X_K)[1/p]$ is formed on $\Sm/K$, and
the next lemma compares it with fibrant replacement on $\Sm/k$.

\begin{lemma}\label{lem:pi-fibrant}
Let $k$ be a field, let $K=k^{\perf}$, and let $G$ be a presheaf of
$S^1$-spectra on $\Sm/K$.
\begin{enumerate}
\item[(i)] $\tau_*$ maps stable local equivalences of presheaves of
$S^1$-spectra on $\Sm/K$ to stable local equivalences on $\Sm/k$.
\item[(ii)] If $G$ is Nisnevich excisive on $\Sm/K$, then $\tau_*G$ is
Nisnevich excisive on $\Sm/k$.
\item[(iii)] For every $Y\in\Sm/k$ there is a zigzag of stable
equivalences, natural in $G$ and $Y$,
\[
 (\tau_*G)_f(Y)\simeq G_f(Y_K),
\]
where the fibrant replacement on the left is formed on $\Sm/k$ and that
on the right on $\Sm/K$.
\end{enumerate}
\end{lemma}

\begin{proof}
(i) Since $K/k$ is algebraic and purely inseparable, the
morphism $\Spec K\to\Spec k$ is integral, surjective
and radicial, hence a universal homeomorphism
by~\cite[Tags 01S4 and 04DF]{Stacks}.
Its base change $Y_K\to Y$ is therefore an integral
universal homeomorphism
by~\cite[Tags 01WL and 0CEU]{Stacks}.
For every étale morphism $U\to Y$ one has a canonical isomorphism
\[
 U\times_Y Y_K
 \cong U\times_k\Spec K
 =U_K.
\]
By~\cite[IV, 18.1.2]{EGA}, $U\mapsto U_K$ is an equivalence between the
small \'{e}tale sites of $Y$ and of $Y_K$. Let $y\in Y$, and let $y_K$
be the unique point of $Y_K$ over $y$. If $u\in U$ lies over $y$ and
$u_K$ is the point of $U_K$ over $u$, then $\kappa(u)/\kappa(y)$ is
finite separable and $\kappa(y_K)/\kappa(y)$ is purely inseparable, so
$\kappa(u_K)=\kappa(u)\otimes_{\kappa(y)}\kappa(y_K)$ and
$[\kappa(u_K):\kappa(y_K)]=[\kappa(u):\kappa(y)]$. Hence
$(U,u)\mapsto(U_K,u_K)$ is an equivalence between the categories of
Nisnevich neighbourhoods of $y$ in $Y$ and of $y_K$ in $Y_K$ (see the proof of Lemma~\ref{polezno} as well). So the
stalk of $\tau_*G$ at $y$ is the stalk of $G$ at $y_K$. Since stable
homotopy groups commute with filtered colimits and stable local
equivalences are detected on stalks, (i) follows.

(ii) Let $U\subset Y$ and $V\to Y$ form an elementary distinguished
square in $\Sm/k$. Then $U_K\subset Y_K$ is open, $V_K\to Y_K$ is
\'{e}tale, $V_K\times_{Y_K}U_K=(V\times_YU)_K$, and the isomorphism
$\bigl(V\times_Y(Y\smallsetminus U)\bigr)_{\mathrm{red}}\cong
(Y\smallsetminus U)_{\mathrm{red}}$ remains an isomorphism after base
change to $K$ and passage to reduced subschemes. So the base changed
square is elementary distinguished in $\Sm/K$. The square obtained by
applying $\tau_*G$ to the original square is the square obtained by
applying $G$ to the base-changed one, hence homotopy cartesian.
Moreover $(\tau_*G)(\emptyset)=G(\emptyset)$ is contractible.

(iii) By (i) $\tau_*G\to\tau_*(G_f)$ is a stable local equivalence on
$\Sm/k$. Since $G_f$ is fibrant, hence Nisnevich excisive, $\tau_*(G_f)$
is Nisnevich excisive by (ii). Consider
\[
 (\tau_*G)_f\longrightarrow\bigl(\tau_*(G_f)\bigr)_f\longleftarrow
 \tau_*(G_f).
\]
The left map is obtained by applying $(-)_f$ to
$\tau_*G\to\tau_*(G_f)$, so it is a stable local equivalence between
fibrant objects. So it is a sectionwise stable equivalence. The right map
is a sectionwise stable equivalence because $\tau_*(G_f)$ is Nisnevich
excisive. Since $\tau_*(G_f)(Y)=G_f(Y_K)$ the assertion follows.
\end{proof}

\begin{corollary}\label{cor:local-perfection-fibrant}
Let $k$ be any field, $p$ its exponential characteristic and
$K=k^{\perf}$. For all $X,Y\in\Sm/k$, the morphism $b_X$ induces a
zigzag of stable equivalences of $S^1$-spectra, natural in $Y$,
\[
 M_{\mathrm{fr},k}(X)[1/p]_f(Y)\;\simeq\;
 M_{\mathrm{fr},K}(X_K)[1/p]_f(Y_K),
\]
where the fibrant replacement on the left is formed on $\Sm/k$ and that
on the right on $\Sm/K$.
\end{corollary}

\begin{proof}
By Theorem~\ref{thm:local-perfection}, $b_X[1/p]$ is a stable local
equivalence. Applying $(-)_f$ on $\Sm/k$ gives a stable local
equivalence between fibrant objects
\[
 M_{\mathrm{fr},k}(X)[1/p]_f\longrightarrow
 \bigl(\tau_*M_{\mathrm{fr},K}(X_K)[1/p]\bigr)_f,
\]
which is therefore a sectionwise stable equivalence. Part (iii) of
Lemma~\ref{lem:pi-fibrant}, applied to $G=M_{\mathrm{fr},K}(X_K)[1/p]$,
identifies the value of the target at $Y$ with
$M_{\mathrm{fr},K}(X_K)[1/p]_f(Y_K)$.
\end{proof}

Finally, we use Theorem~\ref{thm:local-perfection}, through
Lemma~\ref{lem:pi-fibrant} and Corollary~\ref{cor:local-perfection-fibrant},
to extend the cancellation theorem \cite[Theorem A]{AGP} to imperfect
base fields after inverting $p$. We use the notation of
\cite[\S\S1 and 3]{AGP}.
Let $\mathbb G$ be the cone of $1\colon\mathrm{pt}_+\to(\mathbb G_m)_+$ in
pointed simplicial presheaves, and write
$\underline{\mathrm{Hom}}(\mathbb G,-)$ for the internal Hom of presheaves
of $S^1$-spectra. For sectionwise fibrant $A$, the spectrum
$\underline{\mathrm{Hom}}(\mathbb G,A)(U)$ is naturally equivalent to the
homotopy fibre of the restriction $A(U\times\mathbb G_m)\to A(U)$ along
$U\times\{1\}\hookrightarrow U\times\mathbb G_m$. For $X\in\Sm/k$ and
$n\geq0$, the twisted framed motive
$M_{\mathrm{fr},k}(X)(n)=M_{\mathrm{fr},k}(X\times\mathbb G_m^{\wedge n})$
is the diagonal of the multisimplicial presheaf of $S^1$-spectra obtained
by applying $M_{\mathrm{fr},k}(X\times-)$ to the simplicial object
$\mathbb G_m^{\wedge n}$ of $\Fr_0(k)$. Let
\[
 a_n\colon M_{\mathrm{fr},k}(X)(n)\longrightarrow
 \underline{\mathrm{Hom}}\bigl(\mathbb G,M_{\mathrm{fr},k}(X)(n+1)\bigr)
\]
be the structure maps of \cite[\S3]{AGP}. By definition, $a_n$ is the composite of the
external product with $\mathbb G_m^{\wedge1}$
\[
 M_{\mathrm{fr},k}(X)(n)\longrightarrow
 \underline{\mathrm{Hom}}\bigl(\mathbb G_m^{\wedge1},
 M_{\mathrm{fr},k}(X)(n+1)\bigr)
\]
with the map induced by the adjunction unit
$\mathbb G\to\mathbb G_m^{\wedge1}|_{\Sm/k}$. The map $a_n[1/p]$ 
is defined in a similar fashion.
We may and do assume that the fibrant replacement maps
$\alpha_A\colon A\to A_f$ are trivial cofibrations in the injective model
structure. Since $\underline{\mathrm{Hom}}(\mathbb G,-)$ preserves fibrant
objects, every map $\varphi\colon A\to\underline{\mathrm{Hom}}(\mathbb G,A')$
has an extension
\[
 \varphi^\sharp\colon A_f\longrightarrow
 \underline{\mathrm{Hom}}(\mathbb G,A'_f),\qquad
 \varphi^\sharp\circ\alpha_A
 =\underline{\mathrm{Hom}}(\mathbb G,\alpha_{A'})\circ\varphi.
\]
The map $a_n^\sharp$ is the map of~\cite[Theorem A(2)]{AGP}.

\begin{theorem}[Cancellation after inverting $p$]\label{cor:cancellation}
Let $k$ be an infinite field of exponential characteristic $p$. For
every $X\in\Sm/k$ and every $n\geq0$, the map
\[
 \bigl(a_n[1/p]\bigr)^\sharp\colon M_{\mathrm{fr},k}(X)(n)[1/p]_f
 \longrightarrow
 \underline{\mathrm{Hom}}\bigl(\mathbb G,
 M_{\mathrm{fr},k}(X)(n+1)[1/p]_f\bigr)
\]
is a schemewise stable equivalence and the map
\[
 a_n[1/p]\colon M_{\mathrm{fr},k}(X)(n)[1/p]
 \longrightarrow
 \underline{\mathrm{Hom}}\bigl(\mathbb G,
 M_{\mathrm{fr},k}(X)(n+1)[1/p]\bigr)
\]
is a local stable equivalence. For infinite perfect $k$ this is
\cite[Theorem A]{AGP} with $p$ inverted. 
\end{theorem}

\begin{proof}
Since $k$ is infinite, the field
$K=k^{\perf}$ is infinite perfect. Put
$M_n=M_{\mathrm{fr},k}(X)(n)[1/p]$, a presheaf on $\Sm/k$, and
$N_n=M_{\mathrm{fr},K}(X_K)(n)[1/p]$, a presheaf on $\Sm/K$. Superscripts
$k$ and $K$ indicate the base field of the structure maps.
Base change defines morphisms $b_{X,n}\colon M_n\to\tau_*N_n$.
Framed motives are functorial in level zero framed correspondences, and
base change commutes with them. Hence $b_{X,n}$ is the diagonal of the
maps $b_{X\times(\mathbb G_m^{\wedge n})_q}[1/p]$ of
\eqref{eq:local-base-change}, where $(\mathbb G_m^{\wedge n})_q$ runs
through the terms of $\mathbb G_m^{\wedge n}$ and
$\bigl(X\times(\mathbb G_m^{\wedge n})_q\bigr)_K=
X_K\times(\mathbb G_{m,K}^{\wedge n})_q$. Each of these maps is a stable
local equivalence by Theorem~\ref{thm:local-perfection}, and hence so is $b_{X,n}$.

External products with the simplicial object $\mathbb G_m^{\wedge1}$ of $\Fr_0(k)$ commute
with base change, and the adjunction unit
$\mathbb G\to\mathbb G_m^{\wedge1}|_{\Sm/k}$ is defined over the prime
field. Moreover $\tau_*$ commutes with
$\underline{\mathrm{Hom}}(\mathbb G_m^{\wedge1},-)$ and with
$\underline{\mathrm{Hom}}(\mathbb G,-)$, because
$(U\times\mathbb G_m)_K=U_K\times\mathbb G_{m,K}$. Hence
\begin{equation}\label{eq:cancellation-square}
 \underline{\mathrm{Hom}}(\mathbb G,b_{X,n+1})\circ a_n^k[1/p]
 =\tau_*\bigl(a_n^K[1/p]\bigr)\circ b_{X,n}.
\end{equation}
Since $K$ is an infinite perfect
field, \cite[Theorem A(2)]{AGP} says that $(a_n^K)^\sharp$ is a
schemewise stable equivalence (no restriction on the characteristic is
needed). As inverting $p$ preserves schemewise stable equivalences and Nisnevich excisiveness,
$(a_n^K[1/p])^\sharp\colon N_{n,f}\to\underline{\mathrm{Hom}}(\mathbb G,N_{n+1,f})$ is a schemewise stable equivalence.

Let $\varphi\colon A\to\underline{\mathrm{Hom}}(\mathbb G,A')$ and
$\psi\colon B\to\underline{\mathrm{Hom}}(\mathbb G,B')$ be maps of
presheaves of $S^1$-spectra on $\Sm/k$ or on $\Sm/K$, and let
$u\colon A\to B$ and $u'\colon A'\to B'$ be stable local equivalences with
$\underline{\mathrm{Hom}}(\mathbb G,u')\circ\varphi=\psi\circ u$. Assume
that $B$ and $B'$ are sectionwise fibrant and Nisnevich excisive. Then
$\varphi^\sharp$ is a schemewise stable equivalence if and only if $\psi$
is. Indeed, $\underline{\mathrm{Hom}}(\mathbb G,u'_f)\circ\varphi^\sharp$
and $\psi^\sharp\circ u_f$ agree after composition with $\alpha_A$, and
their common target $\underline{\mathrm{Hom}}(\mathbb G,B'_f)$ is fibrant,
so they are homotopic. The maps $u_f$ and $u'_f$ are stable local
equivalences between fibrant objects, hence sectionwise stable
equivalences, and $\underline{\mathrm{Hom}}(\mathbb G,-)$ preserves
sectionwise stable equivalences between sectionwise fibrant presheaves.
Finally,
$\psi^\sharp\circ\alpha_B=\underline{\mathrm{Hom}}(\mathbb G,\alpha_{B'})
\circ\psi$, and $\alpha_B$, $\alpha_{B'}$ are sectionwise stable
equivalences, because $B$ and $B'$ are Nisnevich excisive.

Apply the previous paragraph over $k$ with
$\varphi=a_n^k[1/p]$, $B=\tau_*(N_{n,f})$, $B'=\tau_*(N_{n+1,f})$,
$\psi=\tau_*\bigl((a_n^K[1/p])^\sharp\bigr)$,
$u=\tau_*(\alpha_{N_n})\circ b_{X,n}$ and
$u'=\tau_*(\alpha_{N_{n+1}})\circ b_{X,n+1}$. The maps $u$ and $u'$ are
stable local equivalences by the first paragraph and Lemma~\ref{lem:pi-fibrant}(i).
The presheaves $B$ and $B'$ are sectionwise fibrant and Nisnevich excisive by
Lemma~\ref{lem:pi-fibrant}(ii). By
\eqref{eq:cancellation-square} and the defining identity of
$(a_n^K[1/p])^\sharp$
\[
 \underline{\mathrm{Hom}}(\mathbb G,u')\circ\varphi
 =\tau_*\bigl(\underline{\mathrm{Hom}}(\mathbb G,\alpha_{N_{n+1}})\circ
 a_n^K[1/p]\bigr)\circ b_{X,n}
 =\psi\circ u.
\]
Finally, the value of $\psi$ at $Y\in\Sm/k$ is the value of
$(a_n^K[1/p])^\sharp$ at $Y_K$. As we have already shown
above, $(a_n^K[1/p])^\sharp$ is a schemewise stable equivalence, and hence so is $(a_n[1/p])^\sharp$.

The fact that $a_n[1/p]$ is a stable local equivalence follows from the first assertion and the fact that
$\underline{\mathrm{Hom}}(\mathbb G,\alpha_{M_n})$ is a stable local equivalence for all $n$.
The proof of the latter fact literally repeats the arguments of~\cite[pp.15-16]{AGP} if we also use Theorem~\ref{thm:main}. 
\end{proof}

Consider the bispectrum consisting of the twisted framed motives with $1/p$-coefficients
   $$M_{\mathrm{fr},k}^{\mathbb G}(X)[1/p]=(M_{\mathrm{fr},k}(X)[1/p],M_{\mathrm{fr},k}(X)(1)[1/p],M_{\mathrm{fr},k}(X)(2)[1/p],\ldots)$$
Each structure map is given by $a_n[1/p]$.

\begin{corollary}\label{bispectrum}
Let $k$ be an infinite field of exponential characteristic $p$, $X\in \Sm/k$.
Then the bispectrum
   $$M_{\mathrm{fr},k}^{\mathbb G}(X)[1/p]_f=(M_{\mathrm{fr},k}(X)[1/p]_f,M_{\mathrm{fr},k}(X)(1)[1/p]_f,M_{\mathrm{fr},k}(X)(2)[1/p]_f,\ldots)$$
obtained from $M_{\mathrm{fr},k}^{\mathbb G}(X)[1/p]$ by taking levelwise
Nisnevich local stable fibrant replacements with structure maps given by
$(a_n[1/p])^\sharp$ is a motivically fibrant $(S^1,\mathbb G)$-bispectrum. For infinite perfect $k$ this is
\cite[Theorem~B]{AGP} with $p$ inverted. 
\end{corollary}

\begin{proof}
The proof is like that of \cite[Theorem~B]{AGP} (if we also use Theorem~\ref{thm:main} here).
\end{proof}

\section{Framed motives of motivic Thom spaces}\label{sec:pairs}

We extend Theorem~\ref{thm:local-perfection} to the framed motives of the
sheaves $X/(X-Z)$, also referred to as motivic Thom spaces, keeping the notation of
Section~\ref{sec:perfection-framed-motives}. Let $F$ be a field, $Y\in\Sm/F$
and $S\subset Y$ a closed subset. Following \cite{GPfm,GNP},
an explicit framed correspondence of level $m$ from $U\in\Sm/F$ to
$Y/(Y-S)$ is a quadruple $(Z,W,\varphi,g)$, where $Z\subset\A^m_U$ is a
closed subset finite over $U$, $W$ is an \'{e}tale neighbourhood of $Z$
in $\A^m_U$, $\varphi=(\varphi_1,\ldots,\varphi_m)$ are regular functions
on $W$, and $g\colon W\to Y$ is a morphism such that $Z=Z(\varphi)\cap g^{-1}(S)$.
The equivalence of explicit correspondences, the pointed sets
$\Fr_m(U,Y/(Y-S))$, the groups $\ZF_m(U,Y/(Y-S))$, the suspension
$\Sigma=-\boxtimes\sigma$ and the groups $\ZF(U,Y/(Y-S))$ are
defined as for $S=Y$, in which case one recovers $\Fr_m(U,Y)$. 
Precomposition with framed correspondences, with the
framing of the source first, makes these framed presheaves, and there is
an external product with $\Fr_*(\mathrm{pt},\mathrm{pt})$. 

For a finite pointed set $K$ with
$K'=K\smallsetminus\{*\}$ put $(Y/(Y-S))\otimes K=(Y\times K')/((Y-S)\times K')$. 
The framed motive $M_{\mathrm{fr},F}(Y/(Y-S))$ is the
Segal $S^1$-spectrum of the $\Gamma$-space
$K\mapsto C_*\Fr\bigl(-,(Y/(Y-S))\otimes K\bigr)$. For
$(Y,S)=(X\times\A^n,X\times\{0\})$ it is the framed motive
$M_{\mathrm{fr}}(X\times T^n)$ of the relative motivic sphere
$X_+\wedge T^n$. We write $\mathcal C_F(Y/(Y-S))(V)=
N\bigl(\ZF_F(\Delta^\bullet_F\times V,Y/(Y-S))\bigr)\otimes R$. For
$Y\in\Sm/k$ let $S_K\subset Y_K$ be the preimage of $S$. Base change of
framed correspondences defines morphisms of presheaves on $\Sm/k$
\begin{gather*}
 b_{Y/(Y-S)}\colon M_{\mathrm{fr},k}(Y/(Y-S))\longrightarrow
 \tau_*M_{\mathrm{fr},K}(Y_K/(Y_K-S_K)),\\
 b^{\mathrm{lin}}_{Y/(Y-S)}\colon\mathcal C_k(Y/(Y-S))\longrightarrow
 \tau_*\mathcal C_K(Y_K/(Y_K-S_K)).
\end{gather*}
The framed motive with $1/p$-coefficients $M_{\mathrm{fr},k}(Y/(Y-S))[1/p]$
is defined similarly to~Definition~\ref{frmotive1/p}.

For an explicit framed correspondence
$b=(Z_b,W_b,\psi,h)$ in $\Fr_d(Y',Y)$, the \emph{lifted support} means the
image of the chosen morphism $s_b:Z_b\to W_b$ in the
\'{e}tale neighbourhood data, whose composite with
$W_b\to\A^d_{Y'}$ is the inclusion $Z_b\hookrightarrow\A^d_{Y'}$.

\begin{lemma}\label{lem:pair-compatible}
Let $Y',Y\in\Sm/F$ with closed subsets $S'\subset Y'$ and $S\subset Y$, and let
$b=(Z_b,W_b,\psi,h)\in\Fr_d(Y',Y)$. Write $q_b:W_b\longrightarrow\A^d_{Y'}\longrightarrow Y'$
for the source projection, and let $\widetilde Z_b\subset W_b$ denote the lifted support.
Suppose that $\widetilde Z_b\cap h^{-1}(S)\subset q_b^{-1}(S')$ as subsets of $W_b$.
For $c=(Z_c,W_c,\varphi,g)\in\Fr_n(U,Y'/(Y'-S'))$,
form $W=W_c\times_{Y'}W_b$ using $g$ and $q_b$, and denote its projections by
$p_c:W\to W_c$ and $p_b:W\to W_b$. Equip $W$ with its natural \'etale map to $\A^{n+d}_U$, the framing
\[
 \Phi=
 \bigl(
 \varphi_1\circ p_c,\ldots,\varphi_n\circ p_c,
 \psi_1\circ p_b,\ldots,\psi_d\circ p_b
 \bigr),
\]
and the target map $G=h\circ p_b$. Then the closed subset $D=Z(\Phi)\cap G^{-1}(S)$
is finite over $U$ and is identified by the ambient map $e:W=W_c\times_{Y'}W_b\rightarrow\A^{n+d}_U$
with a closed subset $Z_{b\circ c}\subset\A^{n+d}_U$. After restricting $W$, if necessary, to an \'etale
neighbourhood $W^\circ$ of this support, the data $b\circ c=\bigl(Z_{b\circ c},W^\circ,\Phi|_{W^\circ},G|_{W^\circ}\bigr)$
define an element of $\Fr_{n+d}(U,Y/(Y-S))$.

If $H\in\Fr_d(Y'\times\A^1,Y)$ satisfies the same
support condition with respect to $S'\times\A^1$ and $S$, then
\[
 H\circ(c\boxtimes\id_{\A^1})
 \in
 \Fr_{n+d}(U\times\A^1,Y/(Y-S))
\]
is a framed homotopy from $H_0\circ c$ to $H_1\circ c$.
\end{lemma}

\begin{proof}
Let $e_c:W_c\to\A^n_U$ and
$e_b:W_b\to\A^d_{Y'}$ be the neighbourhood maps.
The natural ambient map
$e:W=W_c\times_{Y'}W_b\to\A^{n+d}_U$ factors as
\[
 W_c\times_{Y'}W_b
 \xrightarrow{\ \id_{W_c}\times_{Y'}e_b\ }
 W_c\times_{Y'}\A^d_{Y'}
 \cong W_c\times_F\A^d_F
 \xrightarrow{\ e_c\times_F\id_{\A^d_F}\ }
 \A^n_U\times_F\A^d_F
 \cong\A^{n+d}_U.
\]
The first arrow is the base change of $e_b$ along
$W_c\times_{Y'}\A^d_{Y'}\rightarrow\A^d_{Y'}$,
and the second arrow, after the indicated identification,
is the base change of $e_c$ along the projection
$\A^n_U\times_F\A^d_F\rightarrow\A^n_U$.
Both arrows are therefore \'etale, so their composite $e$ is \'etale.

Let $\widetilde Z_c\subset W_c$ and
$\widetilde Z_b\subset W_b$ be the lifted supports,
equipped with their reduced induced scheme structures.
Consider the closed subset $D=Z(\Phi)\cap G^{-1}(S)\subset W$.
For $x\in D$ put $w=p_c(x)$ and $z=p_b(x)$. By definition,
\[
 \varphi(w)=0,\qquad \psi(z)=0,\qquad h(z)\in S.
\]
Since $b$ is an ordinary framed correspondence,
$\psi(z)=0$ implies $z\in\widetilde Z_b$.
The hypothesis therefore gives $q_b(z)\in S'$.
The fibre product identity
$g\circ p_c=q_b\circ p_b$ now implies $g(w)\in S'$,
and hence
\[
 w\in Z(\varphi)\cap g^{-1}(S')=\widetilde Z_c.
\]
Consequently, as closed subsets of $W$,
$D\subset C:=\widetilde Z_c\times_{Y'}\widetilde Z_b$,
where the map $\widetilde Z_c\to Y'$ is the restriction of $g$.

The projection $C\to\widetilde Z_c$ is finite,
being a base change of $Z_b\to Y'$.
Since $\widetilde Z_c\cong Z_c$ is finite over $U$,
the scheme $C$ is finite over $U$.
Equip $D$ with its reduced induced scheme structure.
The preceding containment makes $D$ a closed
subscheme of $C$, so $D$ is finite over $U$.

We next check the embedding into the ambient space.
Using the identifications of the lifted supports
with $Z_c$ and $Z_b$, the restriction of $e$ to $C$
is the composite of closed immersions
\[
 C\cong Z_c\times_{Y'}Z_b
 \lhook\joinrel\longrightarrow
 Z_c\times_{Y'}\A^d_{Y'}
 \cong Z_c\times_F\A^d_F
 \lhook\joinrel\longrightarrow
 \A^n_U\times_F\A^d_F
 \cong\A^{n+d}_U.
\]
Thus $e$ identifies $D$ with a reduced closed
subscheme $Z_{b\circ c}\subset\A^{n+d}_U$ finite over $U$.

To obtain an \'etale neighbourhood of precisely this
support, put $E=W\times_{\A^{n+d}_U}Z_{b\circ c}$.
The identification $D\cong Z_{b\circ c}$ defines a
section of the \'etale projection
$E\to Z_{b\circ c}$.
Its image is open in $E$, since a section of an
\'etale morphism is an open immersion by~\cite[Tag~024T]{Stacks}.
It is also closed in $E$, because $D$ is closed in $W$.
Therefore $E\setminus D$ is closed in $E$, and
hence in $W$. Set $W^\circ=W\setminus(E\setminus D)$.
Then $W^\circ$ contains $D$ and
$W^\circ\times_{\A^{n+d}_U}Z_{b\circ c}\cong Z_{b\circ c}$.
Thus $W^\circ$ is an \'etale neighbourhood of
$Z_{b\circ c}$. On this neighbourhood,
\[
 Z(\Phi|_{W^\circ})
 \cap(G|_{W^\circ})^{-1}(S)=D.
\]
The resulting quadruple therefore defines an element
of $\Fr_{n+d}(U,Y/(Y-S))$. Replacing either representative by an equivalent one
gives a common \'etale neighbourhood after taking the corresponding fibre products. On that
neighbourhood the framing functions and target maps agree, so the composites are equivalent.
Further restrictions of $W^\circ$ likewise do not change the represented correspondence.

A disjoint decomposition of either support induces a disjoint decomposition of $D$ by inverse image
under the corresponding projection. Hence the construction respects the additivity relations.
Let $f:U'\to U$ be a morphism in $\Sm/F$. There is a canonical isomorphism
\[
 (W_c\times_{Y'}W_b)\times_U U'
 \cong
 (W_c\times_U U')\times_{Y'}W_b.
\]
Under this identification, the framing functions and
target maps defining $f^*(b\circ c)$ and
$b\circ(f^*c)$ agree. Their supports also agree:
both are the inverse image of
$D=Z(\Phi)\cap G^{-1}(S)$ under the projection
$W\times_U U'\to W$, since inverse images commute
with intersections and zero loci.
The chosen restriction $W^\circ\subset W$ pulls back
to an \'etale neighbourhood of this support.
Consequently,
\[
 f^*(b\circ c)=b\circ(f^*c)
 \quad\text{in}\quad
 \Fr_{n+d}(U',Y/(Y-S)).
\]

Finally, apply the same construction to $H$ and
$c\boxtimes\id_{\A^1}$ over $U\times\A^1$.
The assumed support condition makes the resulting
support finite over $U\times\A^1$.
Compatibility with pullback identifies its fibres
at $0$ and $1$ with $H_0\circ c$ and $H_1\circ c$,
respectively. It is therefore the asserted framed
homotopy.
\end{proof}

\begin{remark}
For $d>0$, the lemma does not assert strict compatibility
of composition with suspension. Indeed, if $c$ has framing
$\varphi$ and $b$ has framing $\psi$, then
$b\circ\Sigma(c)$ and $\Sigma(b\circ c)$ have framing
orders $(\varphi,t,\psi)$ and $(\varphi,\psi,t)$, respectively.
In the transfer arguments below, composition is applied
to chosen cycles and boundaries at finite framing levels.
The resulting identities are then passed to stabilized
homology. These arguments therefore do not require
composition with $b$ itself to define a map on the
stabilization colimits.
\end{remark}

Write $q_b:W_b\to Y'$ for the source projection and
$\widetilde Z_b\subset W_b$ for the lifted support of $b$.
Suppose that a morphism $f:Y'\to Y$ satisfies
\[
 h|_{\widetilde Z_b}
 =
 f\circ q_b|_{\widetilde Z_b},
 \qquad
 f^{-1}(S)\subset S'.
 \]
Then the condition of Lemma~\ref{lem:pair-compatible} holds:
if $z\in\widetilde Z_b$ and $h(z)\in S$, then $q_b(z)\in f^{-1}(S)\subset S'$.

For a pair $(Y,S)$ with $Y\in\Sm/k$, let $\pi_Y:Y^{(p)}=Y\times_{\Spec k,F_{\Spec k}}\Spec k\rightarrow Y$
be the projection. This morphism is a universal homeomorphism, being a base change of the absolute
Frobenius of $\Spec k$~\cite[Tag 0CC8]{Stacks}. Define the closed subset $S^{(p)}=\pi_Y^{-1}(S)\subset Y^{(p)}$.
The relative Frobenius $F_Y=F_{Y/k}:Y\to Y^{(p)}$ satisfies $\pi_Y\circ F_Y=F_Y^{\mathrm{abs}}$,
where $F_Y^{\mathrm{abs}}$ is the absolute Frobenius of
$Y$. Since absolute Frobenius induces the identity on
the underlying topological space, we obtain
\[
 F_Y^{-1}(S^{(p)})
 =(\pi_Y\circ F_Y)^{-1}(S)
 =(F_Y^{\mathrm{abs}})^{-1}(S)
 =S
\]
as subsets of $Y$. We use the notation $T\bigl(Y/(Y-S)\bigr)=Y^{(p)}/\bigl(Y^{(p)}-S^{(p)}\bigr)$.

\begin{theorem}[Perfection invariance for pairs]\label{thm:pairs}
Let $k$ be any field, $p$ its exponential characteristic and
$K=k^{\perf}$. Let $X\in\Sm/k$ and let $Z\subset X$ be a closed subset,
for instance a smooth closed subscheme. Then
\[
 b_{X/(X-Z)}[1/p]\colon M_{\mathrm{fr},k}(X/(X-Z))[1/p]\longrightarrow
 \tau_*M_{\mathrm{fr},K}(X_K/(X_K-Z_K))[1/p]
\]
is a stable local equivalence of presheaves of $S^1$-spectra on $\Sm/k$.
Consequently, for every $V\in\Sm/k$ there is a natural zigzag of stable
equivalences
\[
 M_{\mathrm{fr},k}(X/(X-Z))[1/p]_f(V)\simeq
 M_{\mathrm{fr},K}(X_K/(X_K-Z_K))[1/p]_f(V_K).
\]
\end{theorem}

\begin{proof}
For $Z=X$ the first assertion is Theorem~\ref{thm:local-perfection}. The
second assertion follows from the first exactly as
Corollary~\ref{cor:local-perfection-fibrant} follows from
Theorem~\ref{thm:local-perfection}, since Lemma~\ref{lem:pi-fibrant}
applies to any presheaf on $\Sm/K$. We may assume that $k$ is imperfect.
We run through similar steps of the proof of
Theorem~\ref{thm:local-perfection} with a smooth scheme $Y$ replaced by
$Y/(Y-S)$ with $Y\in\Sm/k$ and $S\subset Y$ a closed subset.
Every construction there changes only framings, ambient coordinates or
target maps, or pulls back along the source. In each case the support of
the result is the zero locus of the new framing intersected with the
preimage of the relevant closed subset, and we indicate why it is finite over the source.

By~\cite[Theorem 1.2]{GNP} there are natural identifications
\[
 H_i\bigl(M_{\mathrm{fr},F}(Y/(Y-S))(U);R\bigr)\cong
 H_i\bigl(\mathcal C_F(Y/(Y-S))(U)\bigr).
\]
Framed motives of motivic Thom spaces are Segal spectra, hence connective. Therefore the first step of
the proof of Theorem~\ref{thm:local-perfection} goes through. Namely, it suffices
to show that $b^{\mathrm{lin}}_{X/(X-Z)}$ is a Nisnevich local quasi-isomorphism.

The homology presheaves
$V\mapsto H_i(\mathcal C_k(Y/(Y-S))(V))$ are $\A^1$-invariant quasi-stable additive
$\ZF_*(k)$-presheaves. Theorem~\ref{thm:main} implies that their
Nisnevich sheafifications are strictly $\A^1$-invariant. For a finite Zariski
cover $Y=\bigcup_aU_a$ put $S_I=S\cap U_I$. The inclusion $j_I:U_I\hookrightarrow Y$
defines a morphism of open pairs $(U_I,U_I-S_I)\rightarrow(Y,Y-S)$
in the sense of \cite[Definition 5.1]{GPfm}. It induces maps
\[
 (j_I)_*:
 \Fr_N(V,U_I/(U_I-S_I))
 \longrightarrow
 \Fr_N(V,Y/(Y-S)),
\]
compatible with suspension, linearization and Suslin's complex.
Over a henselian local source the support $Z(\varphi)\cap g^{-1}(S)$ splits into
local components, each of which maps into some $U_a$. The \v{C}ech
argument and the second step of Theorem~\ref{thm:local-perfection} therefore apply verbatim and give the
Nisnevich local quasi-isomorphism \eqref{eq:perfection-target-descent}
for $U_I/(U_I-S_I)$ and $Y/(Y-S)$.

Let $U\subset Y$ be an affine open admitting an \'etale map $U\to\A^d_k$, and put
\[
 S_U=S\cap U,\qquad
 Q_U=U/(U-S_U),\qquad
 Q_U^{(p)}
 =U^{(p)}/(U^{(p)}-S_U^{(p)}).
\]
We first check that the Frobenius transfer calculation
respects these pairs.

Since $F_U^{-1}(S_U^{(p)})=S_U$, relative Frobenius defines a morphism of open pairs
$(U,U-S_U)\rightarrow(U^{(p)},U^{(p)}-S_U^{(p)})$. It therefore induces
$(F_U)_*:\mathcal C_k(Q_U)\rightarrow\mathcal C_k(Q_U^{(p)})$.
For the transfer $\tr_U\in\Fr_d(U^{(p)},U)$, write $q_{\tr}$ and $h_{\tr}$ for its source projection
and target map. On its lifted support, the construction gives $q_{\tr}=F_U\circ h_{\tr}$.
Thus if a support point $z$ satisfies $h_{\tr}(z)\in S_U$, then $q_{\tr}(z)\in S_U^{(p)}$.
Lemma~\ref{lem:pair-compatible} consequently permits composition with $\tr_U$ at every finite framing level
for correspondences with target $Q_U^{(p)}$.

The correspondences $\Theta_{p,d,U}$ and $\Theta_{p,d,U^{(p)}}$
also satisfy the condition of that lemma. The same condition
holds throughout the homotopies 
\[
 F_U\circ\tr_U\simeq\Theta_{p,d,U^{(p)}},
 \qquad
 \tr_U\circ F_U\simeq\Theta_{p,d,U}
\]
from \eqref{eq:composites} (every homotopy used in these two identities satisfies 
the support condition of Lemma~\ref{lem:pair-compatible} over the entire parameter line).
Lemma~\ref{lem:pair-compatible} therefore
makes these identities available for correspondences
to the indicated quotient pairs.

Appending the negatively framed coordinate defines
$\nu$, and appending $d$ coordinates with framings
$t_1^p,\ldots,t_d^p$ defines $(\Theta_{p,d})_*$.
These operations are external products with
$\langle-1\rangle$ and $\Theta_{p,d}$, respectively.
Their homotopies are obtained by taking external
products with the homotopies of
Lemmas~\ref{lem:elem} and~\ref{lem:unit} over
$\Spec k$. The original target map is unchanged,
and the supports of these external products are
finite over the source, including throughout the
homotopies. Consequently, on stabilized homology
with $R$-coefficients,
\[
 \nu^2=1,\qquad
 (\Theta_{p,d})_*=(p\nu)^d.
\]
In particular, this scalar operation is invertible,
with inverse $p^{-d}\nu^d$.
Now the proof of the fact that $(F_U)_*$ induces an isomorphism on homology after evaluation at every $V\in\Sm/k$
is like that of Theorem~\ref{thm:local-perfection}.

Next, choose a finite Zariski cover of $Y$ by
such affine opens. Every finite intersection admits
a finite affine cover with \'etale coordinates,
obtained by restricting those on one of the original
opens. Applying the version of the
quasi-isomorphism \eqref{eq:perfection-target-descent}
for pairs, established above,
first to these intersections and then to $Y$
globalizes the preceding calculation.
Relative Frobenius is compatible with all the
inclusions, and Frobenius twist $T$ preserves the
covers and their intersections. In the end, we arrive at  the
Nisnevich local quasi-isomorphism
\[
 (F_Y)_*:
 \mathcal C_k\bigl(Y/(Y-S)\bigr)
 \xrightarrow{\simeq_{\Nis}}
 \mathcal C_k\bigl(T(Y/(Y-S))\bigr).
\]

Put $Q=Y/(Y-S)$, $TQ=Y^{(p)}/(Y^{(p)}-S^{(p)})$.
For $c=(Z,W,\varphi,g)\in\Fr_n(\Delta^q_k\times V,Q)$,
base change along $F_k$ gives a correspondence
$T(c)\in\Fr_n(\Delta^q_k\times T(V),TQ)$, where we use the canonical identification
$T(\Delta^q_k\times V)\cong\Delta^q_k\times T(V)$.
Indeed, base change preserves the neighbourhood data
and the support condition $Z=Z(\varphi)\cap g^{-1}(S)$,
and the base changed support remains finite over the base changed source.
These operations commute with the simplicial structure
maps, suspension and the additivity relations. They therefore define 
a morphism of complexes of presheaves
\[
 \beta_Q:\mathcal C_k(Q)
 \longrightarrow\mathcal C_k(TQ)\circ T,
 \qquad c\longmapsto T(c).
\]

The proof of \eqref{eq:exact} identifies
the neighbourhoods, framing functions and target maps
of its two sides. Under that identification, their
supports agree as well: taking $p$-th powers of the
framing functions does not change their zero set,
and
\[
 (F_Y\circ g)^{-1}(S^{(p)})=g^{-1}(S)
\]
because $F_Y^{-1}(S^{(p)})=S$.
Thus the common support is
$Z(\varphi)\cap g^{-1}(S)=Z$.

For the homotopy in \eqref{eq:framepower}, append
$n$ suspension coordinates and let $\lambda$ denote
the homotopy parameter. The target map is
$g\circ\pr_W$, and the common zero set of the
homotopy framing is $Z(\varphi)\times\{0\}\times\A^1_\lambda$.
Intersecting with the inverse image of $S$ gives
\[
 \bigl(Z(\varphi)\cap g^{-1}(S)\bigr)
 \times\{0\}\times\A^1_\lambda
 =
 Z\times\{0\}\times\A^1_\lambda.
\]
This support is finite over
$(\Delta^q_k\times V)\times\A^1_\lambda$.

For the homotopy in \eqref{eq:ambientpower}, both
the framing and the target map are pulled back along
$\mathcal H$. Its support is therefore the inverse
image under $\mathcal H$ of the suspended support,
extended over the parameter line.
Since $\mathcal H$ is finite, this inverse image
is finite over that support, which is itself finite
over the source times the parameter line.
Thus this construction also gives a framed homotopy
for the quotient target.

Finally, the simplex homotopy \eqref{eq:perfection-simplex-homotopy}
acts by pullback in the source, so it applies to correspondences with target $TQ$ as well.
The same argument as for ordinary targets therefore gives
\[
 H_i\bigl((F_Y)_*\bigr)
 =
 F^*\circ H_i(\beta_Q):
 H_i\bigl(\mathcal C_k(Q)\bigr)
 \longrightarrow
 H_i\bigl(\mathcal C_k(TQ)\bigr),
\]
where $F^*:H_i\bigl(\mathcal C_k(TQ)\bigr)\circ T\rightarrow H_i\bigl(\mathcal C_k(TQ)\bigr)$
is induced by pullback along $F_V:V\to T(V)$ in the source.

As we have shown above, the homology presheaves of $\mathcal C_k(TQ)$ are $\A^1$-invariant, additive
and quasi-stable. Hence $F^*$ is a Nisnevich local isomorphism by the argument of
item (b) in the proof of Theorem~\ref{thm:local-perfection}. The preceding step shows that $(F_Y)_*$ is a
Nisnevich local quasi-isomorphism. The displayed factorization implies
that $H_i(\beta_Q)$ is a Nisnevich local isomorphism for every $i$. Thus $\beta_Q$ is a local quasi-isomorphism.

Passage along $\sigma_n$ identifies $(Y_{k_n},S_{k_n})$ with $(T^nY,T^nS)$, where
$T^nS\subset T^nY$ is the closed subset lying over $S$. Under this
identification, base change along $k_n\subset k_{n+1}$ becomes
$\beta_{T^n(Y/(Y-S))}$. The fact that $b^{\mathrm{lin}}_{X/(X-Z)}$ is a
local quasi-isomorphism is proven similarly to the final step of the proof of Theorem~\ref{thm:local-perfection}.
\end{proof}

\begin{remark}\label{rem:pairs}
It is worth mentioning that smoothness of $Z$ is not used. Namely, Theorem~\ref{thm:pairs} holds for
every closed subset $Z\subset X$. For $Z=X$ it is
Theorem~\ref{thm:local-perfection}, and for
$(X\times\A^n,X\times\{0\})$ it concerns the framed motives of the
relative motivic spheres $X_+\wedge T^n$ of \cite{GNP}. 
\end{remark}

\section{Framed motives of relative motivic spheres}\label{sec:gnp}

We now deduce the stable form of the main theorem of \cite{GNP} over
every infinite field after inverting $p$, keeping the notation of Section~\ref{sec:pairs} and
following \cite{GNP}. Let $T=\A^1/\mathbb G_m$, let
$X\in\Sm/k$ and $n\geq0$. We write $M_{\mathrm{fr},k}(X\times T^n)$,
$M_{\mathrm{fr},k}(X\times T^n\times\A^1)$ and
$M_{\mathrm{fr},k}(X\times T^n\times\mathbb G_m)$ for the framed motives associated to the pairs
\[
 (X\times\A^n,\,X\times\{0\}),\qquad
 (X\times\A^n\times\A^1,\,X\times\{0\}\times\A^1),\qquad
 (X\times\A^n\times\mathbb G_m,\,X\times\{0\}\times\mathbb G_m),
\]
with $X\times T^0=X$. The linear framed presheaves
$\ZF_k(X_+\wedge T^{n+1})$, $\ZF_k(X_+\wedge T^n\wedge\A^1_+)$ and
$\ZF_k(X_+\wedge T^n\wedge\mathbb G_{m+})$ of \cite{GNP} are the
presheaves $\ZF_k(-,Y/(Y-S))$ associated to the pairs $B=(X\times\A^{n+1},
X\times\{0\})$, $A=(X\times\A^n\times\A^1,X\times\{0\}\times\A^1)$ and
$G=(X\times\A^n\times\mathbb G_m,X\times\{0\}\times\mathbb G_m)$. We write
$\mathcal C_k(A)$, $\mathcal C_k(G)$ and $\mathcal C_k(B)$ for the
associated complexes of Section~\ref{sec:pairs}. The open immersion
$\mathbb G_m\subset\A^1$ induces injective maps
$\ZF_k(-,G)\to\ZF_k(-,A)$ and levelwise injective maps
$M_{\mathrm{fr},k}(X\times T^n\times\mathbb G_m)\to
M_{\mathrm{fr},k}(X\times T^n\times\A^1)$. The canonical morphism
$\A^1_+\to T$ induces maps $\ZF_k(-,A)\to\ZF_k(-,B)$ and
$M_{\mathrm{fr},k}(X\times T^n\times\A^1)\to M_{\mathrm{fr},k}(X\times
T^{n+1})$. The composites through the
$\mathbb G_m$-terms are trivial, so we obtain maps
\begin{equation}\label{eq:gnp-cofiber}
 q_{X,n}\colon M_{\mathrm{fr},k}(X\times T^n\times\A^1)/
 M_{\mathrm{fr},k}(X\times T^n\times\mathbb G_m)\longrightarrow
 M_{\mathrm{fr},k}(X\times T^{n+1})
\end{equation}
and
\[
 \tau_{X,n}\colon\mathcal C_k(A)/\mathcal C_k(G)\longrightarrow
 \mathcal C_k(B).
\]
Finally, let $\A^1//\mathbb G_m$ be the simplicial object of $\Fr_0(k)$
which is termwise $\A^1$, $\A^1\sqcup\mathbb G_m$,
$\A^1\sqcup\mathbb G_m\sqcup\mathbb G_m,\ldots$, let
$(\A^1//\mathbb G_m)^{\wedge n}$ be its $n$-th monoidal power, and let
$M_{\mathrm{fr},k}(X\times(\A^1//\mathbb G_m)^{\wedge n})$ be the diagonal
of the framed motives of its terms. The morphism
$\alpha^{\wedge n}\colon(\A^1//\mathbb G_m)^{\wedge n}\to
(\A^1,\mathbb G_m)^{\wedge n}$ of \cite[Section~1]{GNP} induces a map
\[
 M_{\mathrm{fr},k}(\id_X\times\alpha^{\wedge n})\colon
 M_{\mathrm{fr},k}(X\times(\A^1//\mathbb G_m)^{\wedge n})\longrightarrow
 M_{\mathrm{fr},k}(X\times T^n).
\]

\begin{theorem}\label{thm:gnp}
Let $k$ be an infinite field of exponential characteristic $p$, and
let $X\in\Sm/k$.
\begin{enumerate}
\item[(1)] For every $n\geq0$ the map $\tau_{X,n}$ is a Nisnevich local
quasi-isomorphism.
\item[(2)] For every $n\geq0$ the sequence
\[
 M_{\mathrm{fr},k}(X\times T^n\times\mathbb G_m)[1/p]\longrightarrow
 M_{\mathrm{fr},k}(X\times T^n\times\A^1)[1/p]\longrightarrow
 M_{\mathrm{fr},k}(X\times T^{n+1})[1/p]
\]
is locally a homotopy cofiber sequence, that
is $q_{X,n}[1/p]$ is a stable local equivalence.
\item[(3)] For every $n\geq1$ the map
$M_{\mathrm{fr},k}(\id_X\times\alpha^{\wedge n})[1/p]$ is a stable local
equivalence.
\end{enumerate}
\end{theorem}

\begin{proof}
(1) The field $K=k^{\perf}$ is infinite and perfect, so
\cite[Theorem 1.3]{GNP} applies over $K$.
The map $\tau_{X_K,n}$ is a
Nisnevich local quasi-isomorphism on $\Sm/K$: tensor the integral
assertion of \cite[Theorem 1.3]{GNP} with the flat $\mathbb Z$-module
$R$. Base change of framed correspondences is compatible
with the maps induced by $\mathbb G_m\subset\A^1$ and by $\A^1_+\to T$, so
the morphisms $b^{\mathrm{lin}}$ of Section~\ref{sec:pairs} give a commutative square
\[
 \begin{array}{ccc}
 \mathcal C_k(A)/\mathcal C_k(G) & \xrightarrow{\ \tau_{X,n}\ } &
 \mathcal C_k(B)\\[2pt]
 \big\downarrow & & \big\downarrow\\[2pt]
 \tau_*\bigl(\mathcal C_K(A_K)/\mathcal C_K(G_K)\bigr) &
 \xrightarrow{\ \tau_*\tau_{X_K,n}\ } & \tau_*\mathcal C_K(B_K).
 \end{array}
\]
By the proof of Theorem~\ref{thm:pairs} $b^{\mathrm{lin}}$
is a Nisnevich local quasi-isomorphism for every pair. Since
$\mathcal C_k(G)\to\mathcal C_k(A)$ and its counterpart over $K$ are
injective, the long exact homology sequences and the five lemma show that
both vertical maps are Nisnevich local quasi-isomorphisms. The lower map
is one as well, because $\tau_*$ preserves stalks, as in the proof of
Lemma~\ref{lem:pi-fibrant}(i). Hence $\tau_{X,n}$ is a Nisnevich local
quasi-isomorphism.

(2) The framed motives involved are Segal
spectra, hence connective, and so is the source of $q_{X,n}$, which is a
model for the cofiber because the map from the $\mathbb G_m$-term is
levelwise injective. By the proof of Theorem~\ref{thm:pairs},
the linearization maps identify the $R$-homology presheaves of the framed
motives associated to the pairs $G$, $A$ and $B$ with the homology presheaves of
$\mathcal C_k(G)$, $\mathcal C_k(A)$ and $\mathcal C_k(B)$. By the five lemma they identify the
$R$-homology presheaves of the source of $q_{X,n}$ with those of
$\mathcal C_k(A)/\mathcal C_k(G)$. By the first statement $q_{X,n}$ therefore induces isomorphisms of the
Nisnevich sheaves associated with the $R$-homology presheaves. 
It follows that the stalks of
the cofiber of $q_{X,n}[1/p]$ are connective $R$-local spectra with zero
$R$-homology, hence contractible. Thus $q_{X,n}[1/p]$ is a stable local equivalence.

(3) Put $L=\A^1//\mathbb G_m$, $L^{\wedge0}=\mathrm{pt}$, and write
\[
 \theta_{X,n}
 =
 M_{\mathrm{fr},k}
 \bigl(\id_X\times\alpha^{\wedge n}\bigr)[1/p]:
 M_{\mathrm{fr},k}(X\times L^{\wedge n})[1/p]
 \longrightarrow
 M_{\mathrm{fr},k}(X\times T^n)[1/p].
\]
We prove by induction on $n\geq0$, simultaneously for
all $X\in\Sm/k$, that $\theta_{X,n}$ is a stable local equivalence. For $n=0$, this map is the identity.

Fix $n\geq0$. To simplify notation, set
\begin{align*}
 A_1&=
 M_{\mathrm{fr},k}
 (X\times L^{\wedge n}\times\mathbb G_m)[1/p],
 &
 B_1&=
 M_{\mathrm{fr},k}
 (X\times L^{\wedge n}\times\A^1)[1/p],\\
 A_2&=
 M_{\mathrm{fr},k}
 (X\times T^n\times\mathbb G_m)[1/p],
 &
 B_2&=
 M_{\mathrm{fr},k}
 (X\times T^n\times\A^1)[1/p].
\end{align*}
The inclusion $\mathbb G_m\hookrightarrow\A^1$
induces morphisms $f_1:A_1\to B_1$ and
$f_2:A_2\to B_2$. Define
\[
 C_1=\operatorname{hocofib}_{\mathrm{sec}}(f_1),
 \qquad
 C_2=\operatorname{hocofib}_{\mathrm{sec}}(f_2).
\]
Here $\operatorname{hocofib}_{\mathrm{sec}}$ denotes
a functorial mapping cone model for the sectionwise
homotopy cofiber of a morphism of presheaves of
$S^1$-spectra. Thus for every $V\in\Sm/k$,
\[
 C_i(V)\simeq
 \operatorname{hocofib}\bigl(A_i(V)\to B_i(V)\bigr),
 \qquad i=1,2.
\]
In particular, these cones are formed after applying
the framed motive construction.

The maps induced by $\alpha^{\wedge n}$ give a commutative square
\[
\begin{array}{ccc}
 A_1 &\xrightarrow{\ f_1\ }& B_1\\
 \big\downarrow && \big\downarrow\\
 A_2 &\xrightarrow{\ f_2\ }& B_2.
\end{array}
\]
After the canonical permutations of the target
factors, its vertical maps are
$\theta_{X\times\mathbb G_m,n}$ and
$\theta_{X\times\A^1,n}$.
The square induces a morphism $b:C_1\rightarrow C_2$.
By the induction hypothesis both vertical maps are stable local equivalences.
The induced map on homotopy cofibers is therefore
a stable local equivalence as well. 

By construction, $L^{\wedge(n+1)}$ is obtained by
taking the simplicial cone of $L^{\wedge n}\times\mathbb G_m\rightarrow L^{\wedge n}\times\A^1$
and passing to the diagonal of the resulting bisimplicial object.
Choose compatible simplicial mapping cone models
for $C_1$ and $C_2$. The inclusions of the components
of a finite disjoint union induce natural maps
\[
 \vee_j M_{\mathrm{fr},k}(Z_j)
 \longrightarrow
 M_{\mathrm{fr},k}\!\left(\sqcup_j Z_j\right).
\]
Applying these maps degreewise to the simplicial
cone and taking realization gives a natural map
\[
 a:C_1\longrightarrow
 M_{\mathrm{fr},k}
 (X\times L^{\wedge(n+1)})[1/p].
\]
The additivity theorem for framed motives~\cite{GPfm}
makes the degreewise comparisons sectionwise stable
equivalences. Taking realization therefore
shows that $a$ is a sectionwise stable equivalence.

For the second cone, put
\[
 Q=
 \left(
 \frac{
 M_{\mathrm{fr},k}(X\times T^n\times\A^1)}
 {M_{\mathrm{fr},k}(X\times T^n\times\mathbb G_m)}
 \right)[1/p].
\]
Before localization, the map in this quotient is
levelwise injective, so its ordinary quotient
represents its sectionwise homotopy cofiber.
Since localization away from $p$ preserves homotopy
cofiber sequences, the canonical comparison $c:C_2\longrightarrow Q$
is a sectionwise stable equivalence.

By statement~(2) the morphism
\[
 d=q_{X,n}[1/p]:
 Q\longrightarrow
 M_{\mathrm{fr},k}(X\times T^{n+1})[1/p]
\]
is a stable local equivalence. The constructions give a commutative square
\[
\begin{array}{ccc}
 C_1 &\xrightarrow{\ b\ }& C_2\\[4pt]
 {\scriptstyle a}\big\downarrow
 &&
 \big\downarrow{\scriptstyle d\circ c}\\[4pt]
 M_{\mathrm{fr},k}(X\times L^{\wedge(n+1)})[1/p]
 &\xrightarrow{\ \theta_{X,n+1}\ }&
 M_{\mathrm{fr},k}(X\times T^{n+1})[1/p].
\end{array}
\]
Indeed, $b$ applies $\alpha^{\wedge n}$ to the first
$n$ factors, and the composite $d\circ c$ uses the
augmentation $\alpha:L\to T$ in the last factor.
The comparison $a$ is natural with respect to these
maps. We see that $\theta_{X,n+1}\circ a=d\circ c\circ b$.

The maps $a$ and $c$ are sectionwise stable equivalences, hence stable local
equivalences. The maps $b$ and $d$ are stable local equivalences by the induction
hypothesis and statement~(2).
Thus $\theta_{X,n+1}\circ a$ is a stable local equivalence. By two-out-of-three,
so is $\theta_{X,n+1}$. This completes the induction and proves~(3).
\end{proof}

\section{Perfection invariance of \texorpdfstring{$E$}{E}-framed motives}\label{sec:E-perfection}

We extend Theorem~\ref{thm:pairs} to the $E$-framed motives of
\cite[Section~9]{GN}. Put
$P=(\mathbb P^1,\infty)$ and $T=\A^1/(\A^1-\{0\})$.
In this section $T$ denotes the motivic sphere and $(-)^{(p)}$
denotes the Frobenius twist. As before, $R=\mathbb Z[1/p]$,
$K=k^{\perf}$ and $\tau\colon\Spec K\to\Spec k$.
The extra point to verify is perfection invariance for
\[
 \mathcal L^r_F(Y/(Y-S))
 :=
 \underline{\operatorname{Hom}}\!
 \left(P_F^{\wedge r},M_{\mathrm{fr},F}(Y/(Y-S))\right),
 \qquad r\geq0.
\]
Here and below this notation denotes the Segal spectrum of the
$\Gamma$-space as in \cite{GN}
\[
 A\longmapsto C_*\underline{\operatorname{Hom}}
 (P_F^{\wedge r},\Fr_F((Y/(Y-S))\wedge A)).
\]
As usual, replacing $M_{\mathrm{fr},F}(Y/(Y-S))$ by $M_{\mathrm{fr},F}(Y/(Y-S))[1/p]$,
we define $\mathcal L^r_F(Y/(Y-S))[1/p]$.

\begin{proposition}\label{lem:P-perfection}
Let $Y\in\Sm/k$, let $S\subset Y$ be a closed subset, and let $r\geq0$.
The canonical map
\[
 \mathcal L^r_k(Y/(Y-S))[1/p]\longrightarrow
 \tau_*\mathcal L^r_K(Y_K/(Y_K-S_K))[1/p]
\]
is a stable local equivalence. Its linearization is a
Nisnevich local quasi-isomorphism.
\end{proposition}

\begin{remark}
For a field $F$, set $Q=Y/(Y-S)$ and
\[
 \Fr_{m,F}^{\langle r\rangle}(V,Q)
 =
 \operatorname{Hom}_{\mathrm{Shv}_{\Nis,*}}
 \left(V_+\wedge P_F^{\wedge(m+r)},Q\wedge T_F^m\right).
\]
By Voevodsky's Lemma~\cite{GPfm} an element is represented by
$c=(A,W,\varphi,g)$, where $A\subset\A_V^{m+r}$ is closed and
finite over $V$, $W$ is an \'{e}tale neighbourhood of $A$,
$\varphi=(\varphi_1,\ldots,\varphi_m)$, $g\colon W\to Y$, and
\begin{equation}\label{eq:E-support}
 A=Z(\varphi)\cap g^{-1}(S)
\end{equation}
as a closed subset of $W$, identified with the lifted support.
There are $m+r$ ambient coordinates but only $m$ framing functions.
The finiteness of the support follows from Voevodsky's
Lemma. It identifies the support $A$ with a closed subset
of $V\times_k(\mathbb P^1_k)^{m+r}$ disjoint from the boundary
where at least one coordinate is $\infty$.
Thus, with its reduced induced scheme structure,
$A$ is proper over $V$ and is also closed in
$\A^{m+r}_V$, hence affine over $V$.
It is therefore finite over $V$ (cf.
\cite[Corollary 3.3 and Remark 3.4]{GPfm}).

If $S=Y$ and $r>0$, then $\Fr_{m,F}^{\langle r\rangle}(V,Y)=\{*\}$.
Indeed, the support condition becomes $A=Z(\varphi)$.
In a fibre over $V$, a nonempty zero locus of $m$
functions in an \'etale neighbourhood in
$\A^{m+r}$ cannot be zero-dimensional when $r>0$.
Thus finiteness of
$A$ over $V$ forces $A=\varnothing$.
For a general closed subset $S\subset Y$, the
additional condition $g^{-1}(S)$ can cut down
$Z(\varphi)$ to a finite support. 
\end{remark}

\begin{proof}
There is nothing to prove if $k$ is perfect. Assume that $k$ is
imperfect, so that $p>1$ and $k$ is infinite.

1. Let $\ZF_{m,F}^{\langle r\rangle}(V,Q)$ be the free abelian group
modulo the usual additivity relation, and put
\[
 \ZF_F^{\langle r\rangle}(-,Q)
 =\varinjlim_m\ZF_{m,F}^{\langle r\rangle}(-,Q),\qquad
 \mathcal C_F^r(Q)
 =N\ZF_F^{\langle r\rangle}(\Delta_F^\bullet\times-,Q)\otimes R.
\]
The suspension $\Sigma$ appends one ambient coordinate $s$ and one
framing function $s$. These are the linear presheaves
$\ZF^{P^{\wedge r}}$ of \cite[Definition 4.7]{GN}.
The connected support decomposition giving
\cite[Lemma 4.8]{GN} yields, naturally in $V$,
\begin{equation}\label{eq:E-layer-linearization}
 H_i(\mathcal L_F^r(Q)(V);R)
 \cong H_i(\mathcal C_F^r(Q)(V)).
\end{equation}

2. Precomposition with ordinary framed correspondences gives
$\ZF_F^{\langle r\rangle}(-,Q)$ an additive framed structure.
Precomposition commutes with $\Sigma$. The homology
presheaves of $\mathcal C_F^r(Q)$ are $\A^1$-invariant and quasi-stable.
Theorem~\ref{thm:main} now makes their Nisnevich sheafifications
strictly $\A^1$-invariant. The target descent argument in
the proof of Theorem~\ref{thm:pairs} applies to
$\mathcal C_k^r(Q)$: over a henselian local source, the finite
support \eqref{eq:E-support} splits into local components, and each
component maps into a member of a given finite Zariski cover of $Y$.
Shrinking the \'{e}tale neighbourhood and using additivity gives
local exactness of the augmented \v{C}ech complex before applying
$C_*$. Strict invariance then turns its singular augmentation
into a Nisnevich local quasi-isomorphism, by the argument preceding
\eqref{eq:perfection-target-descent}.

3. Appending one coordinate with framing $-s$ defines an automorphism
$\nu$ of stabilized homology with $\nu^2=1$. The proof of
\eqref{eq:perfection-scalar} changes only the appended framing
coordinates, so it applies with $r$ extra ambient coordinates.
Write $a=p\nu$. Appending $d$ coordinates with their $p$-th power
framings acts by $a^d$, an automorphism over $R$.

Let $U\subset Y$ admit \'{e}tale coordinates of dimension
$d$. Composition with the level $d$ local Frobenius transfer
$\tr_U$ constructed in the proof of
Theorem~\ref{prop:Fiso} increases both the ambient
dimension and the number of framing functions by $d$:
a correspondence with $m+r$ ambient coordinates and
$m$ framing functions acquires $m+r+d$ ambient
coordinates and $m+d$ framing functions. Thus their difference remains $r$.
The proof of Lemma~\ref{lem:pair-compatible} uses only finiteness of
\eqref{eq:E-support}, and applies to these correspondences without
change. In particular the transfer, its two composite identities
\eqref{eq:composites}, and their homotopies respect the pairs
$(U,S\cap U)$ and $(U^{(p)},(S\cap U)^{(p)})$.
Those identities, with scalar $a^d$, prove that target Frobenius
is a quasi-isomorphism on these local pieces, by the finite-level
cycle and boundary argument in the proof of
Theorem~\ref{thm:pairs}. No compatibility of a positive level
transfer with the stabilization maps is needed.
Part~2 of the proof therefore gives a Nisnevich local
quasi-isomorphism
\begin{equation}\label{eq:E-target-F}
 (F_Y)_*\colon\mathcal C_k^r(Q)
 \longrightarrow\mathcal C_k^r(Q^{(p)}),
 \qquad Q^{(p)}=Y^{(p)}/(Y^{(p)}-S^{(p)}).
\end{equation}

4. Put $N=m+r$. For a correspondence $c$, let
$c^{[p]}$ mean that all $m$ framing functions are raised to the
$p$th power, and let
$\mu_{p,N}\colon\A_V^N\to\A_V^N$ raise all $N$ ambient coordinates
to the $p$th power. The two constructions in
Lemma~\ref{lem:power} give
\begin{align}
 \Sigma^m c^{[p]}&\simeq c\boxtimes\Theta_{p,m},
 \label{eq:E-frame-power}\\
 \Sigma^N\mu_{p,N}^*c&\simeq c\boxtimes\Theta_{p,N}.
 \label{eq:E-ambient-power}
\end{align}
For \eqref{eq:E-frame-power}, append $m$ coordinates and apply
$H_u$ to each pair $(\varphi_i,s_i)$; the support is
$A\times\{0\}\times\A^1_u$. For
\eqref{eq:E-ambient-power}, append $N$ coordinates and pull
$\Sigma^Nc$ back along the finite morphism applying $H_u$
to each pair of ambient coordinates $(x_j,s_j)$.
Finiteness of the original support and of this morphism gives
finiteness over the entire parameter scheme. The target map is
pulled back as well, so \eqref{eq:E-support} holds throughout.
These homotopies are natural in the source. On stabilized
homology their scalar factors are respectively $a^m$ and
$a^{m+r}$.

5. Twisting all data defines a map of complexes
\[
 \beta_Q^r\colon\mathcal C_k^r(Q)\longrightarrow
 \mathcal C_k^r(Q^{(p)})\circ(-)^{(p)}.
\]
Take $c\in\Fr_{m,k}^{\langle r\rangle}
(\Delta_k^q\times V,Q)$ and write $B=\Delta_k^q\times V$.
The cartesian relative Frobenius square for its \'{e}tale
neighbourhood, exactly as in the proof of \eqref{eq:exact}, gives
\begin{equation}\label{eq:E-exact-F}
 (F_Y\circ c)^{[p]}
 =
 \mu_{p,m+r}^*\bigl(c^{(p)}\circ F_B\bigr).
\end{equation}
This identifies neighbourhoods, framing functions and target
maps. The supports also agree, since
$F_Y^{-1}(S^{(p)})=S$ as closed subsets.
For a singular homology class $z$ represented at level $m$,
\eqref{eq:E-frame-power}--\eqref{eq:E-exact-F} imply
\[
 a^m(F_Y)_*z=a^{m+r}\gamma_*z,
 \qquad \gamma(c)=c^{(p)}\circ F_B.
\]
Every class and every required boundary occurs at a finite level,
and $a$ is invertible. Cancellation therefore gives
$(F_Y)_*z=a^r\gamma_*z$.
The simplex homotopy \eqref{eq:perfection-simplex-homotopy}
removes the Frobenius on $\Delta^\bullet_k$. Consequently, on
homology presheaves,
\begin{equation}\label{eq:E-F-factor}
 H_i((F_Y)_*)=
 a^r\,F^*\circ H_i(\beta_Q^r).
\end{equation}
Here $F^*$ is pullback along $F_V\colon V\to V^{(p)}$ on
$H_i(\mathcal C_k^r(Q^{(p)}))$, and $a^r$ acts on the target of $F^*$. 
For $r=0$, the two scalar factors cancel completely,
recovering \eqref{eq:twist-factorization}.
For $r>0$, there are $r$ more ambient coordinates
than framing functions, so cancellation leaves the
automorphism $a^r=(p\nu)^r$ on the target.
This factor is invertible but is not assumed to
be the identity.

By part~2 the latter homology presheaves are $\A^1$-invariant,
additive and quasi-stable. Paragraph (b) of the proof of
Theorem~\ref{thm:local-perfection}, using
Theorem~\ref{prop:Fiso}, makes $F^*$ a Nisnevich local
isomorphism. Part~3 makes the left side of
\eqref{eq:E-F-factor} a Nisnevich local isomorphism, and $a^r$
is an automorphism. Hence $\beta_Q^r$ is a Nisnevich local
quasi-isomorphism.

6. Let $k_j=k^{1/p^j}$. Transport along $k\cong k_j$ identifies
the transition from $k_j$ to $k_{j+1}$ with $\beta_{Q^{(p^j)}}^r$
precomposed by the $j$th Frobenius twist of the source.
Such precomposition preserves local quasi-isomorphisms, by
paragraph (a) of the proof of Theorem~\ref{thm:local-perfection}.
Supports, \'{e}tale neighbourhoods, their refinements and all
the defining relations descend through filtered unions of fields. Thus, naturally in $V$,
\[
 \mathcal C_K^r(Q_K)(V_K)
 \cong\varinjlim_j\mathcal C_{k_j}^r(Q_{k_j})(V_{k_j}).
\]
Exactness of filtered colimits proves that
$\mathcal C_k^r(Q)\to\tau_*\mathcal C_K^r(Q_K)$ is a
Nisnevich local quasi-isomorphism. Equation
\eqref{eq:E-layer-linearization} and connective stable Hurewicz
on stalks prove the asserted stable local equivalence.
\end{proof}

Recall that a motivic Thom $T$-spectrum in the sense of
\cite[Section~2]{GN} is a $T$-spectrum $E$ with each
\begin{equation}\label{eq:E-Thom-presentation}
 E_n=\varinjlim_i V_{n,i}/(V_{n,i}-D_{n,i}),
\end{equation}
where $V_{n,i}$ is smooth and $D_{n,i}$ is a smooth closed
subscheme. Write $E_K$ for its base change, including the
base changed bonding maps. For $Q=X/(X-Z)$,
the $E$-framed motive is
\[
 M_{E,k}(Q)=C_*\Fr_k^E(Q\wedge\mathbb S),\qquad
 \Fr_k^E(Q)=\varinjlim_n
 \underline{\operatorname{Hom}}(P_k^{\wedge n},Q\wedge E_n),
\]
with the transition maps defined as in \cite{GN}. We define $M_{E,k}(Q)[1/p]$
similarly to~Definition~\ref{frmotive1/p}.

\begin{theorem}[Perfection invariance for $E$-framed motives of pairs]\label{thm:E-perfection}
Let $k$ be any field, let $p$ be its exponential characteristic,
and let $K=k^{\perf}$. Let $E$ be a motivic Thom $T$-spectrum
over $k$ in the sense of~\cite{GN}. For every
$X\in\Sm/k$ and every closed subset $Z\subset X$, put
$Q=X/(X-Z)$. The canonical base change map
\begin{equation}\label{eq:E-perfection-map}
 b_{E,Q}[1/p]\colon M_{E,k}(Q)[1/p]\longrightarrow
 \tau_*M_{E_K,K}(Q_K)[1/p]
\end{equation}
is a stable local equivalence of presheaves of $S^1$-spectra
on $\Sm/k$.
The induced map of linear complexes
\[
 N\ZF_k^E(\Delta_k^\bullet\times-,Q)\otimes R
 \longrightarrow
 \tau_*\bigl(N\ZF_K^{E_K}(\Delta_K^\bullet\times-,Q_K)\otimes R\bigr)
\]
is a Nisnevich local quasi-isomorphism.
Consequently, for every $V\in\Sm/k$ there is a natural zigzag
of stable equivalences
\begin{equation}\label{eq:E-perfection-fibrant}
 M_{E,k}(Q)[1/p]_f(V)\simeq
 M_{E_K,K}(Q_K)[1/p]_f(V_K).
\end{equation}
At the base field the unreplaced spectra already satisfy
\begin{equation}\label{eq:E-perfection-field}
 M_{E,k}(Q)(\Spec k)[1/p]\xrightarrow{\simeq}
 M_{E_K,K}(Q_K)(\Spec K)[1/p].
\end{equation}
No bounding constant, symmetric group action or multiplicative
structure on $E$ is required.
\end{theorem}

\begin{proof}
For perfect $k$ the assertions are identities, so assume that $k$
is imperfect. Let $L_nE=(E_0,\ldots,E_n,E_n\wedge T,\ldots)$ be the $n$th layer of $E$. 
By \cite[Corollary 9.3 and Lemma 9.4]{GN} there are natural isomorphisms
\begin{equation*}
 M_{E,k}(Q)\cong\varinjlim_n M_{L_nE,k}(Q),\qquad
 M_{L_nE,k}(Q)\cong\mathcal L_k^n(Q\wedge E_n).
\end{equation*}

Using \eqref{eq:E-Thom-presentation}, set
\[
 Q_{n,i}
 =
 \frac{X\times V_{n,i}}
 {(X\times V_{n,i})-(Z\times D_{n,i})}.
\]
The ambient scheme $X\times V_{n,i}$ is smooth. Its indicated
closed subset need not be smooth, which is allowed in Proposition~\ref{lem:P-perfection}. One has,
\[
 Q\wedge E_n\cong\varinjlim_i Q_{n,i},\qquad
 \mathcal L_k^n(Q\wedge E_n)
 \cong\varinjlim_i\mathcal L_k^n(Q_{n,i}).
\]

Apply Proposition~\ref{lem:P-perfection} with $r=n$ to each $Q_{n,i}$.
Its canonical maps commute with all the transition maps in $i$,
with the layer maps in $n$, and with base change of the bonding
maps of $E$. Taking first the filtered colimit in $i$ and then
that in $n$ proves the assertion about
\eqref{eq:E-perfection-map}. More explicitly, on every
Nisnevich stalk, ordinary $R$-homology commutes with these
filtered colimits, and the corresponding maps are homology
isomorphisms by Proposition~\ref{lem:P-perfection}. The colimit spectra are connective,
so stable Hurewicz gives a stable equivalence on that stalk.
The resulting map is precisely $c\mapsto c_K$ in the
definition of $M_E$, rather than an abstract equivalence
between the two motives.

For the linear assertion, use the linear layer identity
\cite[Lemma 9.6]{GN}:
\[
 N\ZF_k^E(\Delta_k^\bullet\times-,Q)\otimes R
 \cong\varinjlim_n\varinjlim_i\mathcal C_k^n(Q_{n,i}).
\]
The corresponding identity over $K$ is compatible with base
change, and the assertion follows directly from the linear
part of Proposition~\ref{lem:P-perfection} and exactness of filtered
colimits. Equivalently, this is the homology computation of
\cite[Lemma 9.7]{GN} applied to the preceding argument.

Finally, Lemma~\ref{lem:pi-fibrant} and the proof of
Corollary~\ref{cor:local-perfection-fibrant} give
\eqref{eq:E-perfection-fibrant}. The Nisnevich stalk at the
point of $\Spec k$ is evaluation at $k$ itself, so
\eqref{eq:E-perfection-map} gives
\eqref{eq:E-perfection-field} without fibrant replacement.
\end{proof}

\begin{remark}\label{rem:E-perfection}
Taking $Z=X$ gives perfection invariance for $M_{E,k}(X)$.
Taking $E=S_T$ recovers Theorem~\ref{thm:pairs}.
The proof also applies to the Thom $T^2$-spectra of \cite{GN}:
one uses $r=2n$ in the layer formula and the cofinal sequence
of even ordinary framing levels. This includes the usual
$MSL$ and $MSp$ and the $T$-spectrum case includes $MGL$. 
\end{remark}

We next extend \cite[Proposition 9.12]{GN} to imperfect infinite base
fields after inverting $p$, in the same way as
Theorem~\ref{cor:cancellation} extends \cite[Theorem A]{AGP}. 
As for ordinary framed motives (see \cite[Section~9]{GN}), let
\[
 a_n^E\colon M_{E,k}(X)(n)\longrightarrow
 \underline{\mathrm{Hom}}\bigl(\mathbb G,M_{E,k}(X)(n+1)\bigr),
\]
where $M_{E,k}(X)(n):=M_{E,k}(X\times\mathbb G_m^{\wedge n})$, be the composite of the external product with $\mathbb G_m^{\wedge1}$ and
the map induced by the adjunction unit
$\mathbb G\to\mathbb G_m^{\wedge1}|_{\Sm/k}$. The map $a_n^E[1/p]$ is
defined in the same way.

\begin{theorem}[Cancellation for $E$-framed motives after inverting $p$]\label{thm:E-cancellation}
Let $k$ be an infinite field of exponential characteristic $p$, let
$K=k^{\perf}$, and let $E$ be a motivic Thom $T$-spectrum over $k$ in the
sense of~\cite{GN} such that $E_K$ is a Thom $T$-spectrum with bounding constant $d\leq 1$
over $K$. For every $X\in\Sm/k$ and $n\geq0$, the map
\[
 \bigl(a_n^E[1/p]\bigr)^\sharp\colon M_{E,k}(X)(n)[1/p]_f
 \longrightarrow
 \underline{\mathrm{Hom}}\bigl(\mathbb G,
 M_{E,k}(X)(n+1)[1/p]_f\bigr)
\]
is a schemewise stable equivalence. For infinite perfect $k$ this is
\cite[Proposition 9.12]{GN} with $p$ inverted.
\end{theorem}

\begin{proof}
The proof of Theorem~\ref{cor:cancellation} applies with the following
three changes. Put $M_n=M_{E,k}(X)(n)[1/p]$, a presheaf on $\Sm/k$, and
$N_n=M_{E_K,K}(X_K)(n)[1/p]$, a presheaf on $\Sm/K$. Base change defines
morphisms $b_{X,n}\colon M_n\to\tau_*N_n$.

First, $E$-framed motives are functorial in level zero framed
correspondences, and base change commutes with them. Hence $b_{X,n}$ is
the diagonal of the maps $b_{E,Q}[1/p]$ of \eqref{eq:E-perfection-map}
with $Q=(X\times(\mathbb G_m^{\wedge n})_q)_+$, where
$(\mathbb G_m^{\wedge n})_q$ runs through the terms of
$\mathbb G_m^{\wedge n}$. Each of these maps is a stable local equivalence
by Theorem~\ref{thm:E-perfection}, applied with $Z=X\times(\mathbb G_m^{\wedge n})_q$.
Hence so is $b_{X,n}$.

Second, external products with $\mathbb G_m^{\wedge1}$ commute with base
change, the adjunction unit $\mathbb G\to\mathbb G_m^{\wedge1}|_{\Sm/k}$ is
defined over the prime field, and $\tau_*$ commutes with
$\underline{\mathrm{Hom}}(\mathbb G_m^{\wedge1},-)$ and with
$\underline{\mathrm{Hom}}(\mathbb G,-)$. Hence the analogue of
\eqref{eq:cancellation-square} holds:
\[
 \underline{\mathrm{Hom}}(\mathbb G,b_{X,n+1})\circ a_n^{E,k}[1/p]
 =\tau_*\bigl(a_n^{E,K}[1/p]\bigr)\circ b_{X,n}.
\]

Third, $K$ is an infinite perfect field, so \cite[Proposition 9.12]{GN}
says that $(a_n^{E,K})^\sharp$ is a schemewise stable equivalence.
Inverting $p$ preserves schemewise stable equivalences and Nisnevich
excisiveness. Hence
$(a_n^{E,K}[1/p])^\sharp\colon N_{n,f}\to
\underline{\mathrm{Hom}}(\mathbb G,N_{n+1,f})$ is a schemewise stable
equivalence.

With these three inputs, the remaining part of the proof of
Theorem~\ref{cor:cancellation} applies verbatim. This part consists of
the formal criterion for $\varphi^\sharp$ and its application with
$\varphi=a_n^{E,k}[1/p]$, $B=\tau_*(N_{n,f})$, $B'=\tau_*(N_{n+1,f})$,
$u=\tau_*(\alpha_{N_n})\circ b_{X,n}$ and
$u'=\tau_*(\alpha_{N_{n+1}})\circ b_{X,n+1}$. It uses
Lemma~\ref{lem:pi-fibrant}(i),(ii) and evaluation at $Y_K$ for
$Y\in\Sm/k$. It shows that $(a_n^E[1/p])^\sharp$ is a schemewise stable
equivalence.
\end{proof}

\section{Triangulated categories of effective \texorpdfstring{$E$}{E}-framed motives}\label{sec:E-effective}

We extend results of \cite[Section~6]{Corr} to arbitrary
fields after inverting the exponential characteristic, and prove
perfection invariance of the resulting triangulated category.
Throughout this section, $E$ is a commutative symmetric Thom ring
$T$-spectrum with bounding constant $d\leq1$ and contractible
alternating-group action. These are the hypotheses used for the
symmetric monoidal module constructions in \cite[Section~6]{Corr}.
The bounding and symmetry hypotheses are retained here, although
Theorem~\ref{thm:E-perfection} itself does not require them.

As always, put $K=k^{\perf}$, $\tau\colon\Spec K\to\Spec k$, and $E_K=\tau^*E$.
All occurrences of $[1/p]$ below denote localization of motivic $T$-spectra
away from the exponential characteristic. Recall that the effective category of $E_F[1/p]$-modules
$\operatorname{Mod}^{\mathrm{eff}}_{\mathsf{SH}(F)[1/p]}\bigl(E_F[1/p]\bigr)$ is compactly generated by
$E_F[1/p]\wedge\Sigma_T^\infty X_+$, $X\in\Sm/F$.

We first record two statements that do not involve $E$. The second
one, Proposition~\ref{prop:SH-perfection}, is a special case of a
theorem of Elmanto and Khan~\cite{EK}. They proved that if $S$ is a scheme of
exponential characteristic $p$ and $f\colon S'\to S$ is a universal
homeomorphism, then $f^*\colon\mathsf{SH}(S)[1/p]\to\mathsf{SH}(S')[1/p]$
is an equivalence \cite[Corollary 2.1.5]{EK}. More generally,
$\mathsf{SH}$ is invariant under universal homeomorphisms after
inverting a set of primes containing the residue characteristics
\cite[Theorem 2.1.1]{EK}. In particular, $\mathsf{SH}(-)[1/p]$ is
invariant under passage to the perfection. All three statements of
Proposition~\ref{prop:SH-perfection} are consequences of their theorem:
(2) and (3) are the cases $f=\mathrm{Frob}_F$ and $f=\tau$, and (1)
follows from the case $f=F_U$. We give an independent proof. It uses
the explicit Frobenius transfers from the proof of
Theorem~\ref{prop:Fiso}, the $S^1$-stable motivic homotopy category and
continuity, and no result of \cite{EK}.
Proposition~\ref{prop:SH-perfection} is used in the proof of
Theorem~\ref{thm:E-effective}(3).

\begin{lemma}\label{lem:twist-descent}
Let $p>1$. For every $Y\in\Sm/K$ there are $n\geq0$ and $X\in\Sm/k$
together with an isomorphism of $K$-schemes
 $Y^{(p^n)}\cong X_K.$
\end{lemma}

\begin{proof}
Since $K=\bigcup_{n\geq0}k^{1/p^n}$, finite presentation gives
$n\geq0$, a smooth $k_n$-scheme $Y_n$, where $k_n=k^{1/p^n}$, and an
isomorphism $Y_n\times_{k_n}K\cong Y$. Write
$i_n\colon k_n\hookrightarrow K$ and $j\colon k\hookrightarrow K$ for the
inclusions, and put
\[
 \psi_n\colon k_n\xrightarrow{\sim}k,\quad a\longmapsto a^{p^n},
 \qquad
 X:=Y_n\times_{k_n,\psi_n}k\in\Sm/k.
\]
Let $\mathrm{Frob}_K$ be the absolute Frobenius of $K$. Since
$\mathrm{Frob}_K^n\circ i_n=j\circ\psi_n$, associativity of base change
gives
\[
 Y^{(p^n)}
 \cong(Y_n\times_{k_n,i_n}K)\times_{K,\mathrm{Frob}_K^n}K
 \cong Y_n\times_{k_n,\mathrm{Frob}_K^n\circ i_n}K
 =Y_n\times_{k_n,j\circ\psi_n}K
 \cong(Y_n\times_{k_n,\psi_n}k)\times_{k,j}K=X_K,
\]
as required.
\end{proof}

\begin{proposition}[{Perfection invariance of $\mathsf{SH}(-)[1/p]$}]
\label{prop:SH-perfection}
Let $k$ be a field of exponential characteristic $p$.
\begin{enumerate}
\item For every field $F$ of characteristic $p>0$ and every $U\in\Sm/F$,
relative Frobenius induces an isomorphism in $\mathsf{SH}(F)[1/p]$
\begin{equation}\label{eq:E-effective-stable-F}
 \Sigma_T^\infty(F_U)_+\colon\Sigma_T^\infty U_+[1/p]
 \xrightarrow{\simeq}\Sigma_T^\infty U^{(p)}_+[1/p].
\end{equation}
\item For every such $F$, pullback along the absolute Frobenius
$\mathrm{Frob}_F\colon\Spec F\to\Spec F$ is an autoequivalence of
$\mathsf{SH}(F)[1/p]$.
\item Canonical pullback is a symmetric monoidal triangulated
equivalence
\begin{equation}\label{eq:E-effective-SH-perfection}
 \tau^*\colon\mathsf{SH}(k)[1/p]\xrightarrow{\simeq}\mathsf{SH}(K)[1/p],
\end{equation}
and it restricts to an equivalence
$\mathsf{SH}^{\mathrm{eff}}(k)[1/p]\simeq\mathsf{SH}^{\mathrm{eff}}(K)[1/p]$.
\end{enumerate}
\end{proposition}

\begin{proof}
(1). If $p=1$, then $K=k$ and (1)--(2) are vacuous, so let $p>1$. We use
the following model structures (for a smooth $F$-scheme $X$ we write
$\rho(X)$ for the Nisnevich sheaf it represents):
\begin{enumerate}

\item[(i)] On pointed simplicial Nisnevich sheaves $\mathcal M_F$ on $\Sm/F$: the
pointed version of the Brown--Gersten model structure of
\cite[Section~4]{VoeBG}, localized at the mapping cylinders of the
projections $\rho(X\times\A^1)\to\rho(X)$, $X\in\Sm/F$. The
Brown--Gersten structure is associated with the cd-structure on
$\Sm/F$ given by elementary distinguished squares, which is complete,
regular and bounded and generates the Nisnevich topology
\cite[Section~2]{VoeU}. Its weak equivalences are the Nisnevich local
equivalences, and it is finitely generated, cellular and proper
\cite[Theorem~4.5 and Propositions~4.7, 4.8]{VoeBG}. Its generating
cofibrations are the maps
\[
 \partial\Delta^n\times\rho(X)\to\Delta^n\times\rho(X),\qquad
 \Delta^n\times\rho(A)\sqcup_{\partial\Delta^n\times\rho(A)}
 \partial\Delta^n\times\rho(X)\to\Delta^n\times\rho(X),
\]
for $X\in\Sm/F$ and open embeddings $A\subset X$, and its generating
trivial cofibrations are the corresponding maps with $\Lambda^{n,k}$ in
place of $\partial\Delta^n$ \cite[Lemmas~4.1 and 4.4]{VoeBG}. Every open
embedding $A\subset X$ occurs in an elementary distinguished square,
namely the one with $Y=X$, and the open embeddings occurring in
elementary distinguished squares are exactly these. The sheaves
$\rho(X)$ are finitely presented \cite[Lemma~2.8]{VoeBG}, so the localization (i)
of this finitely generated model category is almost finitely generated
\cite[Section~4]{Hov}. The weak equivalences of (i) are the motivic weak
equivalences, as for the Morel--Voevodsky structure \cite{MV}, which is
the localization at the same maps of the Joyal--Jardine structure, a
model structure with the same weak equivalences as the Brown--Gersten
one \cite[Section~4]{VoeBG}.

\item[(ii)] On $S^1$-spectra in (i): the stable model structure of
\cite[Definition~3.3]{Hov} with respect to $S^1\wedge-$. Its
weak equivalences are the stable motivic equivalences. We write $\mathsf{SH}^{S^1}(F)$ for its
homotopy category. Since sheafification is a stable motivic
equivalence, it agrees with the homotopy category of presheaves of
$S^1$-spectra used elsewhere in the paper.

\item[(iii)] On $(S^1,T)$-bispectra, that is, on spectra in (ii)
with respect to $-\wedge T$: the stable model structure of
\cite[Definition~3.3]{Hov}. Its homotopy category is a model of
$\mathsf{SH}(F)$, for the same reasons as for the
$(S^1,\mathbb G_m^{\wedge1})$-bispectra of Corollary~\ref{bispectrum} and
\cite[Sections~4 and~6]{Corr}, since $T\simeq S^1\wedge\mathbb G_m^{\wedge1}$
\cite{MV}.
\end{enumerate}
The sheaf $T=\rho(\A^1)/\rho(\mathbb G_m)$ is cofibrant in (i) as the
cofiber of the generating cofibration
$\rho(\mathbb G_m)_+\to\rho(\A^1)_+$. We write
$\Sigma^\infty\colon\mathsf{SH}^{S^1}(F)\to\mathsf{SH}(F)$ for the derived
functor of $F_0$ \cite[Definition~1.2]{Hov}. It sends
$\Sigma^\infty_{S^1}U_+$ to $\Sigma^\infty_TU_+$. Then
$\Sigma^\infty_{S^1}U_+\wedge T^c$, $U\in\Sm/F$ and $c\geq0$, are
finitely presented and cofibrant in (ii). In particular, they are compact
in $\mathsf{SH}^{S^1}(F)$, and inverting $p$ tensors Hom groups out of
them with $R=\mathbb Z[1/p]$.

Let $V,W\in\Sm/F$. By Voevodsky's lemma \cite[Section~3]{GPfm}
$\Fr_d(V,W)$ is the set of maps $V_+\wedge P_F^{\wedge d}\to
W_+\wedge T_F^{d}$ of pointed Nisnevich sheaves, compatibly with
composition with morphisms of schemes on either side. Let
$\sigma\colon P_F\to T_F$ be the canonical motivic equivalence. For
$c\in\Fr_d(V,W)$ define $\varrho_d(c)$ in $\mathsf{SH}^{S^1}(F)$ by
\[
 \varrho_d(c)=\Sigma^\infty_{S^1}c\circ
 \Sigma^\infty_{S^1}(\id_{V_+}\wedge\sigma^{\wedge d})^{-1}\colon
 \Sigma^\infty_{S^1}V_+\wedge T^d\longrightarrow
 \Sigma^\infty_{S^1}W_+\wedge T^d .
\]
Then $\varrho_d(g\circ c\circ f)=(g_+\wedge\id)\circ\varrho_d(c)\circ
(f_+\wedge\id)$ for morphisms $f\colon V'\to V$ and $g\colon W\to W'$.
In particular, $\A^1$-homotopic correspondences have the same
realization: if $h\in\Fr_d(V\times\A^1,W)$ and
$i_0,i_1\colon V\to V\times\A^1$ are the two sections, then
$\varrho_d(h\circ i_\nu)=\varrho_d(h)\circ(i_{\nu+}\wedge\id)$, and $i_0$
and $i_1$ agree in $\mathsf{SH}^{S^1}(F)$. Moreover, $\varrho_d$ is additive. Write $c=(Z,N,\varphi;g)$ with
$N\to\mathbb A^d_V$ the \'etale neighbourhood of $Z$, and put
$Q_Z=(\mathbb P^1)^d_V/((\mathbb P^1)^d_V-Z)$. Since $Z\subset\A^d_V$ is
closed in $(\mathbb P^1)^d_V$, the inclusion
$(\mathbb P^1)^d_V-\A^d_V\subset(\mathbb P^1)^d_V-Z$ induces a quotient map
$V_+\wedge P_F^{\wedge d}=(\mathbb P^1)^d_V/((\mathbb P^1)^d_V-\A^d_V)\to Q_Z$.
By Voevodsky's lemma the map of sheaves defined by $c$ is 
\[
 V_+\wedge P_F^{\wedge d}\longrightarrow Q_Z\cong\A^d_V/(\A^d_V-Z)
 \cong N/(N-Z)\xrightarrow{h}W_+\wedge T_F^d.
\]
The first isomorphism is Zariski excision. The second is Nisnevich excision, applied after removing from
$N$ the closed subset $\pi^{-1}(Z)-Z$, which does not change $N/(N-Z)$ by Zariski excision, and $h$ is
induced by $(g,\varphi)$.

Suppose now that $Z=Z_1\sqcup Z_2$, and let $j\neq i$. The quotient maps
$q_i\colon Q_Z\to Q_{Z_i}$ define an isomorphism
\[
 (q_1,q_2)\colon Q_Z\xrightarrow{\ \cong\ }Q_{Z_1}\vee Q_{Z_2}.
\]
Indeed, let $S$ be a local scheme with closed point $s$. Every element
of the stalk of $Q_Z$ at $S$ other than the base point is represented by
a unique morphism $f\colon S\to(\mathbb P^1)^d_V$ with $f(s)\in Z$. Then
$f(s)$ lies in exactly one $Z_i$, and $f$ maps to the base point of
$Q_{Z_j}$, because it factors through the open complement of $Z_j$.

By Zariski excision, the summand $Q_{Z_i}$ corresponds to the subsheaf
$((\mathbb P^1)^d_V-Z_j)/((\mathbb P^1)^d_V-Z)$ of $Q_Z$. The excision
isomorphisms and $(g,\varphi)$ are compatible with restriction to the
open subschemes $(\mathbb P^1)^d_V-Z_j$ and $N-Z_j$. So on this summand
$h$ restricts to the map $h_i$ defined by
$c_i=(Z_i,N-Z_j,\varphi|,g|)$. These $c_i$ are the two summands in the
additivity relation. Moreover, $q_i$ composed with the quotient map
$V_+\wedge P_F^{\wedge d}\to Q_Z$ is the quotient map for $c_i$.
After applying $\Sigma^\infty_{S^1}$ the wedge becomes a biproduct. So
the identity of $Q_{Z_1}\vee Q_{Z_2}$ is the sum of the two composites
of a projection with the corresponding inclusion. Therefore
$\varrho_d(c)=\varrho_d(c_1)+\varrho_d(c_2)$.
So $\varrho_d$ is defined on $\ZF_d(V,W)$ and takes equal values on 
elements identified by $\simeq$.

Let $\varepsilon,\varpi\colon T\to T$ be induced by $t\mapsto-t$ and by
$t\mapsto t^p$, respectively. Directly from the definitions,
\[
 \varrho_d(\Theta_{p,d,X})=\id_{X_+}\wedge\varpi^{\wedge d},\qquad
 \varrho_1(\langle-1\rangle_X)=\id_{X_+}\wedge\varepsilon,\qquad
 \varrho_2(\sigma_X\circ\sigma_X)=\id .
\]
Moreover,
$\varrho_2(\langle-1\rangle_X\circ\langle-1\rangle_X)
=\id_{X_+}\wedge\varepsilon\wedge\varepsilon$, because this composite is
the linear framing with matrix $\diag(-1,-1)$.

Applying $\varrho_2$ to \eqref{eq:units}
for $X=\Spec F$ and $a=b=-1$ gives $\varepsilon\wedge\varepsilon
=\id_{T\wedge T}$. Hence $\varepsilon\wedge\id_T$ and
$\id_T\wedge\varepsilon$ are mutually inverse isomorphisms. Applying
$\varrho_1$ to \eqref{eq:thetap} for $X=\Spec F$ gives
$\varpi=p\,\varepsilon$. Therefore $\varpi^{\wedge d}\wedge\id_T$ is an
isomorphism in $\mathsf{SH}^{S^1}(F)[1/p]$ for every $d\geq0$, and so is
$\varpi^{\wedge c}$ for every $c\geq2$.

Put $f_U=\Sigma^\infty_{S^1}(F_U)_+$. Let $U$ be affine with an \'etale
morphism $U\to\A^d_F$, as in the proof of Theorem~\ref{prop:Fiso}.
Applying $\varrho_d$ to \eqref{eq:composites} gives
\[
 (f_U\wedge\id_{T^d})\circ\varrho_d(\tr_U)
 =\id_{U^{(p)}_+}\wedge\varpi^{\wedge d},
 \qquad
 \varrho_d(\tr_U)\circ(f_U\wedge\id_{T^d})
 =\id_{U_+}\wedge\varpi^{\wedge d}.
\]
After smashing with $T$ both right hand sides become isomorphisms.
Hence $f_U\wedge\id_{T^{d+1}}$ has a left and a right inverse, so it is
an isomorphism, and so is $f_U\wedge\id_{T^c}$ for every $c\geq d+1$. A
general $U\in\Sm/F$ is a finite union of such affine opens. Since $U$ is
separated, the intersection of two of them is again of this form, and
relative Frobenius commutes with open immersions. Induction on the
number of opens, the Mayer--Vietoris triangles of Zariski covers and
the five lemma give
\begin{equation}\label{eq:SH-perf-S1}
 f_U\wedge\id_{T^c}\ \text{is an isomorphism in }
 \mathsf{SH}^{S^1}(F)[1/p]\ \text{whenever}\ c\geq\dim U+1.
\end{equation}
Applying $\Sigma^\infty$, under which $T$ becomes invertible, and
cancelling $T^c$ proves~(1).

(2). On Nisnevich sheaves on $\Sm/F$,
pullback $\mathrm{Frob}_F^*$ is the sheafified left Kan extension of
$U\mapsto U\times_{F,\mathrm{Frob}_F}F=U^{(p)}$. Its right adjoint is
restriction, $(\mathrm{Frob}_{F*}G)(U)=G(U^{(p)})$, which preserves
sheaves because base change along $\mathrm{Frob}_F$ takes elementary
distinguished squares to elementary distinguished squares. The functor
$\mathrm{Frob}_F^*$ sends $\rho(X)$ to $\rho(X^{(p)})$ and open
embeddings to open embeddings, and it commutes with products with
constant simplicial sets and with pushouts. Hence it maps the generating
cofibrations and the generating trivial cofibrations of the
Brown--Gersten structure to maps of the same kind. By adjunction,
$\mathrm{Frob}_{F*}$ preserves fibrations and trivial fibrations, which
are characterized by right lifting properties with respect to these
sets \cite[Definition~4.2 and Lemma~4.4]{VoeBG}. Thus
$\mathrm{Frob}_F^*$ is left Quillen for the Brown--Gersten structure.
It takes the projections $\rho(X\times\A^1)\to\rho(X)$ to projections of
the same kind, so by \cite[Theorem~2.2]{Hov} it is left Quillen
for (i). It is strong symmetric monoidal for the smash product, since base change
preserves finite products and the terminal object. Moreover,
$S^1$ is a constant sheaf and $T=\rho(\A^1)/\rho(\mathbb G_m)$ is obtained
from schemes defined over $\mathbb F_p$, on which $\mathrm{Frob}_F$ is the
identity. Hence $(\mathrm{Frob}_F^*,\tau)$, with $\tau$ the canonical isomorphism
$\mathrm{Frob}_F^*(S^1\wedge-)\cong S^1\wedge\mathrm{Frob}_F^*$, is a map
of pairs in the sense of \cite[Section~5]{Hov}, and by
\cite[Proposition~5.3]{Hov} its prolongation to $S^1$-spectra is left
Quillen for (ii) and commutes with the functors $F_n$. Since $\tau$ is an
isomorphism, this prolongation is computed levelwise. Applying
\cite[Proposition~5.3]{Hov} once more, now to the prolonged functor and
the canonical isomorphism
$\mathrm{Frob}_F^*(-\wedge T)\cong\mathrm{Frob}_F^*\wedge T$, shows that
its levelwise prolongation to $(S^1,T)$-bispectra is left Quillen
for (iii) and commutes with the functors $F_n$. Their derived functors are the pullback functors on
$\mathsf{SH}^{S^1}(F)$ and on $\mathsf{SH}(F)$. For a sheaf $G$ let
\[
 \gamma_G\colon G\longrightarrow\mathrm{Frob}_{F*}\mathrm{Frob}_F^*G
 \longrightarrow\mathrm{Frob}_F^*G
\]
be the unit followed by the map whose component at $U$ is
$(\mathrm{Frob}_F^*G)(F_U)$. Since $F_U\circ g=g^{(p)}\circ F_V$ for
every morphism $g\colon V\to U$, the component of $\gamma$ at
$\rho(U)$ is induced by $F_U$. Both functors preserve
colimits and smash products, and $F_{U\times V}=F_U\times F_V$. Hence
$\gamma\colon\id\Rightarrow\mathrm{Frob}_F^*$ is a natural transformation
with $\gamma_{G\wedge H}=\gamma_G\wedge\gamma_H$, it is the identity on
constant sheaves, and $\gamma_T=\varpi$. Applied in each level, 
$\gamma$ gives a natural map of $S^1$-spectra $E\to\mathrm{Frob}_F^*E$.
It induces a natural transformation of left Quillen endofunctors of (ii),
because $\gamma_{S^1}=\id$. It therefore induces a natural transformation
$\gamma\colon\id\Rightarrow\mathrm{Frob}_F^*$ on $\mathsf{SH}^{S^1}(F)$ with
\[
 \gamma_{\Sigma^\infty_{S^1}U_+\wedge T^c}
 =(f_U\wedge\id_{T^c})\circ(\id_{U_+}\wedge\varpi^{\wedge c}).
\]
It follows from above arguments that this is an isomorphism in $\mathsf{SH}^{S^1}(F)[1/p]$
if $c\geq\max(\dim U+1,2)$.

The model category (ii) is almost finitely generated in the sense of
\cite[Definition~4.1]{Hov}: the level structure on spectra in
(i) is, and (ii) is its localization at maps between finitely presented
cofibrant objects, so the argument of \cite[Section~4]{Hov} applies.
Moreover, sequential colimits in (ii) commute with finite products, and
the right adjoint of $-\wedge T$ commutes with sequential colimits
because $T$ is finitely presented. By
\cite[Theorem~4.12 and Corollary~4.13]{Hov}, for $X$ finitely presented
and cofibrant, with a finitely presented cylinder, and for arbitrary
$Y$ in (ii),
\[
 \operatorname{Hom}_{\mathsf{SH}(F)}(\Sigma^\infty X,\Sigma^\infty Y)
 =\varinjlim_m\operatorname{Hom}_{\mathsf{SH}^{S^1}(F)}
  (X\wedge T^m,Y\wedge T^m),
\]
with transition maps $-\wedge\id_T$. This identification is compatible
with $\mathrm{Frob}_F^*$, which commutes with $F_0$ and with $-\wedge T$. The objects
$T^{-a}\wedge\Sigma^\infty_TU_+[j]$, for $U\in\Sm/F$, $a\geq0$ and
$j\in\mathbb Z$, compactly generate $\mathsf{SH}(F)[1/p]$, and hence
\[
 \operatorname{Hom}_{\mathsf{SH}(F)[1/p]}
 \bigl(T^{-a}\wedge\Sigma^\infty_TU_+,T^{-b}\wedge\Sigma^\infty_TV_+[j]\bigr)
 =\varinjlim_m\operatorname{Hom}_{\mathsf{SH}^{S^1}(F)}
 \bigl(\Sigma^\infty_{S^1}U_+\wedge T^{b+m},
       \Sigma^\infty_{S^1}V_+\wedge T^{a+m}[j]\bigr)\otimes R .
\]
On the $m$th term, $\mathrm{Frob}_F^*$ is the map
$h\mapsto\mathrm{Frob}_F^*h=\gamma\circ h\circ\gamma^{-1}$, by
naturality of $\gamma$, as soon as the two components of $\gamma$
involved are invertible; this holds for large $m$. Hence
$\mathrm{Frob}_F^*$ induces bijections of these Hom groups. By~(1), it
sends $T^{-a}\wedge\Sigma^\infty_TU_+$ to
$T^{-a}\wedge\Sigma^\infty_TU^{(p)}_+\cong T^{-a}\wedge\Sigma^\infty_TU_+$,
so the images of the generators again form a set of compact generators.
Since $\mathrm{Frob}_F^*$ is the total left derived functor of a left
Quillen functor between stable model categories, it is a triangulated
functor and it preserves coproducts, being a left adjoint. Therefore
\cite[Lemma~4.8]{GJ} proves~(2).

(3). Put $k_n=k^{1/p^n}\subset K$, and let
$\tau_n\colon\Spec k_n\to\Spec k$ and
$\tau_{n+1,n}\colon\Spec k_{n+1}\to\Spec k_n$ be the canonical
morphisms. The isomorphism $k_{n+1}\to k_n$, $a\mapsto a^p$, induces an
isomorphism $\Spec k_n\to\Spec k_{n+1}$ whose composite with
$\tau_{n+1,n}$ is the absolute Frobenius of $\Spec k_n$. Hence, up to
the equivalence induced by an isomorphism of schemes,
$\tau_{n+1,n}^*$ is $\mathrm{Frob}_{k_n}^*$, and by~(2) it is fully
faithful after inverting $p$. Since $\Spec K=\varprojlim_n\Spec k_n$
with affine transition maps, continuity of stable motivic homotopy
theory \cite[Proposition C.12(4)]{Hoy} gives, for compact objects $A,B$
of $\mathsf{SH}(k)$, a bijection
\[
 \varinjlim_n\operatorname{Hom}_{\mathsf{SH}(k_n)}(\tau_n^*A,\tau_n^*B)
 \xrightarrow{\ \cong\ }
 \operatorname{Hom}_{\mathsf{SH}(K)}(\tau^*A,\tau^*B),
\]
hence also after tensoring with $R$. After tensoring with $R$ all
transition maps are bijective, so $\tau^*$ induces isomorphisms
$\operatorname{Hom}_{\mathsf{SH}(k)[1/p]}(A,B)\cong
\operatorname{Hom}_{\mathsf{SH}(K)[1/p]}(\tau^*A,\tau^*B)$. We apply this
to the compact generators $T^{-a}\wedge\Sigma^\infty_TX_+[j]$ with
$X\in\Sm/k$. For $Y\in\Sm/K$ choose $n$ and $X$ as in
Lemma~\ref{lem:twist-descent}. Applying~(1) $n$ times over $K$ gives
$T^{-a}\wedge\Sigma^\infty_TY_+\cong
\tau^*(T^{-a}\wedge\Sigma^\infty_TX_+)$ in $\mathsf{SH}(K)[1/p]$. Thus the
images of the generators form a set of compact generators, and hence
\cite[Lemma~4.8]{GJ} shows that $\tau^*$ is an equivalence. It is
symmetric monoidal, being the canonical pullback. The same argument
with $a=0$ applies to the effective subcategories, which are compactly
generated by the objects $\Sigma^\infty_TX_+[j]$.
\end{proof}

We introduced the spectrum $\mathbb S[1/p]=(S^0,\mathcal S^1,\mathcal S^2,\ldots)$
before Definition~\ref{frmotive1/p}. The $\Sigma_n$-actions on $\mathcal S^n=(\mathcal S^1)^{\wedge n}$ permuting
the factors, the canonical isomorphisms
$\mathcal S^n\wedge\mathcal S^m\cong\mathcal S^{n+m}$ and the units
$S^0=\mathcal S^0$ and $\delta_p\colon S^1\to\mathcal S^1$ make
$\mathbb S[1/p]$ a commutative symmetric ring spectrum: it is the free
commutative monoid $\operatorname{Sym}(\mathcal S^1)$ on $\mathcal S^1$ in
symmetric sequences, regarded as an algebra over
$\operatorname{Sym}(S^1)=\mathbb S$ via $\operatorname{Sym}(\delta_p)$. Its
underlying $S^1$-spectrum is the one defined above.

In this section we will often evaluate $\Gamma$-spaces at this spectrum
to get spectral categories and modules over them. Their homotopy categories will model triangulated categories with $1/p$-coefficients.
For a field extension $F/k$, set $P_F=(\mathbb P^1_F,\infty)$,
$T_F=\A^1_F/(\A^1_F-0)$, and write $\mathcal M_F$ for pointed
simplicial Nisnevich sheaves, i.e. for the category of pointed motivic spaces. Similarly to \cite[Section~5]{Corr}, define
\begin{equation}\label{eq:E-coeff-spectral}
 \mathcal O^{E_F}_{R,F}(X,Y)_n
 =\operatorname{Hom}_{\mathcal M_F}
    (X_+\wedge P_F^{\wedge n},Y_+\wedge(E_F)_n\wedge\mathcal S^n)
 =\Fr_n^{E_F}(X,Y\otimes\mathcal S^n).
\end{equation}
 We evaluate the corresponding pointed Hom-set degreewise.
We also put
\begin{equation}\label{eq:E-coeff-singular-spectral}
 \mathcal A_{R,F}(X,Y)_n
 =\mathcal O^{E_F}_{\Delta,R,F}(X,Y)_n
 =C_*\Fr_n^{E_F}(X,Y\otimes\mathcal S^n).
\end{equation}
In this way we get spectral categories $\mathcal O^{E_F}_{R,F}$ and $\mathcal A_{R,F}$.

Here are their structure maps and composition. The map
$P_F\to T_F$, followed by the unit of $E_F$, gives the framing
suspension in \cite[Definition 5.1]{Corr}. Combine it with
$\mathcal S^n\wedge S^1\xrightarrow{\id\wedge\delta_p}\mathcal S^{n+1}$ to obtain
\[
 \mathcal O^{E_F}_{R,F}(X,Y)_n\wedge S^1
 \longrightarrow\mathcal O^{E_F}_{R,F}(X,Y)_{n+1}.
\]
If $a:X_+\wedge P_F^{\wedge n}\to Y_+\wedge(E_F)_n\wedge\mathcal S^n$
and $b:Y_+\wedge P_F^{\wedge m}\to Z_+\wedge(E_F)_m\wedge\mathcal S^m$,
first apply $a$, then $b$, reorder the factors with the source
framing first, and apply
\[
 \mu^E_{n,m}\wedge\mu^{\mathbb S[1/p]}_{n,m}:
 (E_F)_n\wedge(E_F)_m\wedge \mathcal S^n\wedge\mathcal S^m
 \longrightarrow(E_F)_{n+m}\wedge\mathcal S^{n+m}.
\]
For $\delta\in\Sigma_n$ the action on $a$ is
\[
 a\longmapsto
 (\id_Y\wedge\delta_E\wedge\delta_B)\circ a\circ
 (\id_X\wedge\delta_P^{-1}).
\]
The units are induced by the unit of $E_F$ in degree zero.
These formulas verify the spectral category axioms exactly as
in \cite[Theorem 5.2]{Corr}: the additional factors use the
associative, unital and commutative multiplication of $\mathbb S[1/p]$.
For the singular category one also uses the diagonal of
$\Delta_F^\bullet$. Cartesian product of smooth schemes gives
both spectral categories a symmetric monoidal structure.

In particular, their categories of right modules have Day
products. The additional enrichment of
\cite[Theorem 5.5]{Corr} uses the coefficient ring
\[
 (\mathbf1_{E_F,R})_n
 =\operatorname{Hom}_{\mathcal M_F}
       (P_F^{\wedge n},(E_F)_n\wedge\mathcal S^n).
\]
The map $\mathbb S[1/p]\to\mathbf1_{E_F,R}$ is induced by the iterated
framing suspension and the unit of $E_F$. It makes the morphism
spectra, and hence the values of derived modules, $\mathbb S[1/p]$-modules.
Their stable homotopy groups are therefore $R$-modules.
Field extension preserves these constructions, their external
products, and the natural spectral functor
\[
 j_F:\mathcal O^{E_F}_{R,F}\longrightarrow\mathcal A_{R,F}.
\]

Equip right $\mathcal O^{E_F}_{R,F}$-modules and right
$\mathcal A_{R,F}$-modules with their stable projective model
structures and then with their Nisnevich localizations --- see~\cite{GP12,GPkm} for basic properties of the model structures. Write
\[
 \mathscr D_{E_F,R}(F)
   =\mathsf{SH}_{S^1}^{\Nis}(\mathcal O^{E_F}_{R,F}),
 \qquad
 \mathscr C_{E_F,R}(F)
   =\mathsf{SH}_{S^1}^{\Nis}(\mathcal A_{R,F}),
\]
and put
\[
 r_{F,X}=\mathcal O^{E_F}_{R,F}(-,X),\qquad
 s_{F,X}=\mathcal A_{R,F}(-,X).
\]
Define $D_{E_F}^{\mathrm{eff}}(F;R)$ to be the full subcategory
of $\mathscr D_{E_F,R}(F)$ whose Nisnevich sheaves of stable
homotopy groups are $\A^1$-invariant. Finally, let
$\mathscr D_{E_F,R}^{\mathrm{mot}}(F)$ be the further
localization at
\[
 r_{F,X\times\A^1}\longrightarrow r_{F,X}
 \qquad(X\in\Sm/F).
\]
The spectral category $\mathcal O^{E_F}_{R,F}$ is Nisnevich
excisive in the sense of~\cite{GP12} and $\mathscr D_{E_F,R}(F)$ is compactly generated similarly to
\cite[Theorem 6.11]{Corr}. It follows from Theorem~\ref{thm:main} that 
$\mathcal A_{R,F}$ is Nisnevich excisive, too.

\begin{theorem}\label{thm:E-effective}
Let $k$ be any field and let $E$ satisfy the hypotheses above.
\begin{enumerate}
\item Suslin localization $C_*$ is left adjoint to
the inclusion $D_E^{\mathrm{eff}}(k;R)\hookrightarrow\mathscr D_{E,R}(k)$.
It identifies the motivic localization $\mathscr D_{E,R}^{\mathrm{mot}}(k)$ of $\mathscr D_{E,R}(k)$
with $D_E^{\mathrm{eff}}(k;R)$, and there is a
triangulated equivalence
\begin{equation}\label{eq:E-effective-singular}
 D_E^{\mathrm{eff}}(k;R)
 \simeq\mathsf{SH}_{S^1}^{\Nis}
       (\mathcal O^E_{\Delta,R,k}).
\end{equation}

\item Canonical base change induces a triangulated equivalence
\begin{equation}\label{eq:E-effective-perfection}
 \mathbf L\tau^*\colon D_E^{\mathrm{eff}}(k;R)
 \xrightarrow{\simeq}D_{E_K}^{\mathrm{eff}}(K).
\end{equation}

\item Motivic realization induces a triangulated equivalence
\begin{equation}\label{eq:E-effective-modules}
 D_E^{\mathrm{eff}}(k;R)
 \simeq
 \operatorname{Mod}^{\mathrm{eff}}_{\mathsf{SH}(k)[1/p]}\bigl(E[1/p]\bigr).
\end{equation}
\end{enumerate}
These categories are compactly generated. The standard generators
on the framed side are $m_{E,k}(X)=C_*r_{k,X}$, $X\in\Sm/k$.
Base change sends $m_{E,k}(X)$ to $m_{E_K,K}(X_K)$, and realization
sends it to $E[1/p]\wedge\Sigma_T^\infty X_+$.
\end{theorem}

\begin{proof}
(1). The proof literally repeats that of \cite[Theorem 6.12]{Corr} if we apply Theorem~\ref{thm:main}.

(2). We first show that the canonical map
\begin{equation}\label{eq:E-spectral-perfection}
 b_Y\colon\mathcal A_{R,k}(-,Y)\longrightarrow \tau_*\mathcal A_{R,K}(-,Y_K)
\end{equation}
is a stable Nisnevich local equivalence. The point is to pass
from the $E$-motive of
Theorem~\ref{thm:E-perfection} to the spectral representable
appearing in \eqref{eq:E-spectral-perfection}.

For a field $F$, there are two commuting stabilizations in $C_*\Fr_m^{E_F}(U,Y\otimes\mathcal S^n)$:
one in the framing index $m$, and one in the simplicial sphere
index $n$. The structure map of $\mathcal A_F(U,Y)$ is their
diagonal composite. Cofinality of the diagonal in
$\mathbb N^2$, as in \cite[Lemma 6.10]{Corr}, gives
\begin{align}
 \widehat\pi_i\mathcal A_F(U,Y)
 &=\varinjlim_n\pi_{i+n}
     C_*\Fr_n^{E_F}(U,Y\otimes\mathcal S^n)\notag
 \cong\varinjlim_{m,n}\pi_{i+n}
     C_*\Fr_m^{E_F}(U,Y\otimes\mathcal S^n)\notag
 =\pi_i\bigl(M_{E_F,F}(Y)[1/p](U)\bigr).
 \label{eq:E-spectral-diagonal}
\end{align}
Here $\widehat\pi_i$ is the homotopy group of the underlying
ordinary spectrum of a symmetric spectrum. The groups in
the last line are those of the ordinary Segal spectrum
$M_{E_F,F}(Y)$. All the identifications are natural in
$U$, $Y$ and the ground field. They also hold on Nisnevich stalks.
Theorem~\ref{thm:E-perfection} now implies~\eqref{eq:E-spectral-perfection}.

Base change of correspondences gives a spectral functor
\[
 \mathcal A_k\longrightarrow\mathcal A_K,
 \qquad X\longmapsto X_K.
\]
Its derived extension of scalars is triangulated, preserves
coproducts, and takes $s_{k,X}$ to $s_{K,X_K}$. It descends
to the Nisnevich module categories, since base change
preserves the defining distinguished squares.
Equation~\eqref{eq:E-spectral-perfection} and
Lemma~\ref{lem:pi-fibrant} give, naturally in $X,Y$,
\begin{equation}\label{eq:E-effective-mapping}
 (s_{k,Y})_f(X)\simeq(s_{K,Y_K})_f(X_K).
\end{equation}
By spectral Yoneda,
\[
 \operatorname{Hom}_{\mathscr C_{E,R}(k)}(s_{k,X},s_{k,Y}[j])
 =\pi_{-j}(s_{k,Y})_f(X),\qquad j\in\mathbb Z.
\]
Thus \eqref{eq:E-effective-mapping} identifies the Hom groups
between all generators and all their shifts under the
actual base-change functor. Since the images are compact
and base change preserves coproducts, the compact-generator
argument of~\cite[Lemma~4.8]{GJ} proves full faithfulness on the entire category.
Via \eqref{eq:E-effective-singular} this is full faithfulness
of the functor in \eqref{eq:E-effective-perfection}.

For any field $F$ of characteristic $p>0$, the earlier transfer calculation gives
\begin{equation}\label{eq:E-effective-relative-F}
 m_{E_F,F}(F_U):m_{E_F,F}(U)
 \xrightarrow{\simeq}m_{E_F,F}(U^{(p)}).
\end{equation}
To see this without changing $E_F$, first take an affine $U$
with \'etale coordinates of dimension $d$. Compose at finite
framing levels with the ordinary transfer of
Theorem~\ref{prop:Fiso}, using the unit $S_T\to E_F$.
Its two homotopies \eqref{eq:composites} give the scalar
$(p\nu)^d$, with $\nu^2=1$, on stabilized homology.
This scalar is invertible over $R$. The finite level cycle
and boundary argument from the proof of
Theorem~\ref{thm:local-perfection} and the layer
argument of Theorem~\ref{thm:E-perfection} show that target
Frobenius induces a stable local equivalence of the
$E_F$-framed motives. The same conclusion is valid for
$s_{F,U}\to s_{F,U^{(p)}}$. This avoids requiring a
positive level transfer to commute strictly with the
stabilization maps. A finite Zariski cover by coordinate
neighbourhoods and descent prove the assertion for arbitrary
$U$. The coefficient spectrum $E_F$ has been fixed throughout.

Let $Y\in\Sm/K$, and choose $n\geq0$ and $X\in\Sm/k$ with
$Y^{(p^n)}\cong X_K$ as in Lemma~\ref{lem:twist-descent}.
Iterating \eqref{eq:E-effective-relative-F}, with $E_K$
fixed throughout, gives
\[
 m_{E_K,K}(Y)
 \cong m_{E_K,K}(Y^{(p^n)})
 \cong m_{E_K,K}(X_K)
 \cong\mathbf L\tau^*m_{E,k}(X).
\]
Every compact generator of the target therefore belongs
to the essential image. Together with preservation of coproducts and~\cite[Lemma~4.8]{GJ}, this proves
\eqref{eq:E-effective-perfection}. If $p=1$, then $K=k$
and the assertion is immediate.

(3). If $k$ is perfect, then $K=k$ and assertion~(3) is
\cite[Corollary 6.13]{Corr}, whose proof applies verbatim with the
coefficient spectrum $\mathbb S[1/p]$. This covers the case $p=1$. Over
the perfect field $K$ we write
\[
 \mathcal R_K\colon D_{E_K}^{\mathrm{eff}}(K)\xrightarrow{\simeq}
 \operatorname{Mod}^{\mathrm{eff}}_{\mathsf{SH}(K)[1/p]}E_K[1/p]
\]
for this equivalence. It sends $m_{E_K,K}(Y)$ to
$E_K[1/p]\wedge\Sigma_T^\infty Y_+$. Assume from now on that $k$ is
imperfect. Step~(a) below compares the module categories over $k$ and over
$K$. Step~(b) constructs $\mathcal R_k$ from $\mathcal R_K$.
Step~(c) identifies the underlying motivic spectra of $\mathcal R_k$ with
$\mathbb G_m$-stabilization over $k$.

(a) Let $F$ be a field. We model $\mathsf{SH}(F)$ by symmetric $T_F$-spectra in the model
category~(i) of the proof of Proposition~\ref{prop:SH-perfection}, that
is, in pointed simplicial Nisnevich sheaves with the $\A^1$-localized
Brown--Gersten structure. We equip them with the stable model structure
of \cite[Section~8]{Hov}. On the same category we also use Jardine's
stable model structure \cite[Theorem~4.15]{Jar00}. Jardine works with
presheaves, but his constructions and proofs apply verbatim to sheaves.
The two structures have the same weak equivalences. In both, the level
equivalences are the levelwise motivic weak equivalences. In both, a map
is a stable equivalence if and only if it induces weak equivalences of
homotopy function complexes into all $\Omega$-spectra, that is, into level
fibrant objects $W$ whose adjoint structure maps
$W_n\to\underline{\operatorname{Hom}}(T,W_{n+1})$ are weak equivalences ---
see \cite[Section~4.2]{Jar00} and \cite[Section~8]{Hov}. An
$\Omega$-spectrum for either structure is level equivalent to one for the
other, and homotopy function complexes are invariant under level
equivalences. Moreover, every (trivial) cofibration of the
Brown--Gersten model structure is one of Jardine's model structure.
It follows from \cite[Proposition~4.19]{Jar00} that Jardine's model structure satisfies the monoid axiom of
\cite{SS00}, and hence the Brown--Gersten model structure also
satisfies the monoid axiom. Furthermore, the Brown--Gersten model structure is cofibrantly
generated monoidal and all objects are small.

By \cite[Theorem~4.1]{SS00}, for every monoid $\Lambda$ the right
$\Lambda$-modules form a cofibrantly generated model category in which a
map is a weak equivalence or a fibration if and only if it is a stable
equivalence or a stable fibration of underlying symmetric $T$-spectra.
Its homotopy category is the localization of the category of
$\Lambda$-modules at the underlying stable equivalences. We write
$\operatorname{Mod}_{\mathsf{SH}(F)[1/p]}\Lambda[1/p]$ for its
localization at multiplication by $p$. The forgetful functor $U$ is
conservative and preserves coproducts, since wedges of stable equivalences
are stable equivalences. For cofibrant $A$ the free--forgetful Quillen
pair gives
\begin{equation}\label{eq:E-effective-free}
 \operatorname{Hom}_{\Lambda}(\Lambda\wedge A,N)\cong
 \operatorname{Hom}_{\mathsf{SH}(F)}(A,UN),
\end{equation}
also after inverting $p$. The symmetric spectra $F_a(X_+)$, $X\in\Sm/F$,
$a\geq0$, are cofibrant and represent $T^{-a}\wedge\Sigma_T^\infty X_+$.
Together with their shifts they compactly generate $\mathsf{SH}(F)$. By
\eqref{eq:E-effective-free}, the modules $\Lambda[1/p]\wedge F_a(X_+)[j]$
form a set of compact generators of
$\operatorname{Mod}_{\mathsf{SH}(F)[1/p]}\Lambda[1/p]$. Those with $a=0$
compactly generate the effective part
$\operatorname{Mod}^{\mathrm{eff}}_{\mathsf{SH}(F)[1/p]}\Lambda[1/p]$.

As in the proof of Proposition~\ref{prop:SH-perfection}(2), $\tau^*$ is
left Quillen for~(i). It is strong symmetric monoidal, and
$\tau^*T_k=T_K$. Hence its levelwise prolongation is a strong symmetric
monoidal left Quillen functor on symmetric $T$-spectra
\cite[Section~8]{Hov}, with right adjoint $\tau_*$,
$(\tau_*Z)_n(V)=Z_n(V_K)$. In particular, $E_K=\tau^*E$ is a commutative
monoid, and the unit $E\to\tau_*E_K$ is a map of monoids. Restriction of
scalars along it is right adjoint to $N\mapsto\tau^*N$, from $E$-modules
to $E_K$-modules. It preserves fibrations and trivial fibrations, because
these are detected by $U$. So $\tau^*$ is left Quillen on modules.
Its derived functor
\[
 \mathbf L\tau^*\colon
 \operatorname{Mod}_{\mathsf{SH}(k)[1/p]}E[1/p]\longrightarrow
 \operatorname{Mod}_{\mathsf{SH}(K)[1/p]}E_K[1/p]
\]
is triangulated, preserves coproducts, and satisfies
$\mathbf L\tau^*(E\wedge A)\cong E_K\wedge\tau^*A$ for cofibrant $A$. Its
right adjoint $\mathbf R\tau_*$ commutes with the forgetful functors.

All model categories of motivic spectra over $F$ used in this section
are built from~(i) and are connected by Quillen equivalences: the
identity functors between the Brown--Gersten and Jardine's
structures, the forgetful functors from symmetric to non-symmetric
spectra, levelwise $S^1$-suspension, and their adjoints. These functors
commute on the nose with $\tau^*$ and with $\tau_*$, and $\tau^*$ is left
Quillen in each Brown--Gersten model structure. Hence the derived functors
$\mathbf L\tau^*$ and $\mathbf R\tau_*$ commute with the derived
comparison equivalences, compatibly with units and counits. In
particular, Proposition~\ref{prop:SH-perfection}(3) holds in the present model.

The levelwise pullback also represents the derived pullback on the
non-cofibrant spectra we need. Each $E_m$ is a sequential colimit of
Thom spaces $V/(V-Z)$. So is $E_m\wedge X_+$, with Thom spaces
$(V\times X)/((V-Z)\times X)$. Thom spaces are cofibrant in~(i). Sequential colimits of motivic
spaces are homotopy colimits, because motivic weak equivalences are closed
under filtered colimits. Hence $\tau^*$ computes $\mathbf L\tau^*$ on
$E_m$ and on $E_m\wedge X_+$. The levels of $E\wedge F_a(X_+)$ are wedges
of spaces $E_m\wedge X_+$. Comparing with a cofibrant replacement, whose
levels are cofibrant, we see that
$U\mathbf L\tau^*(E\wedge F_a(X_+))=\tau^*(E\wedge F_a(X_+))$ represents
$\mathbf L\tau^*U(E\wedge F_a(X_+))$.

Let $A=F_a(X_+)$ and $B=F_b(Y_+)$ with $X,Y\in\Sm/k$. Under
\eqref{eq:E-effective-free}, the map
\[
 \operatorname{Hom}_{E}\bigl(E[1/p]\wedge A,E[1/p]\wedge B[j]\bigr)
 \longrightarrow
 \operatorname{Hom}_{E_K}\bigl(E_K[1/p]\wedge\tau^*A,
   E_K[1/p]\wedge\tau^*B[j]\bigr)
\]
induced by $\mathbf L\tau^*$ is identified with the map
\[
 \operatorname{Hom}_{\mathsf{SH}(k)[1/p]}\bigl(A,E\wedge B[j]\bigr)
 \longrightarrow
 \operatorname{Hom}_{\mathsf{SH}(K)[1/p]}\bigl(\tau^*A,
   \tau^*(E\wedge B)[j]\bigr)
\]
induced by $\mathbf L\tau^*$ on $\mathsf{SH}$. Indeed, $\tau^*$ commutes
with the free module functor and with its unit, and $\tau^*(E\wedge B)$
represents $\mathbf L\tau^*(E\wedge B)$ by the preceding paragraph. The
latter map is bijective by Proposition~\ref{prop:SH-perfection}(3).

Let $Y\in\Sm/K$, and choose $n\geq0$ and $X\in\Sm/k$ with
$Y^{(p^n)}\cong X_K$ as in Lemma~\ref{lem:twist-descent}. Applying
\eqref{eq:E-effective-stable-F} $n$ times over $K$ gives
$\Sigma_T^\infty Y_+\cong\Sigma_T^\infty(X_K)_+$ in $\mathsf{SH}(K)[1/p]$.
By \eqref{eq:E-effective-free}, this gives
$E_K[1/p]\wedge F_a(Y_+)\cong E_K[1/p]\wedge\tau^*F_a(X_+)$. Thus the
images $\mathbf L\tau^*(E[1/p]\wedge F_a(X_+))[j]$, $X\in\Sm/k$, compactly
generate the $E_K[1/p]$-modules, and those with $a=0$ the effective ones.
By \cite[Lemma~4.8]{GJ} $\mathbf L\tau^*$ is a triangulated equivalence
both on all modules and on effective modules. Its right adjoint
$\mathbf R\tau_*$ is its inverse. In particular,
\begin{equation}\label{eq:E-effective-module-perfection}
 \mathbf L\tau^*\colon
 \operatorname{Mod}^{\mathrm{eff}}_{\mathsf{SH}(k)[1/p]}E[1/p]
 \xrightarrow{\simeq}
 \operatorname{Mod}^{\mathrm{eff}}_{\mathsf{SH}(K)[1/p]}E_K[1/p].
\end{equation}

(b) Define
 $\mathcal R_k=\mathbf R\tau_*\circ\mathcal R_K\circ\mathbf L\tau^*\colon
 D_E^{\mathrm{eff}}(k)\longrightarrow
 \operatorname{Mod}^{\mathrm{eff}}_{\mathsf{SH}(k)[1/p]}E[1/p]$.
It is the composite of the equivalence \eqref{eq:E-effective-perfection}
of assertion~(2), the equivalence $\mathcal R_K$ and the inverse of
\eqref{eq:E-effective-module-perfection}. Hence it is a triangulated
equivalence and preserves coproducts. Since
$\mathbf L\tau^*m_{E,k}(X)\cong m_{E_K,K}(X_K)$,
\[
 \mathcal R_k(m_{E,k}(X))\cong
 \mathbf R\tau_*\bigl(E_K[1/p]\wedge\Sigma_T^\infty(X_K)_+\bigr)
 \cong E[1/p]\wedge\Sigma_T^\infty X_+,
\]
because $E_K[1/p]\wedge\Sigma_T^\infty(X_K)_+\cong
\mathbf L\tau^*\bigl(E[1/p]\wedge\Sigma_T^\infty X_+\bigr)$ and
$\mathbf R\tau_*$ is inverse to $\mathbf L\tau^*$. For the same reason,
\begin{equation}\label{eq:E-effective-realization-square}
 \mathbf L\tau^*\mathcal R_k\simeq\mathcal R_K\,\mathbf L\tau^*.
\end{equation}

(c) For a field $F$ let
$\Sigma_F\colon D_{E_F}^{\mathrm{eff}}(F)\longrightarrow\mathsf{SH}(F)[1/p]$
be the composite of three functors. The first is
\eqref{eq:E-effective-singular}. The second is derived
$\mathbb G_m$-stabilization of $\mathcal A_{R,F}$-modules,
\[
 \mathscr X\longmapsto\bigl(\mathscr X,\mathscr X(\mathbb G_m^{\wedge1}),
 \mathscr X(\mathbb G_m^{\wedge2}),\ldots\bigr),
 \qquad
 \mathscr X(\mathbb G_m^{\wedge n})
 =\mathscr X\wedge^{\mathbf L}\mathcal A_{R,F}(-,\mathbb G_m^{\wedge n}).
\]
The third is the canonical equivalence of the homotopy category of
symmetric $(S^1,\mathbb G_m^{\wedge1})$-bispectra with $\mathsf{SH}(F)$. This
construction uses only the symmetric monoidal structure of
$\mathcal A_{R,F}$ and requires no assumption on $F$. The functor
$\Sigma_F$ is triangulated and preserves coproducts. Also,
$\Sigma_F(m_{E_F,F}(X))$ is represented by the bispectrum
$(s_{F,X\times\mathbb G_m^{\wedge n}})_{n\geq0}$, where
$X\times\mathbb G_m^{\wedge n}$ is a simplicial object of $\Fr_0(F)$ and
$s_{F,-}$ is applied termwise. By the construction of $\mathcal R_K$ in the proof of
\cite[Corollary 6.8]{Corr}, $U\mathcal R_K\cong\Sigma_K$. That
proof uses \cite[Theorem 9.13]{GN}. It shows that the levelwise
Nisnevich local fibrant replacement of $\Sigma_K(m_{E_K,K}(Y))$ is
motivically fibrant.

Extension of scalars along $\mathcal A_{R,k}\to\mathcal A_{R,K}$,
$X\mapsto X_K$, is left adjoint to restriction $\tau_*$, where
$(\tau_*\mathscr Y)(V)=\mathscr Y(V_K)$. Restriction commutes with
$\mathbb G_m$-stabilization, because
$(V\times\mathbb G_m^{\wedge n})_K=V_K\times\mathbb G_m^{\wedge n}$. It
also commutes with the identification of bispectra with $\mathsf{SH}$,
by the compatibility of models established in Step~(a). The unit of this adjunction
therefore gives a natural map
\[
 \eta_{\mathscr X}\colon\Sigma_k(\mathscr X)\longrightarrow
 \mathbf R\tau_*\Sigma_K(\mathbf L\tau^*\mathscr X).
\]
For $\mathscr X=m_{E,k}(X)$ the $n$th level of $\eta_{\mathscr X}$ is the
base change map $b_{X\times\mathbb G_m^{\wedge n}}$ of
\eqref{eq:E-spectral-perfection}, applied termwise and realized. Hence it
is a stable local equivalence. Let $G=\Sigma_K(m_{E_K,K}(X_K))$, and let
$G\to G_f$ be its levelwise Nisnevich local fibrant replacement, which is
motivically fibrant. By Lemma~\ref{lem:pi-fibrant}(i),
$\tau_*G\to\tau_*G_f$ is a levelwise stable local equivalence. Moreover,
$\tau_*G_f$ is motivically fibrant, since $\tau_*$ preserves Nisnevich
excision, $\A^1$-invariance and the $\Omega_{\mathbb G_m}$-property (we tacitly use here
$(V\times\A^1)_K=V_K\times\A^1$ and
$(V\times\mathbb G_m)_K=V_K\times\mathbb G_m$). Hence the composite
$\Sigma_k(m_{E,k}(X))\to\tau_*G_f$ is a stable motivic fibrant
replacement, and $\eta_{m_{E,k}(X)}$ is an isomorphism. Both sides of
$\eta$ are triangulated and preserve coproducts, $\mathbf R\tau_*$ being
an equivalence after inverting $p$. Since the $m_{E,k}(X)$ are compact
generators, $\eta$ is an isomorphism. Thus,
\[
 U\mathcal R_k=U\mathbf R\tau_*\mathcal R_K\mathbf L\tau^*
 \cong\mathbf R\tau_*U\mathcal R_K\mathbf L\tau^*
 \cong\mathbf R\tau_*\Sigma_K\mathbf L\tau^*
 \cong\Sigma_k,
\]
the last isomorphism being $\eta^{-1}$.

By Step~(b) $\mathcal R_k$ is a
triangulated equivalence sending $m_{E,k}(X)$ to
$E[1/p]\wedge\Sigma_T^\infty X_+$. By Step~(c), its composite with the
forgetful functor is the motivic realization $\Sigma_k$. This proves assertion~(3).
\end{proof}

\begin{corollary}[Perfection invariance of all $E$-modules]
\label{cor:E-module-perfection}
Let $k$ be a field of exponential characteristic $p$, and let $E$ be a
monoid in motivic symmetric $T$-spectra over $k$ such that the levelwise
pullback $E_K=\tau^*E$ represents $\mathbf L\tau^*$ on the free modules
$E\wedge F_a(X_+)$, $X\in\Sm/k$, $a\geq0$. For instance, $E$ may be one of
the Thom spectra of Theorem~\ref{thm:E-effective}. Derived pullback
induces a triangulated equivalence
\begin{equation}\label{eq:E-all-module-perfection}
 \tau^*\colon
 \operatorname{Mod}_{\mathsf{SH}(k)[1/p]}E[1/p]
 \xrightarrow{\simeq}
 \operatorname{Mod}_{\mathsf{SH}(K)[1/p]}E_K[1/p],
\end{equation}
with inverse $\mathbf R\tau_*$, which restricts to an equivalence of
effective module categories.
\end{corollary}

\begin{proof}
This is proved in Step~(a) of the proof of
Theorem~\ref{thm:E-effective}(3), using the compact generators
$E[1/p]\wedge T^{-a}\wedge\Sigma_T^\infty X_+[j]$ for all modules and
$E[1/p]\wedge\Sigma_T^\infty X_+[j]$ for effective modules. That argument
uses only Proposition~\ref{prop:SH-perfection},
Lemma~\ref{lem:twist-descent}, the module model structures, the
free--forgetful adjunction and the assumption that the levelwise
pullback represents $\mathbf L\tau^*$ on free modules. For the Thom
spectra of Theorem~\ref{thm:E-effective}, this assumption is verified
there.
\end{proof}

We now extend \cite[Corollary 6.9]{Corr}, including full faithfulness of
suspension, to arbitrary fields with $\mathbb Z[1/p]$-coefficients. 
All module categories below are the homotopy categories of the model categories of 
$E$-modules constructed in Step (a) of the proof of Theorem~\ref{thm:E-effective}(3), localized by 
inverting multiplication by $p$. Thus $\operatorname{Mod}_{\mathsf{SH}(k)[1/p]}E[1/p]=\operatorname{Ho}(\operatorname{Mod}E)[1/p]$. 
We keep the notation of this section and put
\[
 \mathcal C_k=\mathscr C_{E,R}(k)=\mathsf{SH}^{\Nis}_{S^1}(\mathcal A_{R,k}),
 \qquad
 s_X=s_{k,X},
 \qquad
 P_X=E[1/p]\wedge\Sigma_T^\infty X_+ .
\]
As in \cite[Section~6]{Corr}, let $\mathbb G_m^{\wedge1}$ be the mapping
cone in $\Delta^{\mathrm{op}}\Fr_0(k)$ associated with
$1\colon\mathrm{pt}\to\mathbb G_m$. Denote by $-\boxtimes\mathbb G_m^{\wedge1}$
the endofunctor
\[
 \mathscr X\longmapsto\mathscr X(\mathbb G_m^{\wedge1})
 =\mathscr X\wedge_{\mathcal A_{R,k}}\mathcal A_{R,k}(-,\mathbb G_m^{\wedge1})
\]
of $\mathcal A_{R,k}$-modules. For $A\in\mathcal C_k$ we write $A(1)$ for
$A\boxtimes^{\mathbf L}\mathbb G_m^{\wedge1}$. Following
\cite[Section~8]{Hov}, we consider the stable model structure on
$\mathbb G_m^{\wedge1}$-symmetric spectra
$Sp^\Sigma(\mathsf{Mod}\,\mathcal A_{R,k},\mathbb G_m^{\wedge1})$, starting
with the Nisnevich local stable model structure on
$\mathcal A_{R,k}$-modules. Its homotopy category is the
$\mathbb G_m$-stabilization
$\mathcal C_k^{\mathrm{st}}=\mathsf{SH}^{\Nis}_{S^1,\mathbb G_m}(\mathcal A_{R,k})$
of $\mathcal C_k$. Since $\mathcal A_{R,k}$ is ``$\mathbb S[1/p]$-linear", no
further inversion of $p$ is needed.

\begin{theorem}
\label{cor:E-stabilized}
Let $k$ be a field of exponential characteristic $p$, and let $E$ be a
commutative symmetric Thom ring $T$-spectrum with bounding constant
$d\leq1$ and contractible alternating group action in the sense of~\cite{GN}. There is a
triangulated equivalence of compactly generated triangulated categories
\begin{equation}\label{eq:E-stable-realization}
 \overline{\mathcal R}_k\colon\mathcal C_k^{\mathrm{st}}
 \xrightarrow{\simeq}
 \operatorname{Mod}_{\mathsf{SH}(k)[1/p]}E[1/p].
\end{equation}
The suspension functor
$\Sigma^\infty_{\mathbb G_m}\colon\mathcal C_k\to\mathcal C_k^{\mathrm{st}}$
is fully faithful, and the canonical morphism
\begin{equation}\label{eq:E-stable-cancellation}
 \operatorname{Hom}(A,B[j])
 \xrightarrow{\simeq}
 \operatorname{Hom}(A(1),B(1)[j]),\qquad j\in\mathbb Z,
\end{equation}
is an isomorphism for all $A,B\in\mathcal C_k$. The essential image of the composite functor
$\overline{\mathcal R}_k\Sigma^\infty_{\mathbb G_m}$ is the full subcategory
$\operatorname{Mod}^{\mathrm{eff}}_{\mathsf{SH}(k)[1/p]}E[1/p]$ of effective $E[1/p]$-modules.
Moreover, the functor
$\mathbf L\tau^{*,\mathrm{st}}\colon\mathcal C_k^{\mathrm{st}}\to
\mathcal C_K^{\mathrm{st}}$ induced on $\mathbb G_m$-stabilizations by the
equivalence $\mathbf L\tau^*\colon\mathcal C_k\to\mathcal C_K$ of
Theorem~\ref{thm:E-effective}(2), carried along
\eqref{eq:E-effective-singular}, is an equivalence, and there is a natural isomorphism
\begin{equation}\label{eq:E-stable-base-change}
 \tau^*\,\overline{\mathcal R}_k\cong
 \overline{\mathcal R}_K\,\mathbf L\tau^{*,\mathrm{st}},
\end{equation}
where $\tau^*$ is the equivalence \eqref{eq:E-all-module-perfection}.
\end{theorem}

\begin{proof}
By Theorem~\ref{thm:E-effective}(3), carried along
\eqref{eq:E-effective-singular}, there is an equivalence, lifting derived
$\mathbb G_m$-stabilization to $E[1/p]$-modules,
\begin{equation}\label{eq:E-stable-effective-input}
 \mathcal R_k\colon\mathcal C_k
 \xrightarrow{\simeq}
 \operatorname{Mod}^{\mathrm{eff}}_{\mathsf{SH}(k)[1/p]}E[1/p],
 \qquad s_X\longmapsto P_X .
\end{equation}
The effective category on the right is a full subcategory of the
category of all $E[1/p]$-modules. Hence, once we have an equivalence
\eqref{eq:E-stable-realization} with
$\overline{\mathcal R}_k\Sigma^\infty_{\mathbb G_m}\cong\mathcal R_k$, the
essential image of $\overline{\mathcal R}_k\Sigma^\infty_{\mathbb G_m}$ is
$\operatorname{Mod}^{\mathrm{eff}}_{\mathsf{SH}(k)[1/p]}E[1/p]$.

If $k$ is perfect, then the corollary is \cite[Corollary 6.9]{Corr},
whose proof applies verbatim with the coefficient spectrum
$\mathbb S[1/p]$. The equivalence $\overline{\mathcal R}_k$ is the stable
realization constructed there, and it satisfies
$\overline{\mathcal R}_k\Sigma^\infty_{\mathbb G_m}\cong\mathcal R_k$. In
this case $K=k$, so the last assertion is trivial.
Assume from now on that $k$ is imperfect (hence $k$ is infinite). 
The spectrum $E_K$ satisfies the hypotheses of
the corollary over the perfect field $K$. So the corollary holds for
$E_K$ over $K$, and we write $\overline{\mathcal R}_K$ for the
corresponding equivalence.

By Theorem~\ref{thm:E-effective}(2) and
\eqref{eq:E-effective-singular}
$\mathbf L\tau^*\colon\mathcal C_k\to\mathcal C_K$ is an equivalence with
$\mathbf L\tau^*(s_X)=s_{X_K}$. Extension of scalars carries
$\mathcal A_{R,k}(-,\mathbb G_m^{\wedge1})$ to
$\mathcal A_{R,K}(-,\mathbb G_m^{\wedge1})$ and commutes with
$-\boxtimes\mathbb G_m^{\wedge1}$, so that
$\mathbf L\tau^*(A(1))\cong(\mathbf L\tau^*A)(1)$.
Since $k$ is infinite and $E_K$ satisfies the
assumptions of Theorem~\ref{thm:E-cancellation}, we can apply Theorem~\ref{thm:E-cancellation} 
to $E$. Let $U,X\in\Sm/k$ and $j\in\mathbb Z$. Since $s_U(1)$ is
represented by $U\times\mathbb G_m^{\wedge1}$, there
are natural isomorphisms
\begin{align*}
 \operatorname{Hom}_{\mathcal C_k}(s_U[j],s_X)&\cong
 \pi_jM_{E,k}(X)[1/p]_f(U),\\
 \operatorname{Hom}_{\mathcal C_k}(s_U(1)[j],s_X(1))&\cong
 \pi_j\underline{\mathrm{Hom}}\bigl(\mathbb G,
 M_{E,k}(X)(1)[1/p]_f\bigr)(U).
\end{align*}
Under these isomorphisms, the map induced by $-(1)$ is induced by
$(a_0^E[1/p])^\sharp$, which is a schemewise stable equivalence by
Theorem~\ref{thm:E-cancellation}. Hence
\eqref{eq:E-stable-cancellation} holds for $A=s_U[j]$ and $B=s_X$.

For fixed compact $A$, both sides of \eqref{eq:E-stable-cancellation}
are triangulated in $B$ and commute with coproducts in $B$, since $A(1)$ is
again compact. So the class of $B$ for which
\eqref{eq:E-stable-cancellation} holds is localizing and contains the
compact generators $s_X$. Hence it is all of $\mathcal C_k$. For fixed
$B$, both sides are triangulated in $A$ and take coproducts to products. So the
class of $A$ for which \eqref{eq:E-stable-cancellation} holds is
localizing and contains the $s_U$. Hence it is all of $\mathcal C_k$.
It follows that $\Sigma^\infty_{\mathbb G_m}$ is fully faithful. 

Extension of scalars along
$\mathcal A_{R,k}\to\mathcal A_{R,K}$ is a left Quillen functor of the
underlying stable model categories which commutes with
$-\boxtimes\mathbb G_m^{\wedge1}$. Its prolongation to $\mathbb G_m$-spectra
is therefore a triangulated coproduct preserving functor
\[
 \mathbf L\tau^{*,\mathrm{st}}\colon\mathcal C_k^{\mathrm{st}}
 \longrightarrow\mathcal C_K^{\mathrm{st}},
 \qquad
 \mathbf L\tau^{*,\mathrm{st}}\Sigma^\infty_{\mathbb G_m}
 \cong\Sigma^\infty_{\mathbb G_m}\mathbf L\tau^*,
\]
compatible with twists. For $X\in\Sm/k$ and $a\geq0$ put
$G_{X,a}=(\Sigma^\infty_{\mathbb G_m}s_X)(-a)$. The objects $G_{X,a}[j]$
form a set of compact generators of $\mathcal C_k^{\mathrm{st}}$, and
$\mathbf L\tau^{*,\mathrm{st}}(G_{X,a})\cong G_{X_K,a}$. For $q\geq a,b$,
the twist is an autoequivalence of $\mathcal C_k^{\mathrm{st}}$ commuting
with $\Sigma^\infty_{\mathbb G_m}$, and $\Sigma^\infty_{\mathbb G_m}$ is fully
faithful. Hence
\[
 \operatorname{Hom}(G_{X,a},G_{Y,b}[j])\cong
 \operatorname{Hom}_{\mathcal C_k}\bigl(s_X(q-a),s_Y(q-b)[j]\bigr),
\]
and similarly over $K$. Under these isomorphisms, the map induced by
$\mathbf L\tau^{*,\mathrm{st}}$ is the map induced by $\mathbf L\tau^*$, which
is bijective. Since $\mathbf L\tau^*$ is an equivalence, the objects
$s_{X_K}$, $X\in\Sm/k$, compactly generate $\mathcal C_K$. Hence the
$G_{X_K,a}[j]$ compactly generate $\mathcal C_K^{\mathrm{st}}$. By
\cite[Lemma~4.8]{GJ}, $\mathbf L\tau^{*,\mathrm{st}}$ is an equivalence.

Set $\overline{\mathcal R}_k=\mathbf R\tau_*\circ\overline{\mathcal R}_K\circ
 \mathbf L\tau^{*,\mathrm{st}}$,
where $\mathbf R\tau_*$ is the inverse of 
\eqref{eq:E-all-module-perfection} of
Corollary~\ref{cor:E-module-perfection}. As a composite of three
equivalences,
$\overline{\mathcal R}_k$ is an equivalence \eqref{eq:E-stable-realization}.
Moreover, by the compatibility 
$\mathbf L\tau^{*,\mathrm{st}}\Sigma^\infty_{\mathbb G_m}\cong\Sigma^\infty_{\mathbb G_m}\mathbf L\tau^*$, 
the perfect case $\overline{\mathcal R}_K\Sigma^\infty_{\mathbb G_m}\cong\mathcal R_K$, 
and the definition of $\mathcal R_k$ in Step (b) of the proof of Theorem 6.3(3),
\[
 \overline{\mathcal R}_k\Sigma^\infty_{\mathbb G_m}\cong
 \mathbf R\tau_*\overline{\mathcal R}_K\Sigma^\infty_{\mathbb G_m}\mathbf L\tau^*
 \cong\mathbf R\tau_*\mathcal R_K\mathbf L\tau^*=\mathcal R_k .
\]
This proves the assertion about the essential image. Finally, $\tau^*$ is
inverse to $\mathbf R\tau_*$ by Corollary~\ref{cor:E-module-perfection}, so
the definition of $\overline{\mathcal R}_k$ gives
\eqref{eq:E-stable-base-change}.
\end{proof}

\subsection*{AI Statement} Generative artificial intelligence (Claude and ChatGPT) was used as a supplementary 
tool in the composition of this article. The AI did not contribute to the foundational research, 
methodologies, or data analysis, but served solely as a general writing aid. The author 
maintains complete editorial control and accountability for the final manuscript.

\end{document}